\documentclass[a4paper,11pt]{article}
\usepackage[utf8]{inputenc}
\usepackage[T1]{fontenc}
\usepackage{graphicx}
\usepackage{amsmath,amsfonts,amssymb,amsthm}
\usepackage{mathtools}
\usepackage{mathrsfs}
\usepackage{bm}
\usepackage{xcolor}% color
\usepackage{fullpage}

\definecolor{myred}{rgb}{0.77, 0.0, 0.1}
\definecolor{crimson}{rgb}{0.86, 0.08, 0.24}
\definecolor{awesome}{rgb}{1.0, 0.13, 0.32}
\definecolor{newgreen}{rgb}{0.0,0.6,0.0}
\usepackage{hyperref}
\hypersetup{
  colorlinks, %
  linkcolor=blue, % myred, magenta or blue
  citecolor=blue, % blue, newgreen, awesome, crimson
  linktoc=page % default : section, alternatives : all, page, none
}

\renewcommand{\leq}{\leqslant}
\renewcommand{\geq}{\geqslant}

\newcommand{\di}{\mathrm{d}}

\newcommand{\eps}{\varepsilon}

\newcommand{\wt}{\widetilde}
\newcommand{\wh}{\widehat}
\newcommand{\ol}{\overline}

\DeclareFontFamily{U}{mathx}{}
\DeclareFontShape{U}{mathx}{m}{n}{<-> mathx10}{}
\DeclareSymbolFont{mathx}{U}{mathx}{m}{n}
\DeclareMathAccent{\widecheck}{0}{mathx}{"71}

\newcommand{\N}{\mathbb N}

\newcommand{\R}{\mathbb R}

\usepackage{xspace}
\newcommand{\ie}{{i.e.}\@\xspace} 
\newcommand{\eg}{e.g.\@\xspace}
\newcommand{\iid}{i.i.d.\@\xspace}

\newcommand{\ceil}[1]{\lceil #1 \rceil}
\newcommand{\floor}[1]{\lfloor #1 \rfloor}
\newcommand{\ints}[2]{\lbrace #1, \dots, #2\rbrace}%{[\![ #1, #2]\!]}
\newcommand{\iset}[1]{\ints{1}{#1}}
\DeclarePairedDelimiterX{\ooint}[2]{\lparen}{\rparen}{#1, #2}
\DeclarePairedDelimiterX{\ccint}[2]{\lbrack}{\rbrack}{#1, #2}
\DeclarePairedDelimiterX{\ocint}[2]{\lparen}{\rbrack}{#1, #2}
\DeclarePairedDelimiterX{\coint}[2]{\lbrack}{\rparen}{#1, #2}

\DeclarePairedDelimiterX{\innerp}[2]{\langle}{\rangle}{#1, #2}
\renewcommand{\ker}{\mathop{\mathrm{ker}}}
\newcommand{\im}{\mathop{\mathrm{im}}}

\newcommand{\tr}{\mathrm{Tr}}
\newcommand{\diam}{\mathop{\mathrm{diam}}}
\newcommand{\rad}{\mathop{\mathrm{rad}}}% 
\DeclareMathOperator{\vol}{Vol}

\newcommand{\conv}{\mathrm{conv}}

\DeclarePairedDelimiter{\parens}{\lparen}{\rparen}
\DeclarePairedDelimiter{\bracks}{\lbrack}{\rbrack}
\DeclarePairedDelimiter{\braces}{\lbrace}{\rbrace}
\DeclarePairedDelimiter{\set}{\lbrace}{\rbrace}
\DeclarePairedDelimiter{\abs}{\lvert}{\rvert}
\DeclarePairedDelimiter{\norm}{\lVert}{\rVert}
\DeclarePairedDelimiter{\opnorm}{\lVert}{\rVert_{\mathrm{op}}}
\DeclarePairedDelimiter{\hsnorm}{\lVert}{\rVert_{\mathrm{HS}}}

\renewcommand{\P}{\operatorname{\mathsf{P}}}
\DeclareMathOperator{\E}{\mathsf{E}}
\DeclareMathOperator{\var}{\mathsf{Var}}
\DeclareMathOperator{\cov}{\mathsf{Cov}}
\DeclarePairedDelimiterXPP{\Probab}[1]{\P}{\lparen}{\rparen}{}{#1}
\DeclarePairedDelimiterXPP{\Expect}[1]{\E}{\lbrack}{\rbrack}{}{#1}
\DeclarePairedDelimiterXPP{\Var}[1]{\var}{\lparen}{\rparen}{}{#1}
\DeclarePairedDelimiterXPP{\Cov}[1]{\cov}{\lparen}{\rparen}{}{#1}
\DeclarePairedDelimiterXPP{\indic}[1]{\bm 1}{\lparen}{\rparen}{}{#1}

\newcommand{\gaussdist}{\mathcal{N}}

\newcommand{\sgn}{\gaussdist(0, I_n)}

\newcommand{\A}{\mathcal{A}}
\newcommand{\B}{\mathcal{B}}

\newcommand{\inrad}{\mathop{\mathrm{inrad}}\nolimits}%
\newcommand{\smin}{s_{\mathrm{min}}}
\newcommand{\smax}{s_{\mathrm{max}}}
\newcommand{\jppv}{j_{\mathrm{PPV}}}
\newcommand{\jdm}{j_{\mathrm{DM}}}
\newcommand{\jkv}{j_{\mathrm{KV}}}

\newtheorem{proposition}{Proposition}[section]
\newtheorem{theorem}{Theorem}[section]
\newtheorem{lemma}{Lemma}[section]
\newtheorem{corollary}{Corollary}[section]
\newtheorem{fact}{Fact}[section]

\theoremstyle{definition}
\newtheorem{definition}{Definition}[section]
\newtheorem{example}{Example}[section]

\theoremstyle{remark}
\newtheorem{remark}{Remark}[section]

\numberwithin{equation}{section}

\title{Quantitative characterization of Gaussian small ball probabilities}
\author{Jaouad Mourtada}
\date{\today}

\begin{document}
\maketitle

\begin{abstract}
  The quantitative study of the Gaussian measure of symmetric convex sets is a central theme at the interface of probability theory and convex geometry.
  From a probabilistic perspective, it amounts to analyzing the tail behavior of suprema of Gaussian processes.
  The typical magnitude and upper tail of Gaussian processes are well-understood, thanks to Talagrand's celebrated Majorizing Measures theorem and to Gaussian concentration.
  In contrast, despite fundamental connections with metric entropy established by Kuelbs and Li, no general sharp characterization of the lower tail is available.
  In this work, we close this gap by establishing a sharp two-sided characterization of Gaussian small ball probabilities, involving both metric entropy and local widths of the underlying convex body.
  This yields an explicit chaining-type description of lower quantiles of Gaussian processes in terms of the metric structure of the convexified index set, as well as a general comparison inequality for convex sets.
  We further obtain analogous quantitative characterizations of two important geometric quantities associated with symmetric convex sets, namely intrinsic volumes (quermassintegrals) and the radii of random sections beyond the classical Dvoretzky range.
  A unifying theme of these results is the central role of convexity in the quantitative understanding of small ball probabilities.    
\end{abstract}

\tableofcontents

\section{Introduction and main results}
\label{sec:introduction}

The behavior of the Gaussian measure of symmetric convex sets is a classical topic at the intersection of Probability, Convex Geometry, and Functional Analysis; see \eg \cite{fernique1975regularite,ledoux1996isoperimetry,ledoux1991probability,bogachev1998gaussian,li2001gaussian,latala2002some,adler2007random,lifshits2012lectures,artstein2021asymptotic,talagrand2021upper} for several monographs and surveys containing treatments of various aspects of this subject.

The primary purpose of this work is to obtain precise quantitative characterizations (in a suitable approximate sense) of such probabilities in terms of geometric features of the underlying sets.
This will allow a description of these probabilities in terms of the metric structure of the sets.
We will also present similar characterizations, involving the same geometric parameters, of two other important quantities in Convex Geometry: the intrinsic volumes (quermassintegrals) of a convex body on the one hand, and the diameter of low-dimensional sections of symmetric convex sets on the other hand.
This will imply, in particular, a precise correspondence between intrinsic volumes and the lower tail of suprema of Gaussian processes.

This section is organized as follows.
In Section~\ref{sec:gauss-proc-prob}, after some definitions, we survey known results that yield a precise characterization of both typical and large values of Gaussian processes.
Section~\ref{sec:gaussian-small-ball} is devoted to Gaussian small ball probabilities.
After reviewing existing general results, we present the first main contribution of this paper, namely a two-sided characterization of general small ball probabilities.
In Section~\ref{sec:chaining-comparison}, we deduce an explicit characterization of small ball probabilities in terms of the metric structure of the convex index set and a general comparison inequality for the distribution of the oscillation of a Gaussian process.
In Section~\ref{sec:intr-volum-wills}, we present a quantitative characterization of the behavior of intrinsic volumes of a convex body, and deduce a precise connection between small ball probabilities and intrinsic volumes.
Finally, in Section~\ref{sec:low-dim-sections} we provide a characterization of the radius of low-dimensional random sections of symmetric convex sets, which we contrast with the behavior of small ball probabilities.
Section~\ref{sec:related-work} discusses additional related work, while Section~\ref{sec:notation} sets the notation.

\subsection{Gaussian processes and large ball probabilities}
\label{sec:gauss-proc-prob}

We start with some definitions on Gaussian processes, and discuss how their large deviations (in other words, the Gaussian measure of ``large'' symmetric convex sets) are well-understood.
% This will allow us to set the notation for the subsequent 

\paragraph{Gaussian processes and convex sets.}

We consider the tail behavior of suprema of Gaussian processes.
Specifically, let $(Z_t)_{t \in T}$ a centered Gaussian process indexed by a finite set $T$.
In order to deal with a non-negative quantity, it is natural to measure the magnitude of (a realization of) the process by its \emph{oscillation}, namely
\begin{equation}
  \label{eq:def-oscillation}
  \sup_{t, t' \in T} \abs*{Z_{t} - Z_{t'}}
  \, .
\end{equation}

To adopt a more geometric language, the tail probabilities of the oscillation~\eqref{eq:def-oscillation} are connected to the Gaussian measure of certain symmetric and convex sets.
Specifically, writing $T = \set{t_1, \dots, t_n}$,
the random vector $(Z_{t_i})_{1 \leq i \leq n}$ has a centered Gaussian distribution and may thus be written as $\Sigma^{1/2} X$ with $X \sim \gaussdist (0, I_n)$ for some positive semi-definite matrix $\Sigma \in \R^{n \times n}$.
Now letting $(e_{i})_{1 \leq i \leq n}$ denote the standard basis of $\R^n$ and $T = \set{\Sigma^{1/2} e_i : 1 \leq i \leq n}$, we have 
\begin{equation*}
  \sup_{t, t' \in T} \abs*{Z_{t} - Z_{t'}}
  = \sup_{1 \leq i, j \leq n} \abs*{\innerp{e_i}{\Sigma^{1/2} X} - \innerp{e_j}{\Sigma^{1/2} X}}
  = \sup_{t,t' \in T} \abs*{\innerp*{t-t'}{X}}
  \, .
\end{equation*}
We thus focus on canonical Gaussian processes $(\innerp{t}{X})_{t \in T}$, with $T$ a compact subset of $\R^n$.

Now, let $K = \conv (T-T)$ be the convex hull of $T-T$, which is a symmetric (in the sense that $K = -K$) \emph{convex body} (convex compact set) in $\R^n$.
Denote by $K^\circ = \set{x \in \R^n : \sup_{t \in K} \innerp{t}{x} \leq 1}$ the \emph{polar set} of $K$, which is a closed convex set with nonempty interior.
Also, let $\gamma_n = \gaussdist (0, I_n)$ denote the standard Gaussian measure on $\R^n$ and $X \sim \gaussdist (0, I_n)$.
Since $\sup_{t,t' \in T} \abs{\innerp{t-t'}{x}} = \sup_{t \in T-T} \innerp{t}{x} = \sup_{t \in K} \innerp{t}{x}$ for any $x \in \R^n$, we have for any $u > 0$,
% letting $X \sim \gaussdist (0, I_n)$,
\begin{equation}
  \label{eq:low-tail-gauss-polar}
  \Probab*{\sup_{t,t' \in T} \abs*{\innerp*{t-t'}{X}} \leq u}
  = \Probab*{\sup_{t \in K} \innerp{t}{X} \leq u}
  = \gamma_n \parens*{u K^\circ}
  \, .
\end{equation}
Thus, the study of the tail probabilities of the oscillation~\eqref{eq:def-oscillation} of Gaussian processes is equivalent to that of the Gaussian measure of (dilates of) symmetric closed convex subsets of $\R^n$.
We note in passing that, although both the set $K^\circ$ and the measure $\gamma_n$ depend on the ambient  space $\R^n$, the tail probability $\gamma_n (u K^\circ)$ is ``intrinsic'' to $K$, in the sense that it is preserved under isometric linear embeddings of $\R^n$ into $\R^m$ for $m \geq n$.

\paragraph{Expected value.}

A first crucial step towards studying the supremum of a Gaussian process is to understand its typical (expected) value.
The \emph{Gaussian width} (also known as \emph{Gaussian mean width} or \emph{$\ell_*$-norm} in the literature~\cite{vershynin2018high,artstein2015asymptotic}) of a compact subset $T$ of $\R^n$ is defined as
\begin{equation}
  \label{eq:def-gaussian-width}
  w (T)
  = \Expect*{\sup_{t \in T} \innerp{t}{X}}
  = \frac{1}{2} \Expect*{\sup_{t,t' \in T} \abs*{\innerp{t}{X} - \innerp{t'}{X}}}
  \, ,
\end{equation}
where the last equality holds by symmetry of the distribution of $X \sim \gaussdist (0, I_n)$.

Characterizing the magnitude of the Gaussian width of a set in terms of its structure has been a central question in theoretical probability.
Early results of Dudley~\cite{dudley1967sizes}, Sudakov~\cite{sudakov1971gaussian,sudakov1976geometric} and Fernique~\cite{fernique1975regularite}
% suggested that the Gaussian width of a set might be understood in terms of its metric structure under Euclidean distance.
established deep connections between Gaussian width and the metric structure of the underlying set.
An explicit characterization of the Gaussian width was established by Talagrand~\cite{talagrand1987regularity} (confirming a conjecture of Fernique) in his celebrated Majorizing Measures theorem.

The modern formulation of the Majorizing Measures theorem involves the notions of \emph{admissible sequences of partitions} and \emph{generic chaining functional}, defined below.
We refer the interested reader to Talagrand's book~\cite{talagrand2021upper}
for a definitive reference on this subject.

\begin{definition}[Admissible sequences, $\gamma_2$ functional]
  \label{def:adm-seq-gen-chaining}
  Let $T \subset \R^n$ be nonempty.
  An \emph{admissible sequence} of partitions of $T$ is a sequence $(\A_j)_{j \geq 0}$ of partitions of $T$ such that, for each $j \in \N$:
  \begin{itemize}
  \item $\A_{j+1}$ is a refinement of $\A_j$ (\ie, each element of $\A_{j+1}$ is contained in an element of $\A_j$);
  \item $\abs{\A_0} = 1$ and $\abs{\A_j} \leq 2^{2^j}$ for $j \geq 1$.
  \end{itemize}
  Given an admissible sequence, for $t \in T$ we denote by $A_j (t)$ the unique element of $\A_j$ that contains $t$.
  In addition, we define the \emph{$\gamma_2$ functional} (or \emph{generic chaining functional}) by\footnote{While there is a slight conflict between the $\gamma_2$ notation for the generic chaining functional and the (equally standard) $\gamma_n$ notation for the standard Gaussian measure on $\R^n$, we will only use $\gamma_2$ for the chaining functional.}
  \begin{equation}
    \label{eq:def-gamma-2}
    \gamma_2 (T)
    = \inf_{(\A_j)_{j \geq 0}} \sup_{t \in T} \sum_{j \geq 0} 2^{j/2} \diam A_j (t)
    \, ,
  \end{equation}
  where the infimum is taken over all admissible sequences of partitions of $T$, and where $\diam (A)$ denotes the diameter of $A$ with respect to the Euclidean distance.
\end{definition}

With this definition in place, the Majorizing Measures theorem (see~\cite[Theorem~2.10.1 p.~59]{talagrand2021upper}) reads as follows.

\begin{theorem}[\cite{talagrand1987regularity,talagrand2021upper}]
  \label{thm:talagrand-maj-meas}
  There exists a universal constant $c > 0$ such that, for every $n \geq 1$ and compact $T \subset \R^n$, 
  \begin{equation}
    \label{eq:maj-meas}
    c^{-1} \gamma_2 (T)
    \leq w (T)
    \leq c \gamma_2 (T)
    \, .
  \end{equation}
\end{theorem}

A remarkable aspect of Theorem~\ref{thm:talagrand-maj-meas} is that it sharply characterizes a probabilistic quantity (the expected supremum of a Gaussian process) in terms of the metric structure of the underlying index set, in an explicit and universal manner.

\paragraph{Large deviations.}

Having discussed the expected value of the supremum of a Gaussian process, we turn to the behavior of its upper tail.

Towards this end, the fundamental property of Gaussian measures is the following concentration inequality, due to Tsirel'son, Ibragimov and Sudakov~\cite{tsirelson1976gaussian}, which
is closely related to the isoperimetric property of half-spaces under Gaussian measure~\cite{borell1975brunn,sudakov1978extremal}.

\begin{theorem}[\cite{tsirelson1976gaussian}]
  \label{thm:gaussian-concentration}
  Let $f : \R^n \to \R$ be an $L$-Lipschitz function for some $L > 0$, and $X \sim \gaussdist (0, I_n)$.
  For any $u \in \R^+$, 
  \begin{equation}
    \label{eq:gaussian-concentration}
    \Probab*{f (X) - \Expect*{f (X)} \geq u}
    \leq \exp \parens*{- \frac{u^2}{2 L^2}}
    \, .
  \end{equation}
\end{theorem}

Since the map $x \mapsto \sup_{t \in K} \innerp{t}{x}$ is $\rad (K)$-Lipschitz on $\R^n$ (where $\rad (K) = \sup_{t \in K} \norm{t} = \sup_{t, t' \in T} \norm{t-t'} = \diam (T)$ is the radius of $K$ under Euclidean norm, which equals the diameter of $T$),
Theorem~\ref{thm:gaussian-concentration} (and basic properties of Gaussian variables, see Section~\ref{sec:proof-cor-upper-tail}) immediately yields a characterization of large deviations of Gaussian processes:

\begin{corollary}
  \label{cor:upper-tail-gaussian}
  For any symmetric convex body $K \subset \R^n$, $K \neq \set{0}$, and any $u \geq 2 w (K)$,
  \begin{equation}
    \label{eq:gaussian-large-balls}
    \frac{u^2}{8 \rad (K)^2}
    \leq - \log \parens*{1 - \gamma_n (u K^\circ)}
    \leq \frac{3 u^2}{\rad (K)^2}
    \, .
  \end{equation}  
  Hence, if $X \sim \gaussdist (0, I_n)$, for any compact set $T \subset \R^n$ not reduced to a single point and $u \geq 4 w (T)$, 
  \begin{equation}
    \label{eq:upper-tail-gaussian}
    \exp \parens*{- \frac{3 u^2}{\diam (T)^2}}
    \leq \Probab[\bigg]{\sup_{t,t' \in T} \abs{\innerp{t-t'}{X}} \geq u}
    \leq \exp \parens*{- \frac{u^2}{8 \diam (T)^2}}
    \, .
  \end{equation}
\end{corollary}

In particular, large values of Gaussian processes are governed by two scalar parameters: the Gaussian width and the radius of the underlying set.
In a more geometric language, the Gaussian measure of ``large'' symmetric closed convex sets $L = u K^\circ$ (that is, large enough constant-factor dilations of symmetric closed convex sets whose Gaussian measure is bounded away from $0$) is exponentially close to $1$, with an exponent governed by the radius of the largest centered Euclidean ball contained in $L$ (indeed, if $L = u K^\circ$, this radius is $u/\rad (K)$).

Additionally, the combination of Theorems~\ref{thm:talagrand-maj-meas} and~\ref{thm:gaussian-concentration} immediately yields the following characterization of upper quantiles of Gaussian processes (see Section~\ref{sec:proof-cor-upper-quantiles} for a proof).
Below, for any real random variable $Y$ and $s \in \R_+^*$, we denote by
\begin{equation}
  \label{eq:def-upper-quantile}
  Q_+ (s, Y)
  = \sup \braces*{q \in \R : \Probab*{Y \geq q} \leq \exp (-s)}  
\end{equation}
the quantile of order $1-e^{-s}$ of $Y$.

\begin{corollary}
  \label{cor:upper-quantiles-gaussian}
  Let $X \sim \gaussdist (0, I_n)$.
  There exists a universal constant $c \geq 1$ such that, for every compact set $T \subset \R^n$ and $s \geq 3$, 
  \begin{equation}
    \label{eq:upper-quantiles-gaussian}
    \frac{1}{c} \braces*{\gamma_2 (T) + \diam (T) \sqrt{s}}
    \leq Q_+ \parens*{s, \sup_{t,t' \in T} \abs*{\innerp{t-t'}{X}}}
    \leq c \braces*{\gamma_2 (T) + \diam (T) \sqrt{s}}
    \, .
  \end{equation}
\end{corollary}

Note that both $\gamma_2 (T)$ and $\diam (T)$ are both ``metric'' quantities associated to $T$, in the sense that they are increasing functions of the pairwise distances between elements of $T$.
Thus, large deviations of the Gaussian process indexed by $T$ can be characterized explicitly in terms of the metric structure of $T$.

It is also worth noting that Theorem~\ref{thm:talagrand-maj-meas} and Corollary~\ref{cor:upper-quantiles-gaussian} hold not only in the case where $T = K$ is convex, but also in the nonconvex case.
In particular, typical and large values of Gaussian processes only depend on the structure of the index set as a metric space, while convexity does not play any special role.

\subsection{Gaussian small ball probabilities}
\label{sec:gaussian-small-ball}

Since typical and large values of Gaussian processes are well-understood as discussed in the previous section, we now turn to the remaining regime, namely small values of Gaussian processes.

The study of small ball probabilities for Gaussian processes has a long history, which can be traced at least to the study~\cite{hoffmann1979lower}.
We refer to the monograph~\cite[Chapter~7]{ledoux1996isoperimetry} and to the surveys~\cite{li2001gaussian,lifshits1999asymptotic,fatalov2003constants} for a broad overview of the topic.
In addition, a rather extensive list of articles on small ball probabilities can be found in~\cite{lifshits2023biblio}.

A major progress in the study of Gaussian small ball probabilities was achieved in the seminal work of Kuelbs and Li~\cite{kuelbs1993metric}, who put forward a general connection with covering numbers.
Specifically, the key result of Kuelbs and Li (which can be found in this form in~\cite[Theorem 7.6 p.~256]{ledoux1996isoperimetry}) reads as follows.
Below, if $A, B \subset \R^n$ we let $N (A, B)$ denote the covering number of $A$ by $B$, namely the smallest number of translates of $B$ that cover $A$ (see Definition~\ref{def:covering-numbers}).
We also let $B_2 = B_2^n$ denote the closed unit Euclidean ball in $\R^n$.

\begin{theorem}[\cite{kuelbs1993metric,ledoux1996isoperimetry}]
  \label{thm:kuelbs-li}
  For any closed symmetric convex set $L \subset \R^n$ with nonempty interior, we have
  \begin{equation}
    \label{eq:kuelbs-li}
    \frac{1}{2 \gamma_n (2 L)}
    \leq N \parens*{B_2^n, \frac{1}{\sqrt{2 \log \parens*{1/\gamma_n (L)}}} L}
    \leq \frac{1}{2 \gamma_n (L/2)^2}
    \, .
  \end{equation}
\end{theorem}

A conclusion from Theorem~\ref{thm:kuelbs-li} is that, if $- \log \gamma_n (2 L) \asymp -\log \gamma_n (L/2) \gg 1$, then $s = - \log \gamma_n (L)$ satisfies $\log N (B_2^n, \frac{1}{\sqrt{2 s}} L) \asymp s$.
Writing $L = u K^\circ$ as in the previous section, and using the duality of covering numbers~\cite{artstein2004duality} (recalled as Theorem~\ref{thm:duality-covering} below), we deduce that if additionally $\log N (K, r B_2) \asymp \log N (K, \frac{u}{\sqrt{s}} B_2^n)$ whenever $r \asymp u/ \sqrt{s}$, then $s = - \log \gamma_n (u K^\circ)$ satisfies:
\begin{equation}
  \label{eq:implicit-sb-kuelbs-li}
  \log N \parens*{K, \frac{u}{\sqrt{s}} B_2^n}
  \asymp s
  \, .
\end{equation}
Thus, whenever both the Gaussian small ball probabilities and the covering numbers behave regularly, Theorem~\ref{thm:kuelbs-li} provides a condition~\eqref{eq:implicit-sb-kuelbs-li} satisfied by the small ball probability involving covering numbers of $K$.
Under additional conditions, all (large enough) numbers $s$ satisfying the condition~\eqref{eq:implicit-sb-kuelbs-li} are of the same order of magnitude, allowing one to estimate the order of magnitude of $- \log \gamma_n (u K^\circ)$.
Based on these observations, Kuelbs and Li~\cite{kuelbs1993metric} determined the behavior for small $u$ of small ball probabilities under a variety of situations, involving certain regularly varying behaviors for small ball probabilities and covering numbers.
These results were extended to weaker regularity assumptions by Li and Linde~\cite{li1999approximation}.

While useful in several situations~\cite{kuelbs1993metric,li1999approximation}, the above conditions do not characterize the behavior of Gaussian small ball probabilities in the general case, for two reasons.
First, they require a regularly varying behavior of covering numbers and small ball probabilities, which fails to hold for general symmetric convex sets $K$.
When this assumption is not met, the condition~\eqref{eq:implicit-sb-kuelbs-li} may fail to hold for $s = - \log \gamma_n (u K^\circ)$, in which case we are left with no information on this quantity.
Second, even when the condition~\eqref{eq:implicit-sb-kuelbs-li} does hold for
$s = - \log \gamma_n (u K^\circ)$, it only constitutes an implicit condition satisfied by $s$, which may not suffice to characterize its order of magnitude.
Indeed, there may be various values of $s > 1$, of highly different orders of magnitude, that satisfy condition~\eqref{eq:implicit-sb-kuelbs-li}.
We provide examples for both cases in Sections~\ref{sec:lack-solution-kuelbs}--\ref{sec:many-solut-kuelbs}.

These limitations are not merely technical, and point towards a more fundamental issue: covering numbers alone do not suffice to characterize the magnitude of the supremum of a Gaussian process.
To start, even the expected value depends, by Talagrand's Majorizing Measures theorem (Theorem~\ref{thm:talagrand-maj-meas}), on the more sophisticated generic chaining functional (Definition~\ref{def:adm-seq-gen-chaining}) that captures spatial inhomogeneity of the set at various scales.

Another important general inequality relating small ball probabilities and covering numbers was obtained by Talagrand~\cite{talagrand1993new} in his study of the Gaussian measure of enlarged convex sets.
This inequality, in the form extracted by Ledoux~\cite[p.~257]{ledoux1996isoperimetry}, can be formulated as follows (we ignore well-understood measurability subtleties here to simplify the statement).

\begin{theorem}[\cite{talagrand1993new,ledoux1996isoperimetry}]
  \label{thm:talagrand-sb-lower-cover}
  Let $T \subset \ell^2$ be a compact subset of the infinite-dimensional sequence space $\ell^2$ \textup(with its usual scalar product $\innerp{\cdot}{\cdot}$ and unit ball $B_2$\textup), and $(X_t)_{t \in T}$ a centered Gaussian process with $\Expect*{X_{t} X_{t'}} = \innerp{t}{t'}$ for any $t, t' \in T$.
  Assume that there is a non-negative function $\psi : \R_+^* \to \R_+^*$ such that for all $r > 0$, one has $N (T, r B_2) \leq \psi (r)$ and 
  $\kappa_1 \psi (r) \leq \psi (r/2) \leq \kappa_2 \psi (r)$ for some constants $1 < \kappa_1 \leq \kappa_2 < + \infty$.
  Then, there exists a constant $\kappa_3 \geq 1$ depending only on $\kappa_1, \kappa_2$ such that for any $u > 0$,
  \begin{equation}
    \label{eq:talagrand-sb-lower-cover}
    \Probab*{\sup_{t, t' \in T} \abs*{X_t - X_{t'}} \leq u}
    \geq \exp \parens[\big]{- \kappa_3 \psi \parens*{u}}
    \, .
  \end{equation}
\end{theorem}

Theorem~\ref{thm:talagrand-sb-lower-cover} is largely complementary to Theorem~\ref{thm:kuelbs-li}.
First, while Theorem~\ref{thm:kuelbs-li} contains an implicit (recursive) condition satisfied by the small ball probability, Theorem~\ref{thm:talagrand-sb-lower-cover} provides an explicit lower bound.
Second, the condition~\eqref{eq:implicit-sb-kuelbs-li} derived from Theorem~\ref{thm:kuelbs-li} requires convexity of $K$, while Theorem~\ref{thm:talagrand-sb-lower-cover} critically holds even for nonconvex $T$.
On the other hand, while the condition of Theorem~\ref{thm:kuelbs-li} is two-sided, Theorem~\ref{thm:talagrand-sb-lower-cover} only provides a one-sided lower bound.

While simple and general, the one-sided lower bound~\eqref{eq:talagrand-sb-lower-cover} on the small ball probability turns out to be loose in many concrete situations.
For instance, if $K = T$ is a fixed $k$-dimensional symmetric convex subset of $\R^n \subset \ell^2$, one can show that $N (K, r)$ diverges at the rate $r^{-k}$ as $r \to 0^+$, and (either directly or from Theorem~\ref{thm:kuelbs-li}) the small ball probability $\gamma_n (u K^\circ)$ vanishes as $u^{c k}$ as $u \to 0$; while Theorem~\ref{thm:talagrand-sb-lower-cover} yields a lower bound in $\exp (- c \kappa_3 u^{-k})$.
The looseness of Theorem~\ref{thm:talagrand-sb-lower-cover} in certain situations comes from the fact that it is expressed in terms of covering numbers alone (much like the Kuelbs-Li condition~\eqref{eq:kuelbs-li}) and that (unlike the Kuelbs-Li condition) it holds without convexity assumptions on $T$.
It seems unlikely that Theorem~\ref{thm:talagrand-sb-lower-cover} could be significantly tightened without introducing additional parameters or structural assumptions on the set, as it applies equally to highly non-convex sets with a rich convex hull.
(We will come back to the important question of the role of convexity in Section~\ref{sec:chaining-comparison}.)

It transpires from this discussion of existing general results that, while informative results hold for general sets, the understanding of small deviations lags significantly behind that of large deviations, as surveyed in Section~\ref{sec:gauss-proc-prob}.
A similar assessment was expressed by Ledoux in his Saint-Flour lecture notes~\cite{ledoux1996isoperimetry}.
Specifically, Chapter~7 of~\cite{ledoux1996isoperimetry}, which is devoted to small ball probabilities and centers on the results of Kuelbs-Li and Talagrand, opens as follows:
\emph{``While\dots the behavior of (the complement of) large balls
for Gaussian measures is relatively well described, small ball probabilities are much
less known. This problem has gone recently a quick development\dots
% and we intend to present in this chapter some significant recent results,
although it seems that there is still a long way to the final word (if there is any).''}

\paragraph{Main result.}

We now turn to the first main contribution of the present work, which is to provide an explicit characterization of general Gaussian small ball probabilities.

As a motivation for what follows, observe that Theorems~\ref{thm:kuelbs-li} and~\ref{thm:talagrand-sb-lower-cover} convey an ambivalent message: on the one hand, they suggest (and partly establish) a strong connection between small ball probabilities and covering numbers of the underlying index set $K$;
while at the same time, their limitations suggest that covering number are insufficient to obtain a full characterization.

The general characterization will indeed involve covering numbers, but (in contrast with the Kuelbs-Li conditions) will also depend on additional parameters: the \emph{local Gaussian widths}, namely the Gaussian width of the ``local sets'' $K \cap r B_2^n$ for $r > 0$.

We will first state the result in a form that allows comparison with the upper tail characterization of Corollary~\ref{cor:upper-tail-gaussian}.
Below, we denote by $N (K, r) = N (K, r B_2^n)$ the Euclidean covering number of $K$ at scale $r > 0$ (Definition~\ref{def:covering-numbers}).
To any symmetric convex body $K \subset \R^n$ and $u \in \ooint{0}{w (K)}$, we associate the following two ``effective radii'' at width $u$:
\begin{align}
  \label{eq:def-radius-widths}
  r_* (K, u)
  &= \inf \braces[\big]{r > 0 : w (K \cap r B_2) \geq u}
    \in \ccint{0}{\rad (K)}
    \, ; \\
  \wt r (K, u)
  &= \inf \braces[\big]{r > 0 : r \sqrt{\log N (K, r)} \geq u}
    \in \ccint{0}{+\infty}
    \label{eq:def-radius-cover}
    \, .
\end{align}
Above, we use the convention that $\inf \varnothing = + \infty$; while it follows from Lemma~\ref{lem:loc-width-bij} below that $r_* (K, u) \in \ooint{0}{\rad (K)}$ and that this number is the unique solution $r > 0$ of the equation $w (K \cap r B_2) = u$.

Covering numbers and local Gaussian widths capture different features of the set $K$.
Indeed, at a given scale $r > 0$, the local width $w (K \cap r B_2)$ only depends on the local set $K \cap r B_2$ and not on the behavior of $K$ at scales larger than $r$.
By contrast, the covering number $N (K, r)$ only depends on the global behavior of the set $K$ at scales larger than $r$, irrespective of the structure of the local set $K \cap r B_2$.

With these definitions in place, we have the following two-sided characterization of Gaussian small-ball probabilities:

\begin{theorem}
  \label{thm:gaussian-small-ball}
  There exist universal constants $c_0, c > 0$ such that the following holds.
  For every $n \geq 1$, symmetric convex body $K \subset \R^n$ and
  $u \in (0, w (K)/c_0)$,
  \begin{equation}
    \label{eq:gaussian-small-ball}
    \frac{c^{-1} u^2}{\min \{ r_* (K, c u)^2, \wt r (K, c u)^2 \}}
    \leq - \log \gamma_n (u K^\circ)
    \leq \frac{c u^2}{\min \set{r_* (K, u/c)^2, \wt r (K, u/c)^2}}
    \, .
  \end{equation}
\end{theorem}

Theorem~\ref{thm:gaussian-small-ball} shows that (\eg, letting $u = 1$ and $K = L^\circ$ where $L \subset \R^n$ is a closed symmetric convex set with nonempty interior) the Gaussian measure of a ``small'' symmetric convex set (that is, a constant-factor contraction of a symmetric convex set whose Gaussian measure is bounded away from $1$) is governed by the behavior of the covering numbers and local Gaussian widths of its polar body at varying scales.

The characterization of the lower tail of suprema of Gaussian processes of Theorem~\ref{thm:gaussian-small-ball} may be compared with the characterization of the upper tail of Corollary~\ref{cor:upper-tail-gaussian}.
Specifically, the radius $\rad (K)$ appearing in large deviations is replaced by the minimum of the two scale-dependent effective radii $r_* (K, u)$ and $\wt r (K, u)$.
% , determined respectively by local widths and covering numbers.
In particular, the behavior of lower deviations is more subtle than that of upper deviations: while the latter are governed by two numbers (the Gaussian width and the radius), the former are governed by two functions (covering numbers and local widths).

The proof of Theorem~\ref{thm:gaussian-small-ball} is provided in Section~\ref{sec:proof-small-ball}.
Interestingly, once the right quantities~\eqref{eq:def-radius-widths} and~\eqref{eq:def-radius-cover} are introduced, the proof is hardly more involved than the proof of the Kuelbs-Li inequalities (Theorem~\ref{thm:kuelbs-li}), with which it shares a common core.
Specifically, the left-most inequality in~\eqref{eq:gaussian-small-ball} is obtained by combining a modified version of the Kuelbs-Li bound (involving the quantity $\wt r (K, u)$ in order to avoid a recursive inequality) with a different bound obtained by applying Gaussian concentration to a suitable subset of $K$.
The right-most inequality also involves similar steps as the corresponding Kuelbs-Li condition, but critically controls the measure of certain sets using the information supplied by local widths, as opposed to the Gaussian isoperimetric inequality~\cite{borell1975brunn} in the approach of Kuelbs and Li~\cite{kuelbs1993metric}.

It may be instructive to reframe the two-sided estimate of Theorem~\ref{thm:gaussian-small-ball} in a form that allows for a more direct comparison with the Kuelbs-Li result.
Specifically, the Kuelbs-Li condition (which fails to hold in general, see Section~\ref{sec:lack-solution-kuelbs}; and even when it does, may fail to characterize the small ball probability) takes the form~\eqref{eq:implicit-sb-kuelbs-li}.
To reformulate the estimate~\eqref{eq:gaussian-small-ball}, define for any $u \in \ooint{0}{w(K)}$ the following ``effective dimensions'' at level $u$:
\begin{align}
  \label{eq:def-dim-width}
  s_* (K, u)
  &= \sup \set*{s > 0 : w \parens*{K \cap \frac{u}{\sqrt{s}} B_2} \geq u}
    \, ; \\
  \wt s (K, u)
  &= \sup \set*{s > 0 : \log N \parens*{K, \frac{u}{\sqrt{s}}} \geq s}
    \, ,
    \label{eq:def-dim-cover}
\end{align}
with the convention in~\eqref{eq:def-dim-cover} that $\sup \varnothing = 0$.
It is straightforward to show (see Lemma~\ref{lem:link-fixed-r-k-u}) that $u^2/r_* (K, u)^2 = s_* (K, u)$ and $u^2/\wt r (K, u)^2 = \wt s (K, u)$.
The characterization of Theorem~\ref{thm:gaussian-small-ball} then takes the following form: for some universal constants $c_0, c \geq 1$, for every symmetric convex body $K \subset \R^n$ and $u \in \ooint{0}{w(K)/c_0}$ one has
\begin{equation}
  \label{eq:gaussian-small-ball-dim-1}
  c^{-1} \, \ol s \parens*{K, c u}
  \leq - \log \gamma_n (u K^\circ)
  \leq c \, \ol s \parens*{K, u/c}
  \, , 
\end{equation}
where for $u > 0$ we let $\ol s (K, u) = \max \set*{s_* \parens*{K, u}, \wt s \parens*{K, u}}$, so that
\begin{equation}
  \label{eq:gaussian-small-ball-dim-2}
  \ol s \parens*{K, u}
  = \sup \braces*{s > 0 : \max \bracks*{w \parens*{\frac{\sqrt{s}}{u} K \cap B_2^n}^2, \log N \parens*{K, \frac{u}{\sqrt{s}}}} \geq s}
  \, .  
\end{equation}

We note that both covering numbers and local widths are indeed required in general, as either $s_* (K, u)$ or $\wt s (K, u)$ may dominate depending on the value of $u$.
At one extreme, consider the case of a fixed $K \subset \R^n$ with nonempty interior, while $u \to 0^+$.
Then the inequality $w (K \cap (u/\sqrt{s}) B_2^n) \leq (u/\sqrt{s}) \sqrt{n}$ implies that $s_* (K, u) \leq n$ for any $u > 0$.
Meanwhile, the fact that $\log N (K, r)$ diverges as $n \log (1/r)$ as $r \to 0$ implies that $\wt s (K, u)$ diverges as $n \log (1/u)$ as $u \to 0$.
Hence, $\wt s (K, u)$ dominates $s_* (K, u)$ in the limit of low scales in finite dimension.
% We discuss the opposite situation next.

\paragraph{Constant-factor lower deviations.}

Consider now the important special regime where $u$ is of order $w (K)$; that is, deviations of a Gaussian process by a constant factor below its expectation.
The next result shows that in this regime, small ball probabilities are solely characterized by local widths.

\begin{corollary}
  \label{cor:small-ball-constant-simple}
  There exists a universal constant $c \geq 1$ such that the following holds.
  For every $n \geq 1$, symmetric convex body $K \neq \set{0}$ in $\R^n$ and $u \in \ooint{0.01 w (K)}{0.99 w (K)}$, one has
  \begin{equation}
    \label{eq:small-ball-constant-simple}
    \frac{c^{-1} w (K)^2}{r_* \parens*{K, 1.01 u}^2}
    \leq - \log \gamma_n \parens*{u K^\circ}
    \leq \frac{c \, w (K)^2}{r_* \parens*{K, 0.99 u}^2}
    \, .
  \end{equation}
\end{corollary}

The proof of Corollary~\ref{cor:small-ball-constant-simple} is provided in Section~\ref{sec:const-fact-lower}.
Specifically, it is obtained by combining lemmas in the proof of Theorem~\ref{thm:gaussian-small-ball} (rather than its statement) together with Sudakov's inequality (recalled as Lemma~\ref{lem:sudakov}); the latter implies that the term depending on local widths dominates the one depending on covering numbers in the regime where $u \asymp w (K)$.

% Hence, it follows from Corollary~\ref{cor:upper-tail-gaussian} and Corollary~\ref{cor:small-ball-constant-simple} that the probability that the supremum of the Gaussian process indexed by $K$ exceeds its expectation by a constant factor is exponentially small in $w (K)^2/\rad(K)^2$, while the probability of a constant factor lower deviation is exponentially small in $w (K)^2/r_* (K, w (K)/c)^2$.
Since for $r_* (K, w (K)/c) \leq \rad (K)$ for $c > 1$, Corollaries~\ref{cor:upper-tail-gaussian} and~\ref{cor:small-ball-constant-simple} imply that constant-factor lower deviations of suprema of Gaussian processes (relative to their expectation) are less likely (up to constants in the exponent) than constant-factor upper deviations.
Example~\ref{ex:operator-norm} below provides a relevant example where a gap between upper and lower deviations occurs.

\begin{example}
  \label{ex:operator-norm}
  Let $n \geq 1$ and $G_n$ be an $n \times n$ random matrix with \iid standard Gaussian entries, and let $\opnorm{A} = \sup_{u \in B_2^n} \norm{A u}$ denote the operator norm of a matrix $A \in \R^{n \times n}$.
  Then
  \begin{equation}
    \label{eq:expect-op-norm}
    \sqrt{\frac{2 n}{\pi}} 
    \leq \Expect*{\opnorm{G_n}}
    \leq 2 \sqrt{n}
    \, .
  \end{equation}
  Then, for some absolute constant $c \geq 1$, one has
  \begin{align}
    \label{eq:operator-norm-upper-dev}
    \exp \parens*{- c n}
    \leq \Probab*{\opnorm*{G_n} \geq 4 \sqrt{n}}
    \leq \exp \parens*{- \frac{n}{c}}
    \, ; \\
    \label{eq:operator-norm-lower-dev}
    \exp \parens*{- c n^2}
    \leq \Probab*{\opnorm*{G_n} \leq \frac{\sqrt{n}}{4}}
    \leq \exp \parens*{- \frac{n^2}{c}}
    \, .
  \end{align}
\end{example}

We refer to Section~\ref{sec:exampl-oper-norm} for a justification of the claims of Example~\ref{ex:operator-norm}.

\subsection{Metric characterization and comparison inequalities}
\label{sec:chaining-comparison}

\paragraph{Chaining characterization.}

The characterization of small ball probabilities of Theorem~\ref{thm:gaussian-small-ball} depends on two geometric quantities that can both be connected to the metric structure of the underlying set $K$.
The first such quantity is the covering numbers, which are defined in terms of the metric structure of $K$.
The second quantity consists in the local Gaussian widths (which, in probabilistic terms, describe the expected modulus of continuity of the underlying Gaussian process), which can be connected to the local metric structure of the set via the Majorizing Measures theorem (see~\cite[Section~2.15]{talagrand2021upper} for an exact correspondence).

It is thus reasonable to expect that the characterization of Theorem~\ref{thm:gaussian-small-ball} could be combined with the Majorizing Measures theorem to obtain an explicit description of the lower tail of Gaussian processes in terms of the metric structure of the (convexified) index set.
% Such a description would complement the metric descriptions of the expected and large values of Gaussian processes of, respectively, Theorem~\ref{thm:talagrand-maj-meas} and Corollary~\ref{cor:upper-tail-gaussian}.
In order to provide such a metric characterization, we will need the following definition:

\begin{definition}[Truncated generic chaining functional]
  \label{def:admissible-sequence}
  Let $T \subset \R^n$ be nonempty.
  For any $p \in \N$, we define the $p$-th \emph{truncated $\gamma_2$ functional} of $T$ by
  \begin{equation}
    \label{eq:truncated-gamma-2}
    \gamma_2^{(p)} (T)
    = \inf_{(\A_j)_{j \geq 0}} \sup_{t \in T} \sum_{j \geq p} 2^{j/2} \diam A_j (t)
    \, ,
  \end{equation}
  where the infimum is taken over all admissible sequences of partitions of $T$ (see Definition~\ref{def:adm-seq-gen-chaining}).
  (The quantity inside the infimum only depends on $(\A_j)_{j \geq p}$.)
  In particular, $\gamma_2^{(0)} (T) = \gamma_2 (T)$.
\end{definition}

Now, for $s \in \R_+^*$ and a real random variable $Y$, we let
% $Q (s, Y) = \sup \braces*{q \in \R : \Probab*{Y \leq q} \leq \exp(-s)}$
\begin{equation*}
  Q (s, Y)
  = Q_- (s, Y)
  = \sup \braces*{q \in \R : \Probab*{Y \leq q} \leq \exp(-s)}
\end{equation*}
denote the quantile of order $\exp(-s)$ of $Y$.
% Alternatively, let
% \begin{equation}
%   \label{eq:def-lower-quantile}
%   Q (s, K)
%   = \sup \braces*{q \in \R : \Probab*{\sup_{t \in K} \innerp{t}{X} \leq q} \leq \exp(-s)}
% \end{equation}
% $Q_{2^{p}} (h_K (X)) = Q_{2^{p}} (\sup_{v \in K} \innerp{v}{X}) = Q_{2^{p}} (\norm{X}_{K^\circ})$
The next result shows that small quantiles of suprema of Gaussian processes can be described in terms of the truncated generic chaining functional.
This complements the known metric descriptions of the expectation and large quantiles of Gaussian processes (Theorem~\ref{thm:talagrand-maj-meas} and Corollary~\ref{cor:upper-quantiles-gaussian}).

\begin{corollary}
  \label{cor:quantiles-chaining}
  There exists a universal constant $c > 0$ such that the following holds.
  For every $n \geq 1$, convex body $K \subset \R^n$ and $p \in \N$, letting $X \sim \gaussdist (0, I_n)$ we have
  \begin{equation}
    \label{eq:quantiles-chaining}
    \frac{1}{c} Q \parens[\bigg]{c 2^p, \sup_{t, t' \in K} \abs*{\innerp{t-t'}{X}}}
    \leq \gamma_2^{(p)} (K)
    \leq c \, Q \parens[\bigg]{2^p/c, \sup_{t, t' \in K} \abs*{\innerp{t-t'}{X}}}
    \, .
  \end{equation}
\end{corollary}

In a more geometric language, Corollary~\ref{cor:quantiles-chaining} implies in particular that \emph{the Gaussian measure of a ``small'' symmetric closed convex set can be explicitly characterized (in an approximate sense) in terms of the metric structure of its polar}.
Here, a closed symmetric convex set $L \subset \R^n$ is said to be small if $\gamma_n (c L) \leq 1/2$ for some large universal constant $c$.
In addition, the characterization is approximate in the same sense as before: namely, it holds up to universal constant factors in the contraction or dilation of $L$ and in the exponent of the probability.
Specifically, Corollary~\ref{cor:quantiles-chaining} relates the mass of dilates $\gamma_n (u L)$ (for $u \in \ooint{0}{u_0/c}$, with $u_0 > 0$ such that $\gamma_n (u_0 L) = 1/2$, hence $u_0 \asymp w (L^\circ) \asymp \gamma_2 (L^\circ)$) to the behavior of the sequence $(\gamma_2^{(p)} (L^\circ))_{p \in \N}$ of truncated generic chaining functionals of its polar body $L^\circ$.

We refer to Section~\ref{sec:chain-char-lower} for the proof of Corollary~\ref{cor:quantiles-chaining}.
While it is likely that Corollary~\ref{cor:quantiles-chaining} could be proven by combining Theorems~\ref{thm:talagrand-maj-meas} and~\ref{thm:gaussian-small-ball}, within the context of this work we found it more economical to deduce it from a connection with the Wills functional (presented in the next section) together with a characterization of the Wills functional involving truncated generic chaining functionals established in~\cite[Corollary~4.1]{mourtada2025universal}.

\paragraph{Comparison inequality.}

From Corollary~\ref{cor:quantiles-chaining} (and Corollary~\ref{cor:upper-quantiles-gaussian}), we deduce a general comparison inequality for the tails of suprema of Gaussian processes indexed by convex sets.

\begin{corollary}
  \label{cor:comparison-small-ball}
  There exists a universal constant $c \geq 1$ such that the following holds.
  Let $n \geq 1$ and $K_1, K_2 \subset \R^n$ be convex bodies.
  Assume that there exists a map $\phi : K_1 \to K_2$ that is expansive \textup(that is, $\norm{\phi (t) - \phi (t')} \geq \norm{t - t'}$ for all $t, t' \in K_1$\textup).
  % Then, for every $s \in \R_+^*$, letting $X \sim \gaussdist (0, I_n)$,
  % \begin{equation}
  %   \label{eq:comparison-quantiles}
  %   Q \parens*{c s, \sup_{t \in K_1} \innerp{t}{X}}
  %   \leq  c \, Q \parens*{s, \sup_{t \in K_2} \innerp{t}{X}}
  %   \, .
  % \end{equation}
  % Equivalently, for every $u \in \R^+$,
  Then, letting $X \sim \gaussdist(0, I_n)$, for every $u \in \R^+$, 
  \begin{equation}
    \label{eq:comparison-cdf}
    \Probab[\bigg]{\sup_{t,t' \in K_2} \abs{\innerp{t-t'}{X}} \leq u}
    \leq \Probab[\bigg]{\sup_{t,t' \in K_1} \abs{\innerp{t-t'}{X}} \leq c u}^{1/c}
    \, .
  \end{equation}
\end{corollary}

We refer to Section~\ref{sec:comp-ineq-conv} for the proof of Corollary~\ref{cor:comparison-small-ball}.
Corollary~\ref{cor:comparison-small-ball} provides a comparison inequality that holds for the full cumulative distribution functions of suprema of Gaussian processes indexed by convex sets.
This complements the Fernique-Sudakov comparison inequality~\cite{fernique1975regularite,sudakov1971gaussian} that holds for the expected value, without convexity assumptions; the latter can be deduced (up to universal constant factors) from the Majorizing Measures theorem.

Perhaps the closest relative to Corollary~\ref{cor:comparison-small-ball} is a classical result of Slepian~\cite{slepian1962one} (see \eg \cite[Theorem~7.2.1 p.~151]{vershynin2018high}), which can equivalently be stated as follows:

\begin{theorem}[\cite{slepian1962one}]
  \label{thm:slepian}
  Let $n \geq 1$ and $T_1, T_2$ be compact subsets of $\R^n$.
  Assume that there exists a map $\phi : T_1 \to T_2$ that is expansive \textup(in the sense of Corollary~\ref{cor:comparison-small-ball}\textup) and norm-preserving \textup(that is, $\norm{\phi (t)} = \norm{t}$ for every $t \in T_1$\textup).
  Then, letting $X \sim \gaussdist(0, I_n)$, for every $u \in \R$, 
  \begin{equation}
    \label{eq:comparison-slepian}
    \Probab[\bigg]{\sup_{t \in T_2} \innerp{t}{X} \leq u}
    \leq \Probab[\bigg]{\sup_{t \in T_1} \innerp{t}{X} \leq u}
    \, .
  \end{equation}
\end{theorem}

While both Corollary~\ref{cor:comparison-small-ball} and Theorem~\ref{thm:slepian} provide comparison inequalities that bear on the cumulative distribution functions, the two statements present significant differences.
On the one hand, the comparison inequality of Theorem~\ref{thm:slepian} holds with constant $c = 1$ and does not require convexity of the underlying index sets $T_1, T_2$.
On the other hand, the assumption under which Slepian's inequality holds is significantly more restrictive, as it is assumed that the expansive map $\phi$ is norm-preserving, which is not required in Corollary~\ref{cor:comparison-small-ball}.

\paragraph{The role of convexity.}

We now conclude this discussion by commenting on the role of convexity of the index set.
Indeed, the metric characterization of the lower tail of Corollary~\ref{cor:quantiles-chaining} and the comparison inequality of Corollary~\ref{cor:comparison-small-ball} are stated only for convex index sets.
In contrast, known characterizations of expected and large values of Gaussian processes (Theorem~\ref{thm:talagrand-maj-meas} and Corollary~\ref{cor:upper-quantiles-gaussian}) hold for general non-convex index sets.
Indeed, as emphasized by Talagrand (see \eg the discussions in~\cite[Section~1.7.1 p.~10 and Section~2.16 p.~84]{talagrand2021upper}), the realization that expected suprema of Gaussian processes only depend on properties of the index set as a metric space, without using additional properties of subsets of Hilbert spaces or requiring convexity, was a key conceptual progress in the understanding of Gaussian processes.

Proposition~\ref{prop:no-comparison-nonconvex} below (proved in Section~\ref{sec:lack-comp-ineq}) shows that unlike for the upper tail, convexity plays an essential role in the description of the lower tail of Gaussian processes.
In particular, the convexity assumptions in Corollaries~\ref{cor:quantiles-chaining} and~\ref{cor:comparison-small-ball} cannot be dropped, as the lower tail of a Gaussian process indexed by a nonconvex set does not depend monotonously on the pairwise distances between elements of this set.

\begin{proposition}
  \label{prop:no-comparison-nonconvex}
  There exists a sequence $(c_n)_{n \geq 1}$ of real numbers with $c_n \to +\infty$ as $n \to \infty$ such that the following holds.
  For any $n \geq 1$, there exist finite symmetric sets $T_{1,n}, T_{2,n} \subset \R^n$ for which there exists an expansive map $\phi_n : T_{1, n} \to T_{2, n}$
  and $u_n \in \R_+^*$ such that
  \begin{equation}
    \label{eq:no-comparison-nonconvex}
    \Probab*{\sup_{t \in T_{2,n}} \innerp{t}{X} \leq u_n}
    > \Probab*{\sup_{t \in T_{1,n}} \innerp{t}{X} \leq c_n u_n}^{1/c_n}
    \, .
  \end{equation}
\end{proposition}

\subsection{Intrinsic volumes and the Wills functional}
\label{sec:intr-volum-wills}

The previous discussion established that convexity plays a key role in the description of small ball probabilities.
Building on this insight and in order to explore more concrete connections between Gaussian processes and convexity, in this section we turn our attention to other classical quantities in Convex Geometry, namely the intrinsic volumes of convex bodies.

\paragraph{Intrinsic volumes.}

Intrinsic volumes (or equivalently \emph{quermassintegrals}, up to changes in normalization and indexing) are fundamental geometric functionals associated to convex bodies.
They play a central role in the Brunn-Minkowski theory of convex bodies, comprehensively exposed in Schneider's book~\cite{schneider2013convex}, and feature prominently in integral and stochastic geometry, notably through kinematic formulae~\cite{schneider2008stochastic,klain1997introduction}.
Their close relationship with curvature measures also connects them to geometric measure theory and differential geometry~\cite{federer1959curvature}.

One standard way to define the intrinsic volumes is through Steiner's formula for the volume of Euclidean enlargements of convex sets (see, \eg,~\cite[p.~208]{schneider2013convex}).
For every convex body $K \subset \R^n$ and $r \in \R^+$, this formula states that
\begin{equation}
  \label{eq:steiner}
  \vol_n \parens*{K + r B_2^n}
  = \sum_{j=0}^n V_j (K) \vol_{n-j} (B_2^{n-j}) r^{n-j}
  \, ,
\end{equation}
where $\vol_n$ denotes the $n$-dimensional Haussdorff (Lebesgue) measure on $\R^n$.
In particular, $V_0 (K) = 1$, $V_n (K) = \vol_n (K)$ and (if $K$ has nonempty interior) $V_{n-1} (K) = \frac{1}{2} \vol_{n-1} (\partial K)$ is half the surface area of the boundary of $K$.
Moreover, Hadwiger's theorem~\cite[Theorem~6.4.14 p.~360]{schneider2013convex} states that the functionals $V_0, V_1, \dots, V_n$ form a basis of the vector space of continuous valuations on convex bodies that are invariant under affine isometries of $\R^n$.

A more explicit representation of intrinsic volumes is provided by Kubota's formula, which expresses $V_j (K)$ in terms of the expected $j$-dimensional volume of the projection of $K$ onto a uniformly random $j$-dimensional subspace of $\R^n$.
A convenient variant of this representation, due to Tsirel'son~\cite{tsirelson1986geometric}, replaces orthogonal projections by Gaussian linear maps.
Specifically, let $X_1, X_2, \dots$ denote an \iid sequence of standard Gaussian random variables in $\R^n$.
For $j \geq 1$, define the Gaussian linear map $G_j : \R^n \to \R^j$ by
\begin{equation}
  \label{eq:def-gaussian-matrix}
  G_j (t)
  = \parens*{\innerp{t}{X_1}, \dots, \innerp{t}{X_j}}
  \, , \quad t \in \R^n
  \, .
\end{equation}
Then (setting $V_j (K) = 0$ for $j > n$), for any convex body $K \subset \R^n$ we have
\begin{equation}
  \label{eq:intrinsic-tsirelson}
  V_j \parens*{K}
  = \frac{(2\pi)^{j/2}}{j!} \cdot \frac{\Expect*{\vol_j (G_j K)}}{\vol_j (B_2^j)}
  \, .
\end{equation}
In particular, this Gaussian representation makes it clear that intrinsic volumes do not depend on the dimension $n$ of the ambient space, in the sense that they are preserved under isometric embeddings into higher-dimensional Euclidean spaces.
An important special case is Sudakov's identity~\cite{sudakov1976geometric}, which relates the first intrinsic volume to the Gaussian width:
\begin{equation}
  \label{eq:sudakov-first-intr-width}
  V_1 \parens[\big]{K/\sqrt{2\pi}}
  = w (K)
  \, .
\end{equation}

A fundamental property of the sequence of intrinsic volumes of a convex body is that it is ultra-log-concave (that is, log-concave with respect to the Poisson measure on $\N$; see Definition~\ref{def:ulc}), in the sense that it contains no internal zeroes and that for every convex body $K \subset \R^n$ and $j \geq 1$,
\begin{equation}
  \label{eq:ulc-intrinsic}
  j V_j (K)^2 \geq (j+1) V_{j-1} (K) V_{j+1} (K)
  \, .
\end{equation}
The inequality~\eqref{eq:ulc-intrinsic}, of isoperimetric nature, was obtained independently by Chevet~\cite{chevet1976processus} and McMullen~\cite{mcmullen1991inequalities}, and follows from the Alexandrov-Fenchel inequalities for mixed volumes (see~\cite[Section~7.3]{schneider2013convex} and~\cite[\S20]{burago1988geometric}).

\paragraph{The Wills functional and small ball probabilities.}

In this section, we will make special use of
the Wills functional~\cite{wills1973gitterpunktanzahl,hadwiger1975wills} defined as\footnote{
  The Wills functional is usually defined as the sum of the intrinsic volumes, without the rescaling by $1/\sqrt{2\pi}$. We adopt this normalization to simplify its connection with Gaussian quantities; all results presented here can be translated to the standard convention with minor adjustments.}
\begin{equation}
  \label{eq:def-wills}
  W (K)
  = \sum_{j=0}^n V_j \parens[\big]{K/\sqrt{2\pi}}
  \, .
\end{equation}
This functional was introduced by Wills~\cite{wills1973gitterpunktanzahl} in connection with lattice point enumeration;
it was later shown in~\cite{mourtada2025universal} to govern the complexity of universal coding under Gaussian noise.

The Wills functional bears some partial connection with the Gaussian width.
Indeed, the ultra-log-concavity condition~\eqref{eq:ulc-intrinsic} and the identity $V_0 (K) = 1$ imply that $V_j (K) \leq V_1 (K)^j/j!$ for all $j \geq 1$.
Together with Sudakov's identity~\eqref{eq:sudakov-first-intr-width}, this yields
\begin{equation}
  \label{eq:wills-width-mcmullen}
  \log \parens*{1 + w (K)}
  \leq \log W (K)
  \leq w (K)
  \, .
\end{equation}
The upper bound was obtained independently by Tsirel'son~\cite{tsirelson1986geometric} and McMullen~\cite{mcmullen1991inequalities}; alternative proofs that do not rely on the Alexandrov-Fenchel inequalities can be found in \cite{alonso2021further,mourtada2025universal}.
In fact, the inequalities~\eqref{eq:wills-width-mcmullen} are the best dimension-free upper and lower bounds on $W (K)$ in terms of $w (K)$ (see \eg the comments around \cite[eq.~(3.6)]{mourtada2025universal}).
% \footnote{The key ideas behind the present work were formed during and following the preparation of this prior work.}

The starting point for the results of this section is a two-sided characterization of the logarithm of the Wills functional established by the author in~\cite{mourtada2025universal}.
In order to state it, we define the following quantities associated to a convex body $K \subset \R^n$:
\begin{align}
  \label{eq:prev-fixed-widths}
  \rho_* (K)
  &= \sup \Big\{ \rho \geq 0 : w \big( K \cap \rho B_2 \big) \geq \rho^2 \Big\}
    \, , \\
  \wt \rho (K)
  &= \sup \Big\{ \rho \geq 0 : \log N (K, \rho) \geq \rho^2 \Big\}
    \label{eq:prev-fixed-cov}
    \, .
\end{align}
For symmetric convex bodies, the characterization takes the following form:

\begin{theorem}[\cite{mourtada2025universal}, Theorem~4.2]
  \label{thm:mourt-wills}
  There exists a universal constant $c > 0$ such that, for every $n \geq 1$ and every symmetric convex body $K \subset \R^n$, 
  \begin{equation}
    \label{eq:wills-width-cover}
    c^{-1} \max \braces*{\rho_* (K)^2, \wt \rho (K)^2}
    \leq \log W (K)
    \leq c \max \braces*{\rho_* (K)^2, \wt \rho (K)^2}
    \, .
  \end{equation}
\end{theorem}

We note in passing that the three quantities $\rho_* (K)^2$, $\wt \rho (K)^2$ and $\log W (K)$ appearing in~\eqref{eq:wills-width-cover}
admit precise information-theoretic interpretations~\cite[Proposition~4.1, Theorem~4.3 and Theorem~2.1]{mourtada2025universal}.
Through these correspondences, Theorem~\ref{thm:mourt-wills} can be reframed as a property of universal coding under Gaussian noise~\cite[Corollary~4.2]{mourtada2025universal}.

Now, observe that the quantitative characterization of small ball probabilities in Theorem~\ref{thm:gaussian-small-ball} and that of the Wills functional in Theorem~\ref{thm:mourt-wills} both involve a combination of covering numbers and local Gaussian widths at various scales.
Combining these two results, and relating the effective radii $r_* (K, u), \wt r (K, u)$ to the fixed points $\rho_* (\lambda K), \wt \rho (\lambda K)$ associated to dilates of $K$, we obtain the following correspondence (proved in Section~\ref{sec:corr-with-wills}).

\begin{corollary}
  \label{cor:link-sb-wills}
  For every $n \geq 1$, symmetric convex body $K \subset \R^n$ and $u \in (0, w(K))$, define
  \begin{equation}
    \label{eq:def-eff-dim-wills}
    s_W \parens*{K, u}
    = \sup \braces*{s > 0 : \log W \parens*{\frac{s}{u} K} \geq s}
    \, .
  \end{equation}
  Then,
  there exist universal constants $c_0, c \geq 1$ such that,
  for every $u \in \ooint{0}{w (K)/c_0}$,
  \begin{equation}
    \label{eq:link-sb-wills}
    c^{-1} s_W \parens*{K, c \, u}
    \leq - \log \gamma_n (u K^\circ)
    \leq c \, s_W \parens*{K, u/c}
    \, .
  \end{equation}
\end{corollary}

Corollary~\ref{cor:link-sb-wills} establishes a correspondence between Gaussian small ball probabilities and the behavior of the Wills functional of dilates of $K$, 
namely with the quantities
\begin{equation}
  \label{eq:wills-mgf-intrinsic}
  \log W (\lambda K)
  % = \log \parens[\bigg]{\sum_{j=0}^n V_j \parens[\big]{K} \parens*{\frac{\lambda}{\sqrt{2\pi}}}^j}
  = \log \parens[\bigg]{\sum_{j=0}^n V_j \parens[\big]{K/\sqrt{2\pi}} \lambda^j}
  \, , \quad \lambda \in \R^+
  \, .
\end{equation}

\paragraph{Quantitative characterization of intrinsic volumes.}

We now turn to the quantitative behavior of individual intrinsic volumes, rather than their sum in the Wills functional.
Since the $j$-th intrinsic volume $V_j$ is positively homogeneous of degree $j$, comparing intrinsic volumes across different indices calls for a suitable normalization.
We therefore define, for $j \geq 1$, the $j$-th \emph{intrinsic width} $\theta_j (K)$ of a convex body $K \subset \R^n$ by
\begin{equation}
  \label{eq:def-intrinsic-width}
  \theta_j (K)
  = \parens*{j! \cdot V_j \parens[\big]{K/\sqrt{2\pi}}}^{1/j}
  \, .
\end{equation}
(Although this terminology is not standard, Paouris, Pivovarov and Valettas~\cite{paouris2019quantitative} consider the closely related quantities $(V_j (K)/V_j(B_2^n))^{1/j} \asymp \theta_j (K) / \sqrt{n}$.)
By Tsirel'son's formula~\eqref{eq:intrinsic-tsirelson}, $\theta_j (K)$ is the radius of the $j$-dimensional Euclidean ball whose volume equals the expected volume of the Gaussian image $G_j K$, namely
\begin{equation}
  \label{eq:intr-width-interpret}
  \vol_j \parens[\big]{\theta_j (K) B_2^j}
  = \Expect[\big]{\vol_j \parens*{G_j K}}
  \, .
\end{equation}
Likewise, Sudakov's identity~\eqref{eq:sudakov-first-intr-width} shows that $\theta_1 (K) = w(K)$.
Finally, the ultra-log-concavity property~\eqref{eq:ulc-intrinsic} implies that the sequence $(j! V_j (K/\sqrt{2\pi}))_{j \geq 0}$ is log-concave.
Since $V_0 (K/\sqrt{2\pi}) = 1$, it follows that the intrinsic widths $\theta_j(K)$ are decreasing in $j \geq 1$.

Theorem~\ref{thm:intrinsic-volumes} below is the main result of this section.
It provides a two-sided quantitative characterization of the behavior of the sequence of intrinsic volumes of a convex body, expressed in terms of its local Gaussian widths and covering numbers.

\begin{theorem}
  \label{thm:intrinsic-volumes}
  There exist universal constants $c_0, c_1, c_2 \geq 1$ such that the following holds.
  For every integer $n \geq 1$, symmetric convex body $K \subset \R^n$, integer $j \geq 1$ and real $u > 0$:
  \begin{enumerate}
  \item if $\theta_j (K) \geq u$, then for every $\xi > 1$,
    \begin{equation}
      \label{eq:intrinsic-upper}
      \log N \parens*{K, \frac{u}{c_1 \xi \sqrt{j}}}
      \geq \frac{j \log \xi}{c_1}
      \quad \text{ or } \quad
      w \parens*{K \cap \frac{u}{c_1 \xi \sqrt{j}} B_2}
      \geq \frac{\sqrt{\log \xi}}{c_2 \xi} \, u
      \, ;
    \end{equation}
  \item conversely, if $j \geq c_0$ and, for every $\xi > 1$,
    \begin{equation}
      \label{eq:intrinsic-lower-cond}
      \log N \parens*{K, \frac{u}{\xi \sqrt{j}}}
      \geq j \log \xi
      \quad \text{ or } \quad
      w \parens*{K \cap \frac{u}{\xi \sqrt{j}} B_2}
      \geq \frac{\sqrt{\log \xi}}{\xi} \, u
      \, ,
    \end{equation}
    then $\theta_{\floor{j/c_0}} (K) \geq u/c_1$.
  \end{enumerate}
\end{theorem}

We refer to Sections~\ref{sec:recovering-an-ultra} and~\ref{sec:two-sided-estimates} for the proof of Theorem~\ref{thm:intrinsic-volumes}.
The idea is to use the ultra-log-concavity property of intrinsic volumes to approximately recover them from the Wills functional, as then resort to the characterization of the latter in Theorem~\ref{thm:mourt-wills}.

\begin{remark}[Non-symmetric case]
  \label{rem:intrinsic-non-sym}
  While the characterization of Theorem~\ref{thm:intrinsic-volumes} is stated for {symmetric} convex bodies, it extends to non-symmetric convex bodies $K \subset \R^n$ with the following modification: in the second condition of~\eqref{eq:intrinsic-upper} and~\eqref{eq:intrinsic-lower-cond}, the quantity $w (K \cap r B_2)$ should be replaced by $w ((K - K) \cap r B_2)$.

  Indeed, the Rogers-Shephard inequality~\cite{rogers1957difference} (see \cite[Theorem 10.1.4 p. 530]{schneider2013convex}) on the volume of difference bodies shows that $V_j (K) \leq V_j (K - K) \leq 4^j V_j (K)$ for every $j \geq 1$, thus $\theta_j (K) \leq \theta_j (K-K) \leq 4 \theta_j (K)$.
  Applying Theorem~\ref{thm:intrinsic-volumes} to the symmetric convex body $K-K$ and using that $\log N (K, r) \leq \log N (K-K, r) \leq 2 \log N (K, r/2)$ gives the desired characterization.
\end{remark}

We note that a recent work of Paouris and Pathak~\cite{paouris2026optimal}, which appeared during the final stages of the preparation of the present work, provides an upper bound on intrinsic volumes under a condition on the behavior of local widths.
Specifically, using our notation and comparable normalization and for symmetric $K$, \cite[Proposition~3.5]{paouris2026optimal} (whose proof also relies on~\cite[Corollary~4.3]{mourtada2025universal} together with additional estimates) reads as follows.
If for some constant $\kappa_0 > 0$,
\begin{equation}
  \label{eq:condition-reese}
  w \parens*{K \cap r B_2^n}
  \leq \kappa_0 \sqrt{n \log (e r)}
  \qquad \text{for all } r \geq 1
  \, ,
\end{equation}
then there exists a constant $\kappa_1 > 0$ depending only on $\kappa_0$ such that, for every $j \in \set{1, \dots, n}$,
\begin{equation}
  \label{eq:upper-intrinsic-reese}
  \theta_j (K)
  \leq \kappa_1 \sqrt{n \log (e n/j)}
  \, .
\end{equation}
This estimate, together with analysis from~\cite{paouris2026optimal} of local widths under suitable affine position, is then used to sharply bound~\cite[Theorem~1.3]{paouris2026optimal} the infimum over invertible linear transformations $T$ of $\R^n$ of the quantity
$\max \braces*{\theta_j (T K)/\vol_n (T K)^{1/n}, \theta_j ((T K)^\circ)/\vol_n ((TK)^\circ)^{1/n}}$.

The result of Paouris and Pathak and Theorem~\ref{thm:intrinsic-volumes} are of a different nature and also quantitatively different.
Indeed, the estimate~\eqref{eq:upper-intrinsic-reese} is an upper bound on intrinsic widths under the growth profile for local widths of the condition~\eqref{eq:condition-reese}; while Theorem~\ref{thm:intrinsic-volumes} provides a two-sided approximate characterization of intrinsic widths for general convex bodies.

A special case of interest of the behavior of intrinsic widths is that of constant factor decrease relative to the Gaussian width.
Specifically, recall that the Alexandrov-Fenchel inequalities imply that
\begin{equation}
  \label{eq:decrease-intrinsic-alexandrov}
  w (K)
  = \theta_1 (K)
  \geq \theta_2 (K)
  \geq \cdots \geq \theta_n (K)
  \, .
\end{equation}
A natural question is to characterize the critical dimension $j \geq 1$ up to which the Alexandrov-Fenchel-type inequality $\theta_j (K) \leq w (K)$ can be reversed, namely $\theta_j (K) \gtrsim w (K)$.
Corollary~\ref{cor:reverse-alexandrov} (proved in Section~\ref{sec:proof-reverse-alexandrov}) shows that this critical dimension can be characterized in terms of the behavior of local widths.

\begin{corollary}
  \label{cor:reverse-alexandrov}
  For any $\kappa_0, \kappa_1 > 1$, there exists a constant $\kappa_2 > 1$ depending on $\kappa_0$ and $\kappa_3, \kappa_4 > 1$ depending on $\kappa_1$ such that the following holds.
  For any $n \geq 1$ and symmetric convex body $K \neq \set{0}$ in $\R^n$:
  \begin{enumerate}
  \item For any $j > 0$, if
    \begin{equation}
      \label{eq:reverse-alexandrov-small-k}
      j
      \leq \frac{w (K)^2}{r_* (K, w (K)/\kappa_0)^2}
      \, ,
    \end{equation}
    then $\theta_{\ceil{j/\kappa_2}} (K) \geq w (K)/\kappa_2$.
  \item Conversely, for any integer $j \geq 1$ such that $\theta_j (K) \geq w (K) / \kappa_1$, we have
    \begin{equation}
      \label{eq:reverse-alexandrov-cons-large-width}
      j
      \leq \frac{\kappa_3 \, w (K)^2}{r_* (K, w (K)/\kappa_4)^2}
      \, .
    \end{equation}
  \end{enumerate}
\end{corollary}

The study of quantitative reversals of the Alexandrov-Fenchel-type inequalities~\eqref{eq:decrease-intrinsic-alexandrov} between intrinsic volumes was initiated by Paouris, Pivovarov and Valettas~\cite{paouris2019quantitative}.
Their key result in this spirit~\cite[Theorem~1.1]{paouris2019quantitative} shows the following (after re-framing it in a dimension-free manner by letting the ambient dimension $n \to \infty$ while re-normalizing): for some parameter $\jppv (K)$ defined in terms of the Gaussian process indexed by $K$ and universal constant $c > 0$, for every $1 \leq j \leq c \jppv (K)$,
\begin{equation}
  \label{eq:ppv-reversal-alexandrov}
  \theta_j (K)
  \geq \parens*{1 + c \sqrt{\frac{j}{\jppv (K)} \log \parens*{\frac{e \, \jppv (K)}{j}}}}^{-1} w (K)
  \, .
\end{equation}
In particular, $c \, \jppv (K)$ constitutes a lower bound on the critical dimension for the reversal $\theta_j (K) \gtrsim w (K)$.
A remarkable feature of the dimension $\jppv (K)$ is that it already often exceeds the Dvoretzky-Milman dimension $\jdm (K) = (w (K)/\rad (K))^2$ (see Section~\ref{sec:low-dim-sections}).
At the same time, the lower bound $\jppv (K)$ on the critical dimension for constant-factor reversals is not always sharp.
For instance, if $K$ is the unit ball of the Schatten-$1$ norm on $n \times n$ matrices (see Example~\ref{ex:operator-norm} and Section~\ref{sec:exampl-oper-norm}), then $\jdm (K) \asymp n$ and $\jppv (K) \asymp n^{4/3}$ (using the Tracy-Widom-type variance estimate from~\cite{ledoux2010beta}), but for $c$ a large enough universal constant, the characterization of Corollary~\ref{cor:reverse-alexandrov} together with computations in Section~\ref{sec:exampl-oper-norm} show that $\theta_{j} (K) \geq w (K)/c$ up to $j \asymp j^* (K) \asymp n^2$.
This being said, \eqref{eq:ppv-reversal-alexandrov} leads to reversals with factors that may be close to $1$, while Corollary~\ref{cor:reverse-alexandrov} only deals with reversals up to constant factors; and it is plausible that the parameter $\jppv (K)$ may be intrinsic in some regimes in the former case.

\paragraph{Small ball probabilities and intrinsic volumes.}

By combining the quantitative characterizations of Gaussian small ball probabilities (Theorem~\ref{thm:gaussian-small-ball}) and intrinsic volumes (Theorem~\ref{thm:intrinsic-volumes}), we may deduce some correspondence between these two quantities.

Specifically, let $X \sim \gaussdist (0, I_n)$.
For any $\eps \in \ooint{0}{1}$ and symmetric convex body $K \neq \set{0}$ in $\R^n$, we let
\begin{equation}
  \label{eq:quantile-eps-sup}
  q \parens*{\eps, K}
  = Q \parens*{\log (1/\eps), \sup_{t \in K} \innerp{t}{X}}  
\end{equation}
denote the quantile of order $\eps$ of $\sup_{t \in K} \innerp{t}{X}$, which is the unique $u > 0$ such that $\gamma_n (u K^\circ) = \eps$.

\begin{corollary}
  \label{cor:intrinsic-quantiles-gaussian}
  For any $n \geq 1$, symmetric convex body $K \neq \set{0}$ in $\R^n$ and real $j > 0$, define
  \begin{equation}
    \label{eq:def-scaled-quantile}
    \wt \theta_{j} (K)
    = \inf_{\eps \in (0, 1)} \braces*{\eps^{-1/j} q \parens*{\eps, K} }
    % = \inf_{u > 0} \braces*{\gamma_n (u K^\circ)^{-1/j} u}
    = \parens*{\sup_{u > 0} \frac{\gamma_n (u K^\circ)}{u^j} }^{-1/j}  
    \, .
  \end{equation}
  Then for some universal constant $c \geq 1$, for every integer $j \geq 1$, 
  \begin{equation}
    \label{eq:intrinsic-quantiles-gaussian}
    c^{-1} \wt \theta_{c j} (K)
    \leq \theta_j (K)
    \leq c \, \wt \theta_{j/c} (K)
    \, .
  \end{equation}
\end{corollary}

We refer to Section~\ref{sec:proof-intr-volum-quant} for the proof of Corollary~\ref{cor:intrinsic-quantiles-gaussian}.

In fact, if $K$ is a symmetric convex body and $j \in \iset{n}$, then $\theta_j (K) = 0$ if and only if $\wt \theta_j (K) = 0$, if and only if the dimension of the linear span of $K$ is strictly less than $j$ (this follows, \eg, from~\eqref{eq:intrinsic-tsirelson} for $\theta_j (K)$ and from Example~\ref{ex:example-quantiles-intrinsic-n} for $\wt \theta_j (K)$).
In addition, as discussed in Section~\ref{sec:extr-cases-theor}, the equivalence $\theta_j (K) \asymp \wt \theta_j (K)$ holds at both extremes $j = 1$ (where it follows from a result of Lata{\l}a and Oleszkiewicz~\cite{latala1999gaussian}) and $j = n$ (where it follows from the combination of the Blaschke-Santal{\'o} inequality~\cite{blaschke1917uber,santalo1949invariante} and its converse from Bourgain and Milman~\cite{bourgain1987new}).
Together, these facts suggest that the equivalence of Corollary~\ref{cor:intrinsic-quantiles-gaussian} may hold in the more precise sense that $\theta_j (K) \asymp \wt \theta_j (K)$ for each $j = 1, \dots, n$, with exact matching indices.
Theorem~\ref{thm:intrinsic-quantile-exact} below (proved in Section~\ref{sec:proof-intrinsic-quantiles-exact}) shows that this is indeed the case.

\begin{theorem}
  \label{thm:intrinsic-quantile-exact}
  There exists a universal constant $c \geq 1$ such that, for every $n \geq 1$, symmetric convex body $K \subset \R^n$ and $1 \leq j \leq n$,
  % denoting by $\theta_j (K)$ the intrinsic width~\eqref{eq:def-intrinsic-width} and $\wt \theta_j (K)$ the rescaled quantile~\eqref{eq:def-scaled-quantile},
  \begin{equation}
    \label{eq:intrinsic-quantile-exact}
    c^{-1} \wt \theta_j (K)
    \leq \theta_j (K)
    \leq c \, \wt \theta_j (K)
    \, .
  \end{equation}
\end{theorem}

Theorem~\ref{thm:intrinsic-quantile-exact} may be adapted to non-symmetric convex bodies $K \subset \R^n$, by using that $\theta_j (K) \asymp \theta_j (K-K)$ (Remark~\ref{rem:intrinsic-non-sym}).
This implies that for any $n \geq 1$, convex body $K \subset \R^n$ and $j \geq 1$,
\begin{equation}
  \label{eq:intrinsic-quantile-non-sym}
  j V_j (K)^{1/j}
  \asymp \inf_{u > 0} \braces*{u \, \gamma_n \parens[\big]{u (K- K)^\circ}^{-1/j}}
  \, .
\end{equation}

Inequalities involving negative moments and intrinsic volumes were studied by Paouris and Pivovarov~\cite{paouris2013small}, in the context of the analysis of small ball probabilities for the volume of random convex sets.
The inequalities~\eqref{eq:intrinsic-quantile-exact} show that working with the ``weak'' version of negative moments given by~\eqref{eq:def-scaled-quantile} yields a two-sided equivalence.

The correspondence~\eqref{eq:intrinsic-quantile-non-sym}
is connected to the Blaschke-Santal{\'o} and reverse Blaschke-Santal{\'o} inequalities, to which it is equivalent for $j=n$ (Example~\ref{ex:example-quantiles-intrinsic-n}).
In~\eqref{eq:intrinsic-quantile-non-sym}, the Lebesgue measure is replaced by the Gaussian measure, while the volume is replaced by general intrinsic volumes; in particular, both terms of the equivalence are independent of the dimension $n$ of the ambient space in which $K$ is embedded.
In fact, both the Blaschke-Santal{\'o} inequality (in its functional form, through an argument of Alonso-Gutiérrez, Hernández Cifre, and Yepes Nicolás~\cite{alonso2021further}) and the Bourgain-Milman reverse inequality, are used in the proof of Theorem~\ref{thm:intrinsic-quantile-exact}.

\subsection{Low-dimensional projections of convex sets}
\label{sec:low-dim-sections}

Kubota's and Tsirel'son's formulas~\eqref{eq:intrinsic-tsirelson} relate intrinsic volumes to the average volumes of random projections and Gaussian linear images of convex bodies, respectively.
In particular, the representation~\eqref{eq:intr-width-interpret} identifies the intrinsic width $\theta_j(K)$ with the radius of a $j$-dimensional Euclidean ball whose volume equals the expected volume of the Gaussian image $G_j K$.
Alongside this volumetric relationship, one may also compare random projections of convex bodies to Euclidean balls in the stricter sense of set inclusion.
This naturally leads to the study of their inradii and circumradii, namely the radii of the largest Euclidean balls they contain and the smallest Euclidean balls containing them.

The study of the radius and inradius of random sections or projections of convex bodies has a long history in the local theory of Banach spaces.
A fundamental milestone was Dvoretzky's theorem on almost Euclidean sections of convex bodies~\cite{dvoretzky1959theorem} and its subsequent probabilistic proof by Milman~\cite{milman1971new} based on random subspaces.
We refer to~\cite[Chapter~5]{artstein2015asymptotic}, \cite[Chapter~11]{vershynin2018high} and~\cite[Chapter~1]{pisier1986probabilistic} for further background.

To make the connection with previous sections explicit, we first consider the radii and inradii of Gaussian linear images of convex bodies.
We then translate the resulting estimates to random sections and compare them with existing results in the literature.

In the former setting, the behavior of the radius is sharply characterized by existing results, whereas that of the inradius remains less well understood.
This parallels the situation for suprema of Gaussian processes (Sections~\ref{sec:gauss-proc-prob} and~\ref{sec:gaussian-small-ball}),
where sharp estimates for the upper tail were already available, while the understanding of the lower tail was incomplete.

\paragraph{Radius of projections.}

As a starting point and for comparison purposes, we first recall the well-understood behavior of the radius (or diameter, for non-symmetric sets) of Gaussian projections.
This description is due to V.~Milman; we refer to~\cite[Proposition~5.7.1 p.~188]{artstein2015asymptotic} and~\cite[Theorem~7.7.1 p.~171]{vershynin2018high} for textbook expositions of this result.
The version given below is similar to the presentation from~\cite[Exercises~7.7.3 \&~7.7.4 pp.~172--173]{vershynin2018high}; for completeness, we give a proof in Section~\ref{sec:proof-diam-proj}.

Recall also for $j \geq 1$ the definition~\eqref{eq:def-gaussian-matrix} of the Gaussian linear map $G_j : \R^n \to \R^j$.

\begin{proposition}
  \label{prop:diameter-gaussian-proj}
  There exists a universal constant $c \geq 1$ such that, for every $j, n \geq 1$ and compact subset $T$ of $\R^n$,
  \begin{equation}
    \label{eq:diam-gaussian-proj}
    \Probab*{c^{-1} \braces[\big]{w (T) + \diam (T) \sqrt{j}}
      \leq \diam (G_j T)
      \leq c \braces[\big]{w (T) + \diam (T) \sqrt{j}}
    }
    \geq 1 - 2 e^{-j/c}
    \, .
  \end{equation}
\end{proposition}

Equivalently (up to changing the constant $c$), for any symmetric convex body $K$ in $\R^n$,
\begin{equation}
  \label{eq:radius-gaussian-proj}
  \Probab*{c^{-1} \braces[\big]{w (K) + \rad (K) \sqrt{j}}
      \leq \rad (G_j K)
      \leq c \braces[\big]{w (K) + \rad (K) \sqrt{j}}
    }
    \geq 1 - 2 e^{-j/c}
    \, .
\end{equation}
There are thus two regimes in the behavior of the radius of $G_j K$.

Letting $\jdm(K) = (w (K)/\rad (K))^2$ denote the Dvoretzky-Milman dimension of $K$, as long as $j \leq c_0 \jdm (K)$ then $\rad (G_j K) \asymp w (K)$ with high probability.
In fact, in this regime (and for $c_0$ small enough), it follows from the Dvoretzky-Milman theorem (see \eg~\cite[Theorem~11.3.3 p.~260]{vershynin2018high}) that $G_j K$ also contains a centered Euclidean ball of radius of order $w (K)$.
Therefore, up to the dimension $c_0 \jdm (K)$, the Gaussian projection $G_j K$ resembles a Euclidean ball of radius $w (K)$, in the sense that $c_1^{-1} w (K) B_2^j \subset G_j K \subset c_1 w (K) B_2^j$ for some $c_1 \geq 1$ depending on $c_0$.

On the other hand, for $j \geq c_0 \jdm (K)$ one has $\rad (G_j K) \asymp \rad (K) \sqrt{j}$ with high probability.
Thus, for $j \gg \jdm (K)$ we have $\rad (G_j K) \gg w(K)$, and thus the dimension $\jdm (K)$ is the critical dimension (up to constant factors) for the bound $\rad (G_j K) \leq c_1 w (K)$, that is, for the inclusion $G_j K \subset c_1 w (K) B_2^j$, a phenomenon put forward by Milman and Schechtman~\cite{milman1997global}.

\paragraph{Inradius of projections.}

We now turn to the reverse inclusion, namely the analysis of Euclidean balls contained in $G_j K$.
Below, if $K \subset \R^n$ is a nonempty closed symmetric convex set, we define its \emph{inradius} $\inrad (K) \in \R^+ \cup \set{+ \infty}$ as the radius of the largest ball contained in $K$:
\begin{equation}
  \label{eq:def-inradius}
  \inrad (K)
  = \sup \set*{r \geq 0 : r B_2^n \subset K}
  \, .
\end{equation}

In contrast with the radius, the behavior of the inradius $\inrad (G_j K)$ has not been sharply characterized in the literature.
We discuss existing estimates below.

Our main result in this section is Theorem~\ref{thm:inner-radius-gaussian-proj} below, which provides a quantitative characterization of the inradius of Gaussian projections of symmetric convex bodies.
This complements the characterization~\eqref{eq:radius-gaussian-proj} of the radius of Gaussian projections.

\begin{theorem}
  \label{thm:inner-radius-gaussian-proj}
  There exists a universal constant $c \geq 1$ such that the following holds.
  For any $n \geq 1$, symmetric convex body $K \neq \set{0}$ in $\R^n$, scale $u \in \ooint{0}{w(K)}$ and integer $j \geq 1$:
  \begin{enumerate}
  \item If $j \leq c^{-1} u^2/r_* (K, u)^2$, then
    \begin{equation}
      \label{eq:inrad-gaussian-lower}
      \Probab*{\inrad \parens*{G_j K} \geq \frac{u}{c}}
      \geq 1 - c e^{- j/c}
      \, ;
    \end{equation}
  \item If $j \geq c u^2/r_* (K, u)^2$, then
    \begin{equation}
      \label{eq:inrad-gaussian-upper}
      \Probab*{\inrad \parens*{G_j K} \leq c u}
      \geq 1 - c e^{-j/c}
      \, .
    \end{equation}
  \end{enumerate}
\end{theorem}

Above, we recall that $r_* (K, u)$ denotes the effective radius~\eqref{eq:def-radius-widths} associated with local Gaussian widths.
Hence, the behavior of the inradius of Gaussian projections a symmetric convex body is governed by its local widths.
We note in passing that, in analogy with the situation for upper and lower deviations for Gaussian processes (Sections~\ref{sec:gauss-proc-prob} and~\ref{sec:gaussian-small-ball}), while convexity plays no role in the description of the diameter of projections (Proposition~\ref{prop:diameter-gaussian-proj}), it is of course essential when the inradius is considered (Theorem~\ref{thm:inner-radius-gaussian-proj}), as the inradius may be $0$ without convexity.

The proof of Theorem~\ref{thm:inner-radius-gaussian-proj} is provided in Sections~\ref{sec:proof-lower-bound-inradius}--\ref{sec:proof-upper-inrad}.
The lower bound~\eqref{eq:inrad-gaussian-lower} on the inradius follows readily from well-established tools to the lower inclusion in the Dvoretzky-Milman theorem (in particular from the Gaussian approach, as developed by Maurey-Pisier~\cite{pisier1986probabilistic} and Gordon~\cite{gordon1985some}), together with a simple inclusion argument.
The idea of the proof of the upper bound~\eqref{eq:inrad-gaussian-upper} on the inradius is the following: if the inradius of $G_j K$ were large, then this must be due to the contribution of points in $K$ with small norm, as points with large enough norm are projected to points with large norm as well.
In other words, $G_j (K \cap r B_2^n)$ must also contain a large Euclidean ball in $\R^j$.
However, with high probability, the Gaussian width of a Gaussian projection of a convex set cannot exceed the Gaussian width of the original set by significantly more than a $\sqrt{j}$ factor, which yields a bound on the width (and thus inradius) of $G_j (K \cap r B_2^n)$.
We refer to Section~\ref{sec:proof-upper-inrad} for details.

In the special case where $u$ is of order $w (K)$, Theorem~\ref{thm:inner-radius-gaussian-proj} characterizes the critical dimension for the lower inclusion in the Dvoretzky-Milman theorem, up to constant factors for the radius and dimension.
Specifically, the condition $G_j K \supset c_0^{-1} w(K) B_2^j$ holds with high probability up to a dimension $j \asymp w (K)^2/r_* (K, w(K)/c_1)^2$.
Since $r_* (K, w (K)/c_1) \leq \rad (K)$, the critical dimension for the lower inclusion in (the variant in terms of Gaussian projections) of the Dvoretzky-Milman theorem exceeds that for the upper inclusion.
This latter fact was already shown by Klartag and Vershynin~\cite{klartag2007small}, which we discuss in more detail below.

It is also instructive to compare the estimates of Theorem~\ref{thm:inner-radius-gaussian-proj} for the inradius of $G_j K$ with those of Theorem~\ref{thm:intrinsic-volumes} for the intrinsic volumes, that is, for the expected volume of $G_j K$.
There is an obvious relation between these two quantities: if $\inrad (G_j K) \geq u$, then $\vol_j (G_j K) \geq \vol_j (u B_2^j)$.
Thus, if $\inrad (G_j K) \geq u$ with probability (say) at least $1/2$, then $\theta_j (K) \geq 2^{-1/j} u$.
Conversely, by comparing Corollary~\ref{cor:reverse-alexandrov} with Theorem~\ref{thm:inner-radius-gaussian-proj}, we see that the condition that $\theta_j (K) \gtrsim w (K)$ holds up to the same projection dimension $j$ (up to absolute constant factors in $j$ and inequality defining the condition) as the \emph{a priori} stronger condition that $\inrad (G_j K) \gtrsim w (K)$ with high probability.
Thus, in this regime, the volume of the Gaussian projection $G_j K$ is explained by the presence of a large Euclidean ball inside this set.
However, this phenomenon does not persist in general for higher-dimensional projections (that is, for $u \ll w(K)$): indeed, the inradius of the Gaussian projection is solely driven by local Gaussian widths, while the quantitative behavior of intrinsic volumes also depends on covering numbers.

\paragraph{Random projections and sections of convex bodies.}

Via standard arguments, results on Gaussian linear images of symmetric convex bodies can be translated in terms of random projections or random sections of symmetric closed convex sets.

Specifically, let $1 \leq j \leq n/2$ and let $F_j$ denote the orthogonal of the kernel of $G_j$, which is the image of the adjoint $G_j^* : \R^j \to \R^n$.
Then $F_j$ follows the uniform (invariant under the action of orthogonal matrices) distribution over $j$-dimensional subspaces of $\R^n$.
In addition, for some universal constant $c \geq 1$, with probability at least $1 - 2 e^{-n/c}$ all nonzero singular values of $G_j$ are between $c^{-1} \sqrt{n}$ and $c \sqrt{n}$
(we refer to Section~\ref{sec:proof-radius-section} for more detail on the claims of this paragraph).
Under this event, we thus have
\begin{equation}
  \label{eq:corres-gaussian-orth-proj}
  \rad (G_j K)
  \asymp \sqrt{n} \rad \parens*{P_{F_j} K}
  \quad \text{ and } \quad
  \inrad (G_j K)
  \asymp \sqrt{n} \inrad_{F_j} \parens*{P_{F_j} K}
  \, ,
\end{equation}
where for a subspace $F$ of $\R^n$, we let $P_F$ denote the orthogonal projection on $F$, and for a symmetric convex subset $K'$ of $F$, we let $\inrad_F (K') = \sup \set{r \geq 0 : (r B_2^n) \cap F \subset K'}$.

In addition, for a subspace $F$ of $\R^n$, the elementary duality relation between sections and projections
\begin{equation*}
  \parens*{P_{F} K}^{\circ F}
  = K^\circ \cap F
\end{equation*}
(where for $K' \subset F$ a symmetric convex body, we let $(K')^{\circ F} = (K')^\circ \cap F$ denote the polar of $K'$ in $F$) shows that
\begin{equation}
  \label{eq:corres-sections-proj}
  \rad (K^\circ \cap F)
  = \parens*{\inrad_{F} \parens*{P_F K}}^{-1}
  \quad \text{ and } \quad
  \inrad_F (K^\circ \cap F)
  = \parens*{\rad (P_F K)}^{-1}
  \, .
\end{equation}

Using the above reductions (see Section~\ref{sec:proof-radius-section}), Theorem~\ref{thm:inner-radius-gaussian-proj} yields the following description of the diameter of random low-dimensional central sections of symmetric convex sets.

\begin{corollary}
  \label{cor:diameter-section}
  There exists a universal constant $c \geq 1$ such that the following holds.
  For any $n \geq 1$, any closed symmetric convex set $L \subset \R^n$ with nonempty interior and such that $L \neq \R^n$, any $u \in \ooint{0}{w (L^\circ)}$ and any integer $j$ with $1 \leq j \leq n/2$, letting $F_j$ be a uniformly distributed random subspace of dimension $j$ of $\R^n$:
  \begin{enumerate}
  \item If $j \leq c^{-1} u^2/r_* (L^\circ, u)^2$, then
    \begin{equation}
      \label{eq:diam-sec-upper}
      \Probab*{\diam (L \cap F_j) \leq \frac{c \sqrt{n}}{u}}
      \geq 1 - c e^{-j/c}
      \, ;
    \end{equation}
  \item If $j \geq c u^2 / r_* (L^\circ, u)^2$, then
    \begin{equation}
      \label{eq:diam-sec-lower}
      \Probab*{\diam (L \cap F_j) \geq \frac{\sqrt{n}}{c u}}
      \geq 1 - c e^{-j/c}
      \, .
    \end{equation}
  \end{enumerate}
\end{corollary}

\paragraph{Comparison with existing estimates.}

The literature on the radius of random sections or inradius of random projections of convex bodies is considerable, including sharp estimates for special convex bodies, or special affine positions of general convex bodies.
Some works also seek almost isometric estimates, with multiplicative factors arbitrarily close to $1$.
To keep the discussion focused, we restrict attention to estimates up to universal constant factors that hold for general symmetric convex bodies in arbitrary affine positions.

Below, we let $L \subset \R^n$ be a symmetric closed convex set with non-empty interior (distinct from $\R^n$), so that $K = L^\circ$ is a symmetric convex body.
% We also let $F_j$ denote a uniformly random $j$-dimensional subspace of $\R^n$.
We also let $F_j$ denote a random subspace of $\R^n$, following the uniform distribution over $j$-dimensional subspaces of $\R^n$.

As mentioned above, the key result in this spirit is the Dvoretzky-Milman theorem (see \eg \cite[Theorem~5.2.10 p.~169]{artstein2015asymptotic}; or \cite[Theorem~11.3.3 pp.~260--261]{vershynin2018high}).
It implies that for some absolute constant $c \geq 1$, if $j \leq c^{-1} \jdm (L^\circ)$ then $\rad (L \cap F_j)$ is of order $\sqrt{n}/w (L^\circ)$ with high probability.

Beyond this regime, a rather different approach yields a precise description of the behavior of $\rad (L \cap F_j)$ in the ``proportional regime'', where both the dimension $j$ and codimension $n-j$ are of order $n$ (see~\cite[Chapter~7]{artstein2015asymptotic}).
Specifically, Milman's $M^*$ estimate~\cite{milman1985random}, together with refinements from Pajor and Tomczak-Jaegermann~\cite{pajor1986subspaces} and Gordon~\cite{gordon1988escape}
and a localization argument, yields the following result~\cite[Lemma~7.5.1 p.~244]{artstein2015asymptotic}: for some absolute constant $c \geq 1$, with probability at least $1 - 4 \exp (- (n-j)/c)$ we have
\begin{equation}
  \label{eq:bound-m-star-estimate}
  \rad (L \cap F_j) \leq \widecheck{r}_j
  \, ,
\end{equation}
where $\widecheck r_j$ is the solution $r > 0$ of the equation
\begin{equation}
  \label{eq:fixed-point-m-star-estimate}
  w \parens*{L \cap r B_2^n}
  = c^{-1}  r \sqrt{n-j}
  \, .
\end{equation}
In addition, this upper bound on the radius can be complemented by a lower bound from Giannopoulos, Milman, and Tsolomitis~\cite{giannopoulos2005asymptotic}, see~\cite[Theorem~7.5.2 p.~245]{artstein2015asymptotic}, which holds for any subspace $F_j$ of dimension $j$.
This lower bound is comparable to the upper bound above in the proportional regime where $\min \set{j, n-j} \geq n/c$ for some absolute constant $c > 1$, but becomes asymptotically smaller in both the low-dimensional regime $j \ll n$ and the low-codimensional regime $n-j \ll n$.

In contrast, the characterization of Corollary~\ref{cor:diameter-section} is two-sided, but holds up to constant factors in $j$; it is thus only suitable when $j \leq n/c$ for some universal constant $c > 1$.
In order to make the comparison with the localized $M^*$ estimate clearer, we may rewrite the quantities involved in Corollary~\ref{cor:diameter-section}.
For $1 \leq j \leq n$ and $\kappa > 0$, let $u \in (0, w (L^\circ))$ be such that $j = \kappa^{-1} u^2/r_* (L^\circ, u)^2$ (see Lemma~\ref{lem:loc-width-bij}).
Then, letting $r = r_* (L^\circ, u)$, we have $u = r \sqrt{\kappa j}$, and from the definition~\eqref{eq:def-radius-widths} together with Lemma~\ref{lem:loc-width-bij}, $w (L^\circ \cap r B_2^n) = u = r \sqrt{\kappa j}$.
Thus, Corollary~\ref{cor:diameter-section} shows that for some absolute constant $c > 1$, with probability at least $1 - c e^{-j/c}$ we have
\begin{equation}
  \label{eq:sharp-section-radius}
  \frac{\sqrt{n/j}}{c \, \wh r_{j/c}}
  \leq \rad (L \cap F_j)
  \leq \frac{c \sqrt{n/j}}{\wh r_{c j}}
  \, ,
\end{equation}
where we let $\wh r_{j'}$ denote the solution $r$ of the equation
\begin{equation}
  \label{eq:fixed-point-section}
  w (L^\circ \cap r B_2^n) = r \sqrt{j'}
  \, .
\end{equation}

It should be noted that, while the localized $M^*$ estimate in~(\ref{eq:bound-m-star-estimate}--\ref{eq:fixed-point-m-star-estimate}) and the two-sided estimate in~(\ref{eq:sharp-section-radius}--\ref{eq:fixed-point-section}) derived from Corollary~\ref{cor:diameter-section} (which can be viewed, in complement with Proposition~\ref{prop:diameter-gaussian-proj}, as an extension of the Dvoretzky-Milman estimate beyond the dimension $\jdm (L^\circ)$) both involve fixed points associated with local widths, they involve different bodies ($L$ versus $L^\circ$) and dimensional factors.
They can thus be loosely viewed as dual to each other.
In the regime where $j \leq n/c$ (for some large enough absolute constant $c > 1$) considered in this work, Corollary~\ref{cor:diameter-section} yields a sharp result; in contrast, the localized $M^*$ estimate is more appropriate in the low-codimensional regime $n-j \leq n/c$.

Recently, Guédon, Litvak, Tatarko and Vritsiou~\cite{guedon2026random} obtained further estimates for the radii of random sections, particularly in the low-codimensional regime, which complement and sometimes improve the $M^*$-estimate.
Their upper bound involves radii and Gaussian widths of projections of $L$ onto suitably chosen subspaces, and is in general incomparable with the $M^*$-estimate.
We also refer to~\cite{guedon2026random} for a broader discussion of the related literature.

\paragraph{Connection with small ball estimates.}

A one-sided improvement of the Dvoretzky-Milman theorem was established by Klartag and Vershynin~\cite{klartag2007small}, who put forward a connection with small ball probabilities.
Specifically, for a symmetric convex body $K \subset \R^n$, define
\begin{equation}
  \label{eq:klartag-vershynin}
  \jkv (K)
  = \min \braces*{n, - \log \Probab*{\sup_{t \in K} \innerp{t}{X} \leq c_0^{-1} w (K)}}
  \, , \quad
  X \sim \gaussdist (0, I_n)
  \, ,
\end{equation}
where $c_0 > 1$ is a suitable universal constant.
Klartag and Vershynin showed that, for some universal constant $c_1>1$, whenever $1 \leq j \leq \jkv(L^\circ)/c_1$, one has
$\rad(L\cap F_j)\leq c_1\sqrt{n}/w(L^\circ)$ with high probability.
Moreover, the Gaussian concentration inequality (Theorem~\ref{thm:gaussian-concentration}) implies that $\jkv (K) \geq c_2 \jdm (K)$ for some universal constant $c_2 > 0$.
Thus, the upper bound on the radius provided by the Dvoretzky-Milman theorem extends to higher dimensions than the corresponding lower bound on the inradius.

The approach of Klartag and Vershynin~\cite{klartag2007small} may be viewed as a one-sided analogue of Milman's probabilistic proof~\cite{milman1971new} of Dvoretzky's theorem, which derives two-sided estimates for random sections from concentration of measure.
Their proof relies on a Gaussian small ball inequality of Lata{\l}a and Oleszkiewicz~\cite{latala2005small}, itself derived from the B-theorem of Cordero-Erausquin, Fradelizi and Maurey~\cite{cordero2004b}, a deep concavity property of the Gaussian measure.

Finally, Klartag and Vershynin observed that the connection between small ball probabilities (or equivalently, negative moments of Gaussian processes) and the radius of random sections fails to hold beyond the critical dimension $\jkv (L^\circ)$, by noting that it may fail for ``flat'' bodies beyond a certain dimension.

The characterizations of small ball probabilities and of the radius of low-dimensional sections allow one to recover these facts.
Specifically, the connection between small ball probabilities and the critical dimension for the diameter of random sections follows from the fact (Corollaries~\ref{cor:small-ball-constant-simple} and~\ref{cor:diameter-section}) that both are governed by the geometric parameter $w (K)^2/r_* (K, w (K)/c)^2$ for $K = L^\circ$.
In addition, the discrepancy between small ball probabilities and the radius of random sections beyond the critical dimension can be explained quantitatively by the dependence of the former on both covering numbers and local widths (Theorem~\ref{thm:gaussian-small-ball}), while the latter only depend on local widths.
For sufficiently elongated convex bodies, covering number may dominate local widths, which leads to a discrepancy.

\subsection{Related work}
\label{sec:related-work}

We briefly discuss further related work, beyond that covered in the preceding sections.

\paragraph{Small ball probabilities.}

There is an extensive literature on small ball probabilities; for a broader account, we refer to~\cite[Chapter~7]{ledoux1996isoperimetry}, the surveys~\cite{li2001gaussian,lifshits1999asymptotic} and the bibliography~\cite{lifshits2023biblio}.
Beyond general estimates as in Section~\ref{sec:gaussian-small-ball},
part of the literature is devoted to the small deviation behavior of specific Gaussian processes;
a landmark result in this direction is Talagrand's sharp characterization of small ball probabilities for the two-parameter Brownian sheet~\cite{talagrand1994small}.

\paragraph{Dimension-dependent estimates.}

Paouris, Tikhomirov and Valettas~\cite{paouris2022hypercontractivity} developed a general approach to Gaussian small ball probabilities and lower deviations of norms, exploiting hypercontractivity and convexity properties of the Gaussian measure.

In particular, they showed that, after a suitable linear change of coordinates, every $n$-dimensional norm admits small ball estimates with probabilities bounded by $\exp(-c n^{1/4-\varepsilon})$, for any $\varepsilon>0$ and sufficiently small relative thresholds (depending on $\varepsilon$).
For suitably positioned unconditional norms, the exponent can be improved to $n^{1-\varepsilon}$, with the $\ell^\infty$ norm providing an extremal example.
These estimates also yield one-sided Dvoretzky-type theorems for random sections and reverse Alexandrov inequalities for intrinsic volumes, valid over the above polynomial ranges of dimensions.

In contrast with the results considered here, these estimates rely on a suitable choice of affine position and are expressed in terms of the ambient dimension $n$, rather than geometric parameters of convex bodies in arbitrary positions.

\paragraph{Intrinsic volumes, the Wills functional, and Gaussian processes.}

Further connections between the Wills functional and Gaussian processes were identified by Vitale~\cite{vitale1996wills,vitale2001intrinsic}, beyond those discussed in Section~\ref{sec:intr-volum-wills}.
Specifically, Vitale~\cite{vitale1996wills} showed that McMullen's inequality~\eqref{eq:wills-width-mcmullen} yields an alternative proof of the classical upper-tail deviation bound for suprema of Gaussian processes, usually obtained through Gaussian concentration.
% In particular, Vitale~\cite{vitale1996wills} recovered the classical Gaussian upper-tail deviation bound for suprema of Gaussian processes from McMullen's inequality
In~\cite{vitale2001intrinsic}, he further showed that, for an infinite-dimensional $K \subset \ell^2$ with $w (K) < + \infty$, the quantity $\lambda^{-1} \log W (\lambda K)$ converges as $\lambda \to + \infty$ to the maximal discontinuity of the Gaussian process indexed by $K$.

In a different direction, Chevet~\cite{chevet1976processus}, extending earlier work of Fernique~\cite{fernique1976corps} in the low-rank case, established connections between the second moment of the supremum of a Gaussian process, intrinsic volumes, and the expected squared Euclidean norm of a maximizing index of the process $(\innerp{t}{X})_{t \in K}$.

\subsection{Organization of the paper}
\label{sec:organization-paper}

The rest of this paper is organized as follows.
Section~\ref{sec:proof-small-ball} establishes the characterizations of small ball probabilities (Theorem~\ref{thm:talagrand-maj-meas} and Corollary~\ref{cor:small-ball-constant-simple}) presented in Section~\ref{sec:gaussian-small-ball}.
Next, Section~\ref{sec:conn-with-wills} is devoted to the proofs of results of Section~\ref{sec:intr-volum-wills} concerning intrinsic volumes, the Wills functional, and their connections to small ball probabilities.
Section~\ref{sec:chain-char-quant} then provides the proofs of the results of Section~\ref{sec:chaining-comparison} on the relationship between small ball probabilities and the metric structure of the index set.
The results of Section~\ref{sec:low-dim-sections} on radii and inradii of random projections and sections of symmetric convex sets are proved in Section~\ref{sec:low-dimens-sect}.
Finally, Section~\ref{sec:basic-prop-supr} contains proofs of basic properties of Gaussian processes referenced elsewhere in the paper, while Section~\ref{sec:examples} presents examples illustrating claims made in the present section.

\subsection{Notation}
\label{sec:notation}

Throughout, we let $\innerp{\cdot}{\cdot}$ denote the standard scalar product on $\R^n$, defined by $\innerp{t}{x} = \sum_{i=1}^n t_i x_i$ for $t = (t_{i})_{1 \leq i \leq n} \in \R^n$ and $x = (x_i)_{1 \leq i \leq n} \in \R^n$, and let $\norm{\cdot} = \norm{\cdot}_2$ denote denote the associated Euclidean norm, defined by $\norm{t} = \sqrt{\innerp{t}{t}}$.
We let $B_2^n = \set{t \in \R^n : \norm{t} \leq 1}$ denote the unit Euclidean ball in $\R^n$; when the dimension is the ambient dimension $n$, we also let $B_2 = B_2^n$.

By \emph{convex body} in $\R^n$ we mean a nonempty convex and compact subset of $\R^n$, without assuming that it has a nonempty interior.
We will generally use the letter $T$ for general (not necessarily convex) sets, and $K$ for convex bodies; we will also typically use the letter $L$ for closed convex sets with non-empty interior, which play a dual role to sets denoted by $K$.
We say that a subset $T \subset \R^n$ is \emph{symmetric} if it is origin symmetric, \ie, if $T = -T = \set*{-t : t \in T}$.
The diameter $\diam (T) = \sup_{t, t' \in T} \norm{t - t'} \subset \R^n$ of a set $T$ is defined with respect to the Euclidean distance.
If $T\subset \R^n$ is symmetric, its radius is defined by $\rad (T) = \sup \set{\norm{t} : t \in T} = \diam (T)/2$.
Thus, denoting by $\conv (T)$ the convex hull of a convex set $T$, we have $\diam (T) = \rad (T-T) = \rad (\conv (T-T))$.
The inradius $\inrad$ is defined in~\eqref{eq:def-inradius}.

We also let $\gamma_n = \gaussdist (0, I_n)$ (where $I_n$ is the $n\times n$ identity matrix) denote the standard Gaussian measure on $\R^n$.
We will generally denote by $X$ a standard Gaussian random vector, and for $T \subset \R^n$, we recall the definition~\eqref{eq:def-gaussian-width} of its Gaussian width $w (T)$.

We use $c_0, c_1, \dots$ to denote universal constants, whose value may change from statement to statement.
Given two non-negative quantities $a, b$, the notation $a \lesssim b$ means that there exists a universal constant $c > 0$ such that $a \leq c b$; while $a \asymp b$ means that $a \lesssim b$ and $b \lesssim a$.

\subsection{Acknowledgements}
\label{sec:acknowledgements}

Part of the research presented in this work was conducted while the author was participating in the ``Mathematics of Data'' program at the Institute for Mathematical Sciences, National University of Singapore, in January 2024.
The author is grateful to the IMS and the program organizers for their hospitality and for providing a stimulating research environment.

\paragraph{Statement on AI use.}

All results presented in this note, with the exception of Theorem~\ref{thm:intrinsic-quantile-exact}, were obtained prior to Spring 2026 and without the use of AI.

Theorem~\ref{thm:intrinsic-quantile-exact} (motivated by the prior Corollary~\ref{cor:intrinsic-quantiles-gaussian}) was obtained in the final stages of the preparation of this manuscript using GPT--6.
The proof was then simplified and re-organized by the author, who takes full responsibility for its correctness.
Finally, AI tools were used for language editing purposes in the abstract and parts of the introduction.

\section{Two-sided characterization of Gaussian small ball probabilities}
\label{sec:proof-small-ball}

The purpose of this section is to establish the quantitative characterization of small ball probabilities from Theorem~\ref{thm:gaussian-small-ball}, as well as Corollary~\ref{cor:small-ball-constant-regime}.

\subsection{Preliminaries}
\label{sec:preliminaries}

We first collect some tools that will be used throughout the proof.

\begin{definition}[Covering numbers]
  \label{def:covering-numbers}
  If $T \subset \R^n$ is bounded and $L \subset \R^n$ is symmetric convex with non-empty interior, we define $N (T, L)$ as the smallest number of translates of $L$ centered in $T$ that cover $T$:
  \begin{equation}
    \label{eq:def-covering-general}
    N (T, L)
    = \inf \braces*{N \geq 1 : \exists t_1, \dots, t_N \in \R^n, \ T \subset \bigcup_{i=1}^N \parens*{t_i + L}}
    \, .
  \end{equation}
  In addition, for any $r > 0$, we let $N (T, r) = N (T, r B_2)$, the \emph{covering number} of $T$ at scale $r$.
\end{definition}

If $T$ is compact convex, in the definition of $N (T, r)$ we may relax the assumption that the centers $t_1, \dots, t_N \in \R^n$ belong to $T$ (since projecting them on $T$ only reduces their distances to elements of $T$).
Therefore, if $T_1, T_2$ are compact convex with $T_1 \subset T_2$, then $N (T_1, r) \leq N (T_2, r)$ for every $r > 0$.
% In addition, a standard compactness argument shows that the map $r \mapsto N (T, r)$ is right-continuous if $T$ is convex compact.

We will use the fundamental duality relation between covering numbers~\cite[Theorem~1]{artstein2004duality}:

\begin{theorem}[\cite{artstein2004duality}, Theorem~1]
  \label{thm:duality-covering}
  There exists a universal constant $c \geq 1$ such that the following holds.
  For every $n \geq 1$ and every symmetric convex body $K \subset \R^n$,
  \begin{equation}
    \label{eq:duality-covering}
    N \parens*{K, c}^{1/c}
    \leq N \parens*{B_2, K^\circ}
    \leq N \parens*{K, 1/c}^c
    \, .
  \end{equation}
\end{theorem}

We will also use the following important facts about the Gaussian measure of shifted symmetric convex sets.
The lower bound below is classical and admits a direct proof (see \eg~\cite[Lemma~15.2.5]{talagrand2021upper}), while the upper bound is due to Anderson~\cite{anderson1955integral} and readily follows from the (subsequently established) Prékopa-Leindler inequality.

\begin{lemma}
  \label{lem:anderson}
  For any closed symmetric convex set $L \subset \R^n$ and any $x \in \R^n$, 
  \begin{equation}
    \label{eq:anderson}
    \gamma_n (L) e^{-\norm{x}^2/2}
    \leq \gamma_n (x + L)
    \leq \gamma_n (L)
    \, .
  \end{equation}
\end{lemma}

\subsection{Upper bound}
\label{sec:proof-small-ball-upper}
%\label{sec:proof-small-ball-upper-width}

We start by showing a first upper bound, which is obtained by Gaussian concentration applied to a suitable subset of $K$.

\begin{lemma}
  \label{lem:upper-sb-width}
  Let $K \neq \set{0}$ be a symmetric convex body in $\R^n$.
  For any $\eps > 0$ and $u \in (0, \frac{w (K)}{1+\eps})$, 
  \begin{equation}
    \label{eq:upper-sb-width}
    \gamma_n (u K^\circ)
    \leq \exp \parens*{- \frac{\eps^2 u^2}{2 r_* (K, (1+\eps) u)^2}}
    \, .
  \end{equation}
\end{lemma}

\begin{proof}[Proof of Lemma~\ref{lem:upper-sb-width}]
  Since $(1+\eps) u < w (K)$, by Lemma~\ref{lem:loc-width-bij},
  letting $r = r_* (K, (1+\eps) u) \in (0, \rad (K))$ we have
  $w (K \cap r B_2) \geq (1+\eps) u$.
  Since in addition the map $x \mapsto - \sup_{t \in K \cap r B_2} \innerp{t}{x}$ is $r$-Lipschitz, Gaussian concentration (Theorem~\ref{thm:gaussian-concentration}) gives, letting $X \sim \gaussdist (0, I_n)$,
  \begin{align*}
    \gamma_n (u K^\circ)
    &= \Probab*{\sup_{t \in K} \innerp{t}{X} \leq u}
      \leq \Probab*{\sup_{t \in K \cap r B_2} \innerp{t}{X} \leq u} \\
    &\leq \Probab*{\sup_{t \in K \cap r B_2} \innerp{t}{X} - w \parens*{K \cap r B_2} \leq - \eps u} 
    \leq \exp \parens*{- \frac{\eps^2 u^2}{2 r^2}}
      \, .
      \qedhere
  \end{align*}
\end{proof}

We now establish a second, different upper bound on small ball probabilities, involving covering numbers.

\begin{lemma}
  \label{lem:upper-sb-covering}
  For some universal constant $c \geq 1$, one has for every $n\geq 1$, any symmetric convex body $K \neq \set{0}$ in $\R^n$ and any $u > 0$,
  \begin{equation}
    \label{eq:upper-sb-covering}
    \gamma_n (u K^\circ)
    \leq \exp \parens*{- \frac{2 u^2}{\wt r (K, c u)^2}}
    \, .
  \end{equation}
\end{lemma}

\begin{proof}[Proof of Lemma~\ref{lem:upper-sb-covering}]
  The proof
  starts similarly to the approach of Kuelbs and Li~\cite{kuelbs1993metric}, but then proceeds differently by introducing the quantity $\wt r (K, c_3 u)$.

  Let $\rho > 0$ to be fixed later, and define the following packing number:
  \begin{equation}
    \label{eq:def-packing-dual}
    N
    = N_p (\rho B_2, 2 u K^\circ)
    = \sup \set*{N \geq 1 : \exists x_1, \dots, x_N \in \rho B_2, \ \forall i \neq j, \ x_j \not\in \parens*{x_i + 2 u K^\circ} }
    .
  \end{equation}
  It is a standard fact~\cite[Lemma 4.2.8]{vershynin2018high} that packing and covering numbers (Definition~\ref{def:covering-numbers}) are related as follows:
  \begin{equation}
    \label{eq:covering-packing}
    N (\rho B_2, 2 u K^\circ)
    \leq N_p (\rho B_2, 2 u K^\circ)
    \leq N (\rho B_2, u K^\circ)
    \, .
  \end{equation}
  Now, by definition~\eqref{eq:def-packing-dual} of $N$, let $x_1, \dots, x_n \in \rho B_2$ such that for every $i \neq j$, one has $x_j \not \in x_i + 2 u K^\circ$, that is (as $K^\circ$ is symmetric convex), $\parens*{x_i + u K^\circ} \cap \parens*{x_j + u K^\circ} = \varnothing$.
  We deduce that
  \begin{align*}
    1
    &\geq \gamma_n \parens*{\bigcup_{i=1}^N \parens*{x_i + u K^\circ}}
      = \sum_{i=1}^N \gamma_n \parens*{x_i + u K^\circ} 
      \geq \sum_{i=1}^N e^{-\norm{x_i}^2/2} \gamma_n \parens*{u K^\circ}
      \geq N e^{-\rho^2/2} \gamma_n (u K^\circ)
      \, ,
  \end{align*}
  where the penultimate inequality follows from Lemma~\ref{lem:anderson} and the last one from the fact that $x_i \in \rho B_2$.
  Re-arranging and using the lower bound in~\eqref{eq:covering-packing}, we deduce that
  \begin{equation}
    \label{eq:proof-prelim-cover}
    - \log \gamma_n (u K^\circ)
    \geq \log N_p \parens*{\rho B_2, 2 u K^\circ} - \frac{\rho^2}{2} 
    \geq \log N \parens*{\rho B_2, 2 u K^\circ} - \frac{\rho^2}{2}
    \, .
  \end{equation}
  Now, letting $c_2 \geq 1$ denote the universal constant of the duality relationship between covering numbers (Theorem~\ref{thm:duality-covering}), we have
  \begin{equation*}
    N \parens*{\rho B_2, 2 u K^\circ}
    = N \parens*{B_2, \parens*{\frac{\rho}{2 u} K}^\circ}
    \geq N \parens*{\frac{\rho}{2 u} K, c_2}^{1/c_2}
    = N \parens*{K, \frac{2 c_2 u}{\rho}}^{1/c_2}
    \, .
  \end{equation*}
  Therefore, inequality~\eqref{eq:proof-prelim-cover} yields
  \begin{equation}
    \label{eq:proof-sb-cover-diff}
    - \log \gamma_n (u K^\circ)
    \geq \frac{1}{c_2} \log N \parens*{K, \frac{2 c_2 u}{\rho}} - \frac{\rho^2}{2}
    \, .
  \end{equation}

  Up to this point, the proof follows the same lines as that of~\cite{kuelbs1993metric}, see also~\cite[p.~256]{ledoux1996isoperimetry}.
  Specifically, after~\eqref{eq:proof-prelim-cover}, Kuelbs and Li set $\rho = \sqrt{- \log \gamma_n (u K^\circ)}$, leading to the upper bound in Theorem~\ref{thm:kuelbs-li}.

  We proceed differently.
  Set $c = 2 c_2^{3/2} \geq 1$.
  If $\wt r \parens*{K, c u} = + \infty$, the bound~\eqref{eq:upper-sb-covering} holds trivially as $\gamma_n (u K^\circ) \leq 1$.
  We thus assume that $\wt r \parens*{K, c u} < +\infty$.
  Let $r' > \wt r \parens*{K, c u}$, so that there exists $r \leq r'$ such that $r \sqrt{\log N (K, r)} \geq 2 c_2^{3/2} u$, namely
  \begin{equation}
    \label{eq:proof-upper-cov-fixed}
    \log N (K, r)
    \geq \frac{4 c_2^3 u^2}{r^2}
    \, .
  \end{equation}
  Also, let $\rho = 2 c_2 u/r$.
  With this choice of $\rho$, inequality~\eqref{eq:proof-sb-cover-diff} together with~\eqref{eq:proof-upper-cov-fixed} gives (as $c_2 \geq 1$)
  \begin{align*}
    - \log \gamma_n (u K^\circ)
    &\geq \frac{1}{c_2} \log N \parens*{K, r} - \frac{2 c_2^2 u^2}{r^2} 
    \geq \frac{1}{c_2} \frac{4 c_2^3 u^2}{r^2} - \frac{2 c_2^2 u^2}{r^2}
      = \frac{2 c_2^2 u^2}{r^2}
      \geq \frac{2 u^2}{(r')^2}
      \, .
  \end{align*}
  Since $r' > \wt r \parens*{K, c u}$ was arbitrary,  
  inequality~\eqref{eq:upper-sb-covering} ensues.
\end{proof}

\subsection{Lower bound}
\label{sec:proof-small-ball-lower}

We now turn to a lower bound on small ball probabilities that matches the combination of the previous two upper bounds.

\begin{lemma}
  \label{lem:lower-sb-width-cover}
  There exists a universal constant $c \geq 1$ such that %the following holds.
  for any $\eps \in \ocint{0}{1/2}$ and $\kappa_0 > 0$, any integer $n \geq 1$, symmetric convex body $K \neq \set{0}$ in $\R^n$ and $u \in \ooint{0}{w (K)}$,
  one has
  \begin{equation}
    \label{eq:lower-sb-width-cover}
    \gamma_n (u K^\circ)
    \geq \eps \exp \parens*{- \frac{c u^2}{\kappa_0^2 \min \braces[\big]{\eps^2 r_* \parens*{K, (1-\eps)^2 u}^2, \wt r (K, u/\kappa_0)^2}}}
    \, .
  \end{equation}
\end{lemma}

We note that the argument used to establish Lemma~\ref{lem:lower-sb-width-cover} features some connections to the one used by Kuelbs and Li~\cite{kuelbs1993metric}, in that it involves on combination of covering and isoperimetry.
However, a critical difference is that we control the measure of certain enlarged sets using the information provided by local widths, as opposed to the Gaussian isoperimetric inequality in the approach of Kuelbs and Li.

\begin{proof}[Proof of Lemma~\ref{lem:lower-sb-width-cover}]
  Let $0 < r < \min \braces*{r_* \parens*{K, (1-\eps)^2 u}, \wt r (K, u/\kappa_0)/\eps}$,
  so that (regardless of whether $\wt r (K, u/\kappa_0)$ is finite or infinite)
  \begin{equation}
    \label{eq:proof-sb-lower-cond-r}
    w \parens*{K \cap r B_2} < (1-\eps)^2 {u}
    \quad \text{ and } \quad
    \log N \parens*{K, r'}
    < \frac{u^2}{\kappa_0^2 r'^2} \ \forall r' \leq \eps r
    \, .
  \end{equation}
  The first condition in~\eqref{eq:proof-sb-lower-cond-r} together with Markov's inequality shows that, if $X \sim \gaussdist(0, I_n)$,
  \begin{equation*}
    \Probab*{\sup_{t \in K \cap r B_2} \innerp{t}{X} \geq (1-\eps) u}
    \leq 1 - \eps
    \, ,
  \end{equation*}
  or equivalently,
  \begin{equation*}
    \gamma_n \parens*{(1-\eps) u (K \cap r B_2)^\circ}
    \geq \eps
    \, .
  \end{equation*}
  Now, since $K^\circ + r^{-1} B_2 \supset K^\circ$ and $K^\circ + r^{-1} B_2 \supset r^{-1} B_2$, by polarity between symmetric convex bodies we have $(K^\circ + r^{-1} B_2)^\circ \subset (K^\circ)^\circ \cap (r^{-1} B_2)^\circ = K \cap r B_2$ and thus $K^\circ + r^{-1} B_2 \supset (K \cap r B_2)^\circ$.
  Thus, the previous inequality implies that
  \begin{equation}
    \label{eq:proof-sb-lower-large-mass}
    \gamma_n \parens*{(1-\eps) u \bracks*{K^\circ + \frac{1}{r} B_2}}
    \geq \eps
    \, .
  \end{equation}

  Now, let $\rho, \nu > 0$ (to be chosen later) and $N = N \parens*{\rho B_2, \nu K^\circ}$, so that there exist $x_1, \dots, x_N \in \R^n$ such that
  \begin{equation*}
    \rho B_2 \subset \bigcup_{i=1}^N \parens*{x_i + \nu K^\circ}
    \, .
  \end{equation*}
  Therefore, for $\sigma > 0$ (using that $\nu K^\circ + \sigma K^\circ = (\nu+\sigma) K^{\circ}$ by convexity of $K^\circ$),
  \begin{equation*}
    \sigma K^\circ + \rho B_2
    \subset \bigcup_{i=1}^N \parens*{x_i + (\nu + \sigma) K^\circ}
    \, .
  \end{equation*}
  Hence, setting
  $\sigma = (1-\eps) u$, $\rho = (1-\eps) u/r$ and $\nu = \eps u > 0$, we obtain from Lemma~\ref{lem:anderson} that
  \begin{equation*}
    \gamma_n \parens*{(1-\eps) u \bracks*{K^\circ + \frac{1}{r} B_2}}
    \leq \gamma_n \parens*{\bigcup_{i=1}^N \parens*{x_i + u K^\circ}}
    \leq \sum_{i=1}^N \gamma_n \parens*{x_i + u K^\circ}
    \leq N \gamma_n \parens*{u K^\circ}
    \, .
  \end{equation*}
  When combined with~\eqref{eq:proof-sb-lower-large-mass}, this yields
  \begin{equation}
    \label{eq:proof-sb-lower-eps-N}
    \gamma_n (u K^\circ)
    \geq \frac{\eps}{N}
    \, .
  \end{equation}

  Finally, letting $c \geq 1$ denote a universal constant for Theorem~\ref{thm:duality-covering},
  \begin{align*}
    N
    &= N \parens*{\frac{(1-\eps) u}{r} B_2, \eps u K^\circ}
      = N \parens*{B_2, \frac{\eps r}{1-\eps} K^\circ} 
      \leq N \parens*{K, \frac{\eps r}{c (1-\eps)}}^{c}
      % \leq N \parens*{K, \frac{2 \eps r}{c}}^{c}
      \, .
  \end{align*}
  Now, using that
  $\frac{\eps r}{c (1-\eps)} \geq \frac{\eps r}{c}$ with
  $\frac{\eps r}{c} \leq \eps r$,
  % \jao{def (so one could take $c_6 = 1$, but maybe there is something to gain here)},
  we deduce from~\eqref{eq:proof-sb-lower-cond-r} that
  \begin{equation*}
    N
    \leq \exp \braces*{c \log N \parens*{K, \frac{\eps r}{c }}}
    \leq \exp \braces*{c \frac{u^2}{\kappa_0^2 (\eps r/c)^2}}
    = \exp \parens*{\frac{c^3 u^2}{\kappa_0^2 \eps^2 r^2}}
    \, .
  \end{equation*}
  Plugging this inequality into~\eqref{eq:proof-sb-lower-eps-N} gives:
  \begin{equation*}  
    \gamma_n (u K^\circ)
    \geq \eps \exp \parens*{- \frac{c^3 u^2}{\kappa_0^2 \eps^2 r^2}}
    \, .
  \end{equation*}
  Since this holds for every $r < \min \braces*{r_* \parens*{K, (1-\eps)^2 u}, \wt r (K, u/\kappa_0)/\eps}$,
  this concludes the proof.
\end{proof}

\subsection{Conclusion of the proof of Theorem~\ref{thm:gaussian-small-ball}}
\label{sec:concluding-proof}

Let $c_1 \geq 1$ be the universal constant of Lemma~\ref{lem:upper-sb-covering} and $c_2 \geq 1$ the universal constant of Lemma~\ref{lem:lower-sb-width-cover}.
By combining Lemma~\ref{lem:upper-sb-width} with $\eps = 1$ and Lemma~\ref{lem:upper-sb-covering}, for every $u \in (0, w (K)/2)$ we have
\begin{equation}
  \label{eq:proof-sb-conc-1}
  - \log \gamma_n (u K^\circ)
  \geq \frac{u^2}{2 \min \braces[\big]{r_* \parens*{K, 2 u}^2, \wt r (K, c_1 u)^2}}
  \, .
\end{equation}
On the other hand, by Lemma~\ref{lem:lower-sb-width-cover} with $\eps = 1/2$ and $\kappa_0 = 1$, for every $u \in (0, w (K))$ one has
\begin{equation}
  \label{eq:proof-sb-conc-2}
  2 \gamma_n (u K^\circ)
  \geq \exp \parens*{- \frac{4 c_2 u^2}{\min \set{r_* (K, u/4)^2, \wt r (K, u)^2}}}
  \, .
\end{equation}
Now, Lemma~\ref{lem:mini-sb} implies that if $u \leq w (K)/11$, then $\gamma_n (u K^\circ) \leq 1/4$ and thus $2 \gamma_n (u K^\circ) \leq \sqrt{\gamma_n (u K^\circ)}$, hence~\eqref{eq:proof-sb-conc-2} implies that
\begin{equation}
  \label{eq:proof-sb-conc-3}
  - \log \gamma_n (u K^\circ)
  \leq \frac{8 c_2 u^2}{\min \braces{r_* (K, u/4)^2, \wt r (K, u)^2}}
  \, .
\end{equation}

Inequalities~\eqref{eq:proof-sb-conc-1} and~\eqref{eq:proof-sb-conc-3} together imply that~\eqref{eq:gaussian-small-ball} holds with $c_0 = 11$ and $c = \max \set{c_1, 8 c_2}$.

\subsection{Proof of Corollary~\ref{cor:small-ball-constant-simple}%
  % Constant-factor deviations
}
\label{sec:const-fact-lower}

We now turn to the proof of Corollary~\ref{cor:small-ball-constant-simple}.
More precisely, we will prove the following:

\begin{corollary}
  \label{cor:small-ball-constant-regime}
  For every $\eps \in (0, 3/4)$ and $\kappa_0 \geq 2$,
  there exists a constant $\kappa_1 = \kappa_1 (\eps, \kappa_0) \geq 1$ depending on $\eps, \kappa_0$ such that the following holds.
  For every $n \geq 1$, symmetric convex body $K \neq \set{0}$ in $\R^n$ and $u \in \ooint{w (K)/\kappa_0}{w(K)/(1+\eps)}$, one has
  \begin{equation}
    \label{eq:small-ball-constant-regime}
    \frac{\eps}{2} \exp \parens*{- \frac{\kappa_1 w(K)^2}{r_* \parens*{K, (1-\eps) u}^2}}
    \leq \gamma_n (u K^\circ)
    \leq \exp \parens*{- \frac{w(K)^2}{\kappa_1 r_* \parens*{K, (1+\eps) u}^2}}
    \, .
  \end{equation}
\end{corollary}

Note that Corollary~\ref{cor:small-ball-constant-simple} follows immediately from Corollary~\ref{cor:small-ball-constant-regime} by letting $\eps = 0.01$ and $\kappa_0 = 100$; the factor $\eps/2 = 0.005$ appearing in the left-hand side of~\eqref{eq:small-ball-constant-regime} can be absorbed in the exponential, by using that
$w (K)/r_* (K, 0.99 u) \geq w (K)/\rad (K) \geq \sqrt{2/\pi}$.

The upper bound in Corollary~\ref{cor:small-ball-constant-regime} follows from Lemma~\ref{lem:upper-sb-width}, it thus remains to prove the lower bound.
The key observation is that in the regime $u \asymp w (K)$, the term depending on covering numbers is dominated by the one depending on local widths, and can therefore be omitted.

We shall use the following fundamental fact about suprema of Gaussian processes, known as \emph{Sudakov minoration}.
Specifically, Lemma 2.10.2 in~\cite{talagrand2021upper}, when combined with the correspondence between covering and packing numbers (\eg, \cite[Lemma 4.2.8]{vershynin2018high}), reads as follows:

\begin{lemma}[Sudakov minoration]
  \label{lem:sudakov}
  There exists a universal constant $c_0 \geq 1$ such that the following holds.
  For every $n \geq 1$, every subset $T \subset \R^n$ and every $r > 0$, 
  \begin{equation}
    \label{eq:sudakov}
    r \sqrt{\log N (T, r)}
    \leq c_0 w (T)
    \, .
  \end{equation}
\end{lemma}

We shall also use the following standard fact:

\begin{lemma}
  \label{lem:local-cov-max}
  Let $K \subset \R^n$ be a symmetric convex body.
  For every $c,r > 0$ and $t \in K$, 
  \begin{equation}
    \label{eq:local-cov-max}
    N \parens*{K \cap B (t, c r), r}
    \leq N \parens*{K \cap \frac{c r}{2} B_2, \frac{r}{2}}
    \, .
  \end{equation}
\end{lemma}

\begin{proof}[Proof of Lemma~\ref{lem:local-cov-max}]
  Since $K$ is symmetric one has $-t \in K$, and since $K$ is convex, we deduce that $\frac{1}{2} \bracks*{(-t) + K \cap B (t, c r)} \subset K$.
  In addition, one has $\frac{1}{2} \bracks*{(-t) + K \cap B (t, c r)} \subset \frac{c r}{2} B_2$, thus $\frac{1}{2} \bracks*{(-t) + K \cap B (t, c r)} \subset K \cap \frac{c r}{2} B_2$, namely $K \cap B (t, c r) \subset t + 2 \bracks*{K \cap \frac{c r}{2} B_2}$.
  This yields
  \begin{equation*}
    N \parens*{K \cap B (t, c r), r}
    \leq N \parens*{t + 2 \bracks*{K \cap \frac{c r}{2} B_2}, r}
    = N \parens*{K \cap \frac{c r}{2} B_2, \frac{r}{2}}
    \, .
    \qedhere
  \end{equation*}
\end{proof}

We can now state the main observation of this section: at scales $u = w (K)/c$ (for some absolute constant $c > 1$), the ``effective radius'' $r_*$ defined by local widths is smaller than the effective radius $\wt r$ defined by covering numbers.

\begin{lemma}
  \label{lem:radii-constant-factor}
  Let $c_0 \geq 1$ be the constant in Lemma~\ref{lem:sudakov} (Sudakov minoration).
  For every $\kappa_0 > 0$, 
  \begin{equation}
    \label{eq:radii-constant-factor}
    r_* \parens*{K, \frac{w (K)}{4 c_0 \kappa_0}}
    \leq c_0 \kappa_0 \, \wt r \parens*{K, \frac{w (K)}{\kappa_0}}
    \, .
  \end{equation}
\end{lemma}

\begin{proof}[Proof of Lemma~\ref{lem:radii-constant-factor}]
  If $\wt r (K, w (K)/\kappa_0) = + \infty$, then the bound~\eqref{eq:radii-constant-factor} trivially holds.
  We thus assume that $\wt r (K, w (K)/\kappa_0) < + \infty$.
  Let $r' > \wt r (K, w (K)/\kappa_0)$ be arbitrary; by definition of $\wt r$, there exists $r \leq r'$ such that
  \begin{equation}
    \label{eq:proof-radii-constant-1}
    r \sqrt{\log N (K, r)}
    \geq \frac{w (K)}{\kappa_0}
    \, .
  \end{equation}
  In addition, it is a straightforward consequence of the definition of covering numbers that, for every $c > 0$ (to be fixed later),
  \begin{equation}
    \label{eq:proof-radii-constant-2}
    N \parens*{K, r}
    \leq N \parens*{K, c r} \sup_{t \in K} N \parens*{K \cap B (t, c r), r}
    \leq N \parens*{K, c r} N \parens*{K \cap \frac{c r}{2} B_2, \frac{r}{2}}
    \, ,
  \end{equation}
  where the last inequality follows from Lemma~\ref{lem:local-cov-max}.
  Thus, combining~\eqref{eq:proof-radii-constant-1},~\eqref{eq:proof-radii-constant-2} and Sudakov minoration (Lemma~\ref{lem:sudakov}, for $K$ at scale $c r$ and $K \cap \frac{c r}{2} B_2$ at scale $r/2$), we obtain
  \begin{align*}
    \frac{w (K)}{\kappa_0}
    &\leq r \sqrt{\log N (K, r)}
      \leq r \sqrt{\log N (K, c r)} + r \sqrt{\log N \parens*{K \cap \frac{c r}{2} B_2, \frac{r}{2}}} \\
    &\leq \frac{c_0 w (K)}{c} + 2 c_0 w \parens*{K \cap \frac{c r}{2} B_2}
      \, .
  \end{align*}
  Letting $c = 2 c_0 \kappa_0 > 0$ and re-arranging gives
  \begin{equation*}
    w \parens*{K \cap c_0 \kappa_0 r B_2}
    \geq \frac{w (K)}{4 c_0 \kappa_0}
    \, .
  \end{equation*}
  This means that $c_0 \kappa_0 r \geq r_* \parens*{K, w (K)/(4 c_0 \kappa_0)}$ and thus $c_0 \kappa_0 r' \geq r_* \parens*{K, w (K)/(4 c_0 \kappa_2)}$.
  Since $r' > \wt r (K, w(K)/\kappa_0)$ was arbitrary, this concludes the proof.
\end{proof}

We now proceed to the proof of the lower bound in~\eqref{eq:small-ball-constant-regime}, which follows from the combination of Lemma~\ref{lem:lower-sb-width-cover} and Lemma~\ref{lem:radii-constant-factor}.
Let $c_1 \geq 1$ denote the constant $c$ of Lemma~\ref{lem:lower-sb-width-cover}, and $c_0$ the universal constant of Lemma~\ref{lem:sudakov}.

For $\eps \in \ooint{0}{3/4}$, let $\eps' \in (0,1/2)$ such that
$1-\eps = (1-\eps')^2$.
Thus $\eps = \eps' (2-\eps')$, and since $0< \eps' <1/2$ we deduce that $\frac{3}{2} \eps' < \eps < 2 \eps'$.
By Lemma~\ref{lem:lower-sb-width-cover}, for any $\kappa_2 > 0$,
\begin{align}
  \gamma_n (u K^\circ)
  &\geq \eps' \exp \parens*{- \frac{c_1 u^2}{\kappa_2^2 \min \braces[\big]{\eps'^2 r_* \parens*{K, (1-\eps')^2 u}^2, \wt r (K, u/\kappa_2)^2}}} \nonumber \\
  &\geq \frac{\eps}{2} \exp \parens*{- \frac{4 c_1 u^2}{\kappa_2^2 \min \braces[\big]{\eps^2 r_* \parens*{K, (1-\eps) u}^2, \wt r (K, u/\kappa_2)^2}}}
    \, .
    \label{eq:proof-constant-1}
\end{align}
Now from Lemma~\ref{lem:radii-constant-factor} and using that $w (K) / u \leq \kappa_0$ by assumption,
\begin{equation*}
  \wt r \parens*{K, \frac{u}{\kappa_2}}
  = \wt r \parens*{K, \frac{w(K)}{\kappa_2 w(K)/u}}
  \geq \frac{1}{c_0 \parens*{\kappa_2 w(K)/u}} r_* \parens*{K, \frac{u}{4 c_0 \kappa_2}}
  = \frac{1}{c_0 \kappa_0 \kappa_2} r_* \parens*{K, \frac{u}{4 c_0 \kappa_2}}
  \, .
\end{equation*}
Letting $\kappa_2 = (1-\eps)^{-1}/(4 c_0)$, we deduce that
\begin{equation}
  \label{eq:proof-constant-2}
  \wt r \parens*{K, \frac{u}{\kappa_2}}
  \geq \frac{4 (1-\eps)}{\kappa_0} r_* \parens*{K, (1-\eps) u}
  \geq \kappa_0^{-1} r_* \parens*{K, (1-\eps) u}
  \, .
\end{equation}
Plugging~\eqref{eq:proof-constant-2} into~\eqref{eq:proof-constant-1} and further using that $\kappa_2 \geq 1/(4 c_0)$ and $u \leq w (K)$, we obtain
\begin{equation}
  \label{eq:proof-constant-3}
  \gamma_n (u K^\circ)
  \geq \frac{\eps}{2} \exp \parens*{- \frac{64 c_0^2 c_1 \max\set{\eps^{-2}, \kappa_0^2} w(K)^2}{r_* \parens*{K, (1-\eps) u}^2}}
  \, .
\end{equation}

Finally, Lemma~\ref{lem:upper-sb-width} together with $u \geq w (K)/\kappa_0$ gives
\begin{equation}
  \label{eq:proof-constant-4}
  \gamma_n (u K^\circ)
  \leq \exp \parens*{- \frac{\eps^2 w(K)^2}{2 \kappa_0^2 r_* (K, (1+\eps) u)^2}}
  \, .
\end{equation}
From~\eqref{eq:proof-constant-3} and~\eqref{eq:proof-constant-4}, we conclude that~\eqref{eq:small-ball-constant-regime} holds with $\kappa_1 =
\max \{ 2 \kappa_0^2 \eps^{-2}, 64 c_0^2 c_1 \max\set{\eps^{-2}, \kappa_0^2} \}
$.

\section{Small ball probabilities, intrinsic volumes and the Wills functional}
\label{sec:conn-with-wills}

We now turn to the proofs of the results stated in Section~\ref{sec:intr-volum-wills}, concerning intrinsic volumes, the Wills functional, and their connections with small ball probabilities (Corollary~\ref{cor:link-sb-wills}, Theorem~\ref{thm:intrinsic-volumes}, Corollary~\ref{cor:reverse-alexandrov}, Corollary~\ref{cor:intrinsic-quantiles-gaussian} and Theorem~\ref{thm:intrinsic-quantile-exact}).

\subsection{Equivalences between fixed points}
\label{sec:equiv-fixed-points}

We first gather some technical lemmas, allowing to manipulate the quantities $r_*, \wt r, s_*, \wt s$, and establishing relationships between them.

\begin{fact}
  \label{fac:width-increasing}
  If $K_1, K_2$ are symmetric convex bodies in $\R^n$ with $K_1 \subset K_2$ strictly, then
  $w (K_1) < w (K_2)$.
\end{fact}

\begin{proof}
  Letting $X \sim \gaussdist (0, I_n)$, we have for $i=1, 2$, 
  \begin{equation*}
    w (K_i)
    = \int_0^{+\infty} \Probab*{\sup_{t \in K_i} \innerp{t}{X} \geq u} \di u
    = \int_0^{+\infty} \bracks*{1 - \gamma_n (u K_i^\circ)} \di u
    \, .
  \end{equation*}
  But by duality, for every $u > 0$ we have $u K_2^\circ \subset u K_1^\circ$ strictly.
  Now since $K_1$ is bounded, $K_1^\circ$ has nonempty interior and thus (by convexity and linear separation) $u K_1^\circ \setminus u K_2^\circ$ has nonempty interior, thus $\gamma_n (u K_1^\circ) - \gamma_n (u K_1^\circ) = \gamma_n (u K_1^\circ \setminus u K_2^\circ) > 0$.
  This proves the claimed inequality.
\end{proof}

\begin{lemma}
  \label{lem:loc-width-bij}
  For any symmetric convex body $K \subset \R^n$, the function $r \mapsto w (K \cap r B_2)$ induces a continuous and increasing bijection from $\ccint{0}{\rad (K)}$ to $\ccint{0}{w (K)}$.
  In particular, for any $u \in \ooint{0}{w(K)}$, one has $r = r_* (K, u) \in \ooint{0}{\rad (K)}$ and $w (K \cap r B_2) = u$.

  In addition, if $K$ has non-empty interior, the map $r \mapsto w (K \cap r B_2)/r$ induces a continuous and decreasing bijection from $\coint{\inrad (K)}{+\infty}$ to $\ocint{0}{w(B_2^n)}$.
\end{lemma}

\begin{proof}[Proof of Lemma~\ref{lem:loc-width-bij}]
  We first prove continuity.
  On the one hand, $w (K \cap r B_2)$ is increasing in $r$.
  On the other hand, for $0 < r < r'$, using that $K \subset \frac{r'}{r} K$ as $K$ is convex and contains $0$, we have $w (K \cap r' B_2) \leq w (\frac{r'}{r} K \cap r'B_2) = \frac{r'}{r} w (K \cap r B_2)$.
  These two properties imply that $r \mapsto w (K \cap r B_2)$ is continuous on $\R_+^*$.
  Continuity at $0$ follows from the bound $w (K \cap r B_2) \leq w(rB_2) \leq r \sqrt{n}$.

  Since $w (K \cap r B_2)$ equals $0$ at $r=0$ and $w (K)$ for $r \geq \rad (K)$, to prove the first claim it remains to show that this function is strictly increasing in $r \in \ccint{0}{\rad (K)}$.
  If $0 \leq r < r' \leq \rad (K)$, letting $t_0 \in K$ such that $\norm{t_0} = \rad (K)$ one has $\frac{r'}{\rad (K)} t_0 \in (K \cap r' B_2) \setminus (K \cap r B_2)$, thus $K \cap r B_2 \subset K \cap r' B_2$ strictly and (by Fact~\ref{fac:width-increasing}) $w(K \cap r B_2^n) < w (K \cap r' B_2)$.

  We now turn to the second claim.
  The assumption on $K$ implies that $\inrad (K) > 0$, and $r^{-1} w (K \cap rB_2^n)$ is continuous in $r$ since $w (K \cap r B_2^n)$ is.
  Also, for $r = \inrad (K)$, one has $w (K \cap r B_2^n)/r = w(B_2^n)$; while as $r \to + \infty$, one has $w (K \cap r B_2)/r \sim w (K)/r \to 0$.
  It thus remains to show that the map is strictly decreasing.
  This follows from the fact that $w (K \cap rB_2^n)/r = w (r^{-1} K \cap B_2^n)$ and that (as in the previous paragraph) if $\inrad (K) \leq r < r'$, then $r^{-1} K \cap B_2^n \supset r'^{-1} K \cap B_2^n$ strictly, together with Fact~\ref{fac:width-increasing}.
\end{proof}

The following claim is straightforward.
Recall the definitions of $r_*,\wt r$ in~(\ref{eq:def-radius-widths},\ref{eq:def-radius-cover}) and those of $s_*, \wt s$ in~(\ref{eq:def-dim-width},\ref{eq:def-dim-cover}).

\begin{lemma}
  \label{lem:link-fixed-r-k-u}
  For any convex body $K \subset \R^n$ and $u \in \ooint{0}{w(K)}$, one has
  \begin{align}
    \label{eq:link-r-k-u}
    s_* (K, u)
    = \frac{u^2}{r_* (K, u)^2}
    \quad \text{ and } \quad
    \wt s (K, u)
    = \frac{u^2}{\wt r (K, u)^2}
    \, .
  \end{align}
\end{lemma}

\begin{proof}[Proof of Lemma~\ref{lem:link-fixed-r-k-u}]
  Denoting $\wt R (K, u) = \set[\big]{r > 0 : r \sqrt{\log N (K, r)} \geq u}$, we have $\wt r (K, u) = \inf \wt R (K, u)$.
  Now for $s > 0$, we have $u / \sqrt{s} \in \wt R (K, u)$ if and only if
  \begin{equation*}
    \frac{u}{\sqrt{s}} \sqrt{\log N \parens*{K, \frac{u}{\sqrt{s}}}} \geq u
    \, ,
  \end{equation*}
  which is equivalent to $\log N (K, u/\sqrt{s}) \geq s$.
  Thus
  \begin{equation*}
    \set*{s > 0 : \log N \parens*{K, \frac{u}{\sqrt{s}}} \geq s}
    = \set*{\frac{u^2}{r^2} : r \in \wt R (K, u)}
    \, ,
  \end{equation*}
  and taking the supremum of each set gives $\wt s (K, u) = u^2 / \wt r (K, u)^2$.
  The equality $s_* (K, u) = u^2 / r_* (K, u)^2$ is immediate.
\end{proof}

It is convenient to introduce the following ``effective dimension'' of a symmetric convex body $K$ at a scale $r > 0$:
\begin{equation}
  \label{eq:def-effective-dim-both}
  \ol d (K, r)
  = \max \set*{w \parens*{r^{-1} K \cap B_2}^2, \log N (K, r)}
  \, .
\end{equation}
We then have the following:

\begin{lemma}
  \label{lem:max-k-eff-dim-fixed}
  For every convex body $K \subset \R^n$ and $u \in \ooint{0}{w(K)}$, we have
  \begin{equation}
    \label{eq:max-k-eff-dim-fixed}
    \max \braces*{s_* (K, u), \wt s (K, u)}
    = \sup \braces*{s > 0 : \ol d \parens*{K, \frac{u}{\sqrt{s}}} \geq s}
    \, .    
  \end{equation}
\end{lemma}

\begin{proof}[Proof of Lemma~\ref{lem:max-k-eff-dim-fixed}]
  For $u, s > 0$, we have
  \begin{equation*}
    w \parens*{K \cap \frac{u}{\sqrt{s}} B_2}
    \geq u
    \quad \text{ if and only if } \quad
    w \parens[\bigg]{\frac{\sqrt{s}}{u} K \cap B_2}^2
    \geq s
    \, .
  \end{equation*}
  Therefore, using that $\max \set{\sup A, \sup B} = \sup (A \cup B)$ for $A, B \subset \R$,
  \begin{align*}
    \max \braces*{s_* (K, u), \wt s (K, u)}
    &= \sup \braces*{s > 0 :
      w \parens[\bigg]{\frac{\sqrt{s}}{u} K \cap B_2}^2 \geq s
      \text{ or }
      \log N \parens*{K, \frac{u}{\sqrt{s}}} \geq s
      } \\
    &= \sup \braces*{s > 0 : \ol d \parens*{K, \frac{u}{\sqrt{s}}} \geq s}
      \, .
      \qedhere
  \end{align*}
\end{proof}

We also refer to Section~\ref{sec:addit-lemm-small} for basic properties of the fixed point $s_W$.

\subsection{Proof of Corollary~\ref{cor:link-sb-wills}
  %Correspondence with the Wills functional
}
\label{sec:corr-with-wills}

It was shown in~\cite[Corollary~A.1]{mourtada2025universal} (as a consequence of \cite[Theorem~4.2]{mourtada2025universal}, recalled in Theorem~\ref{thm:mourt-wills} above) that there exists a universal constant $c \geq 1$ such that, for every symmetric convex set $K \subset \R^n$,
\begin{equation}
  \label{eq:wills-eff-dim}
  \frac{1}{c} \inf_{r > 0} \braces*{\ol d (K, r) + r^2}
  \leq \log W (K)
  \leq c \inf_{r > 0} \braces*{\ol d (K, r) + r^2}
  \, .
\end{equation}

This readily implies the following.

\begin{lemma}
  \label{lem:wills-fixed-eff-dim}
  There exists a universal constant $c \geq 1$ such that the following holds.
  For any $n \geq 1$, symmetric convex body $K \subset \R^n$, and $u, s > 0$:
  \begin{enumerate}
  \item if $\ol d (K, u/\sqrt{s}) \geq s$, then $\log W (\frac{s}{u} K) \geq s/c$;
  \item conversely, if $\log W (\frac{s}{u} K) \geq s$, then $\ol d (K, u/\sqrt{c s}) \geq s/c$.
  \end{enumerate}
\end{lemma}

\begin{proof}[Proof of Lemma~\ref{lem:wills-fixed-eff-dim}]
  Let $c_0 \geq 1$ be the constant $c$ from the inequalities~\eqref{eq:wills-eff-dim}.
  For any $\lambda > 0$, using that $\ol d (\lambda K, \lambda r) = \ol d (K, r)$ for $r > 0$, we have
  \begin{equation}
    \label{eq:proof-wills-lambda-eff-dim}
    \frac{1}{c_0} \inf_{r > 0} \braces*{\ol d (K, r) + \lambda^2 r^2}
    \leq \log W (\lambda K)
    \leq c_0 \inf_{r > 0} \braces*{\ol d (K, r) + \lambda^2 r^2}
    \, .
  \end{equation}

  Assume that $\ol d (K, u/\sqrt{s}) \geq s$.
  Then for any $0 < r \leq u/\sqrt{s}$, noting that $\ol d (K, r)$ decreases in $r$, we have $\ol d(K, r) \geq s$.
  Also, for any $r \geq u/\sqrt{s}$ we have $(s/u)^2 r^2 \geq (s/u)^2 (u/\sqrt{s})^2 = s$.
  We thus have
  \begin{equation*}
    \inf_{r > 0} \braces[\bigg]{\ol d (K, r) + \parens*{\frac{s}{u}}^2 r^2}
    \geq s
    \, ,
  \end{equation*}
  and combined with~\eqref{eq:proof-wills-lambda-eff-dim} this implies that $\log W (\frac{s}{u} K) \geq s/c_0$.

  Conversely, assume now that $\log W (\frac{s}{u} K) \geq s$.
  Then from~\eqref{eq:proof-wills-lambda-eff-dim}, we have for every $r > 0$,
  \begin{equation*}
    \ol d \parens*{K, r} + \parens*{\frac{s}{u}}^2 r^2
    \geq \frac{s}{c_0}
    \, .
  \end{equation*}
  Letting $r = u / \sqrt{2 c_0 s}$ and re-arranging, this yields $\ol d (K, u/\sqrt{2 c_0 s}) \geq s/(2 c_0)$.
\end{proof}

We are now in position to prove Corollary~\ref{cor:link-sb-wills}.

\begin{proof}[Proof of Corollary~\ref{cor:link-sb-wills}]
  First, note that for every $u > 0$, by definition of $\lambda (K, u)$ we have
  \begin{equation}
    \label{eq:proof-link-lambda-k-wills}
    \lambda (K, u) u
    = s_W (K, u)
    = \sup \braces*{s > 0 : \log W \parens*{\frac{s}{u} K} \geq s}
    \, .
  \end{equation}

  Now, let $c_2 \geq 1$ be the constant $c$ from Lemma~\ref{lem:wills-fixed-eff-dim}.
  Thus, for any $s > 0$ such that $\ol d (K, u/\sqrt{s}) \geq s$, we have $\log W (\frac{s/c_2}{u/c_2} K) = \log W (\frac{s}{u} K) \geq s/c_2$.
  Taking the supremum over $s > 0$ satisfying the first condition and using Lemma~\ref{lem:max-k-eff-dim-fixed}, we deduce that
  \begin{equation}
    \label{eq:proof-sb-wills-1}
    \frac{1}{c_2} \max \braces*{s_* \parens*{K, u}, \wt s \parens*{K, u}}   
    \leq s_W \parens*{K, \frac{u}{c_2}}
    \, .
  \end{equation}
  In addition, still from Lemma~\ref{lem:wills-fixed-eff-dim}, if $\log W (\frac{s}{u} K) \geq s$ then $\ol d (K, (u/c_2)/\sqrt{s/c_2}) = \ol d (K, u/\sqrt{c_2 s}) \geq s/c_2$.
  Again, taking the supremum over $s$ and using Lemma~\ref{lem:max-k-eff-dim-fixed}, this gives
  \begin{equation}
    \label{eq:proof-sb-wills-2}
    \frac{1}{c_2} s_W \parens*{K, u}
    \leq \max \braces*{s_* \parens*{K, \frac{u}{c_2}}, \wt s \parens*{K, \frac{u}{c_2}}}
    \, .
  \end{equation}
  In summary, from~\eqref{eq:proof-link-lambda-k-wills},~\eqref{eq:proof-sb-wills-1} and~\eqref{eq:proof-sb-wills-2}, we have for every $u > 0$:
  \begin{equation}
    \label{eq:proof-sb-wills-link-k}
    \lambda \parens*{K, c_2 u} u
    \leq \max \braces*{s_* \parens*{K, u}, \wt s \parens*{K, u}}
    \leq \lambda \parens*{K, \frac{u}{c_2}} u
    \, .
  \end{equation}

  Now, from Theorem~\ref{thm:gaussian-small-ball} together with Lemma~\ref{lem:link-fixed-r-k-u}, for some universal constants $c_0, c_1 \geq 1$, the following holds.
  For every $u \in \ooint{0}{w (K)/c_0}$, we have
  \begin{equation}
    \label{eq:proof-sb-wills-sb-k}
    \frac{1}{c_1} \max \braces*{s_* \parens*{K, c_1 u}, \wt s \parens*{K, c_1 u}}
    \leq - \log \gamma_n (u K^\circ)
    \leq c_1 \max \braces*{s_* \parens*{K, \frac{u}{c_1}}, \wt s \parens*{K, \frac{u}{c_1}}}
    \, .
  \end{equation}
  Combining~\eqref{eq:proof-sb-wills-link-k} and~\eqref{eq:proof-sb-wills-sb-k} completes the proof.
\end{proof}

\subsection{Recovering an ultra-log-concave sequence from its generating function}
\label{sec:recovering-an-ultra}

We now proceed towards the proof of Theorem~\ref{thm:intrinsic-volumes}.
Towards this end, in this section we show how one may approximately recover an ultra-log-concave sequence from its log-generating function.

\begin{definition}
  \label{def:ulc}
  A sequence $(v_j)_{j \geq 0} \in \R_+^{\N}$ of non-negative numbers is \emph{ultra-log-concave} if:
  \begin{enumerate}
  \item it contains no internal zeros: for every integers $0 \leq j \leq k \leq l$, if $v_j, v_l > 0$ then $v_k > 0$;
  \item the sequence $(j! v_j)_{j \geq 0}$ is log-concave, namely $j v_j^2 \geq (j+1) v_{j-1} v_{j+1}$ for every $j \geq 1$.
  \end{enumerate}
\end{definition}

\begin{lemma}
  \label{lem:ulc-seq-from-mgf}
  Let $(v_k)_{k \geq 0}$ be an ultra-log-concave sequence such that $v_0 = 1$.
%  For any $k \geq 1$, let $\theta_k = (k! v_k)^{1/k}$.
  Define the function $\phi : \R^+ \to \R$ by $\phi (\lambda) = \log \parens[\big]{\sum_{k \geq 0} v_k \lambda^k}$ for $\lambda \in \R^+$.
  Then
  % the sequence $(\theta_k)_{k \geq 0}$ is decreasing, and
  for every $k \geq 1$, we have
  \begin{equation}
    \label{eq:ulc-seq-from-mgf}
    v_k
    \leq \inf_{\lambda > 0} \frac{\exp \phi (\lambda)}{\lambda^k}
    \leq (e/k)^k k! v_k
%    \leq e^k v_k
    \, .
  \end{equation}
\end{lemma}

% In fact, in the right-most side of~\eqref{eq:ulc-seq-from-mgf}, the factor $e^k$ could be replaced by $e \sqrt{k}$ for $k \geq 1$, but this improvement would be inconsequential for our purposes.

\begin{proof}[Proof of Lemma~\ref{lem:ulc-seq-from-mgf}]
  The lower bound is straightforward: for every $\lambda > 0$, one has
  $\exp \phi (\lambda) = \sum_{j \geq 0} v_j \lambda^j \geq v_k \lambda^k$.
  We now turn to the upper bound.

  If $v_k = 0$, then the lack of internal zeros implies that $v_j = 0$ for every $j \geq k$, thus $\exp \phi (\lambda) = \sum_{j=0}^{k-1} v_j \lambda^j  = O_v (\lambda^{k-1})$ as $\lambda \to \infty$, and thus $\lim_{\lambda \to +\infty} \exp \parens{\phi (\lambda)}/\lambda^k = 0$ and~\eqref{eq:ulc-seq-from-mgf} holds.

  Assume now that $v_k > 0$, and let $\rho = k v_{k}/v_{k-1} = (k! v_k)/((k-1)! v_{k-1}) > 0$.
  Now, ultra-log-concavity of $(v_j)_{j \geq 0}$ implies that the sequence $((j! v_j)/((j-1)! v_{j-1}))_{j \geq 1}$ is decreasing, thus for every $j \geq 0$ one has
  \begin{equation*}
    j! v_j
    \leq k! v_k \rho^{j-k}
    \, .
  \end{equation*}
  Summing the above inequality yields
  \begin{equation*}
    \exp \parens*{\phi (\lambda)}
    = \sum_{j \geq 0} v_j \lambda^j
    \leq k! v_k \rho^{-k} \sum_{j \geq 0} \frac{\rho^j \lambda^j}{j!}
    = k! v_k \rho^{-k} \exp \parens*{\lambda \rho}
    \, .
  \end{equation*}
  Letting $\lambda = k/\rho$ above and bounding $k! \leq k^k$, we conclude that
  \begin{equation*}
    \inf_{\lambda > 0} \frac{\exp \phi (\lambda)}{\lambda^k}
    \leq k! v_k k^{-k} \exp (k)
%    \leq e^k v_k
    \, .
    \qedhere
  \end{equation*}
\end{proof}

\begin{lemma}
  \label{lem:fixed-log-mgf-theta}
  Let $(v_{j})_{j \geq 0}$ be an ultra-log-concave sequence such that $v_0 = 1$ and define $\phi$ as in Lemma~\ref{lem:ulc-seq-from-mgf}.
  Then, the sequence $\parens*{\theta_j}_{j \geq 1}$ defined by $\theta_j = (j! v_j)^{1/j}$ for $j \geq 1$ is decreasing.
  In addition, for every integer $j \geq 1$ and $u \in \R_+^*$, if $\theta_j \geq u$ then for every $\xi > 1$,
  \begin{equation}
    \label{eq:fixed-log-mgf}
    \phi \parens*{\frac{\xi j}{u}}
    \geq j \log \xi
    \, .
  \end{equation}
  Conversely, if~\eqref{eq:fixed-log-mgf} holds for every $\xi > 1$, then $\theta_j \geq u/e$.
\end{lemma}

\begin{proof}[Proof of Lemma~\ref{lem:fixed-log-mgf-theta}]
  First, the fact that $(\theta_j)_{j \geq 1}$ decreases readily follows from ultra-log-concavity of $(v_j)_{j \geq 0}$ and from the fact that $v_0 = 1$.
  % We first show that $\parens*{\theta_j}_{j \geq 1}$ decreases.
  % The inequality $\theta_j \leq \theta_{j+1}$ is obvious if $v_j = 0$, we thus assume that $v_j > 0$; thus $v_{j} > 0$ for $j=0, \dots, j$.
  % For any $j \geq 1$, ultra-log-concavity implies that the sequence $(j! v_{j}/((j-1)! v_{j-1}))_{1 \leq j \leq j+1}$ is decreasing, thus
  % \begin{equation*}
  %   \parens*{\frac{(j+1)! v_{j+1}}{j! v_j}}^{j}
  %   \leq \prod_{j=1}^{j} \frac{j! v_j}{(j-1)! v_j}
  %   = j! v_j
  % \end{equation*}

  Now, assume that $\theta_j \geq u$, and let us show that~\eqref{eq:fixed-log-mgf} holds.
  Using that $j! \leq j^j$, we deduce that for every $\xi \geq 1$,
  \begin{equation*}
    \phi \parens*{\frac{\xi j}{u}}
    = \log \parens[\bigg]{\sum_{j \geq 0} \frac{\theta_j^j}{j!} \parens[\bigg]{\frac{\xi j}{u}}^j}
    \geq \log \parens[\bigg]{\frac{\theta_j^j}{j!} \parens*{\frac{\xi j}{u}}^j}
    \geq \log \parens[\bigg]{\frac{u^j}{j^j} \parens*{\frac{\xi j}{u}}^j}
    = j \log \xi
    \, ,
  \end{equation*}
  which is the desired inequality.

  Conversely, assume now that~\eqref{eq:fixed-log-mgf} holds for every $\xi > 1$.
  This inequality then holds for every $\xi > 0$, as it is trivial for $\xi \leq 1$ given that $\phi \geq \log v_0 = 0$.
  This can be rewritten as $\exp \phi ({\xi} j/u) \geq \xi^j$,
  thus
  \begin{equation*}
    \inf_{\lambda > 0} \frac{\exp \phi (\lambda)}{\lambda^j}
    = \inf_{\xi > 0} \frac{\exp \phi \parens*{{\xi} j/u}}{({\xi} j/u)^j}
    \geq \parens*{\frac{u}{j}}^j
    \, .
  \end{equation*}
  But by Lemma~\ref{lem:ulc-seq-from-mgf}, the left-hand side above is at most $(e/j)^j j! v_j = (e \theta_j / j)^j$, thus $\theta_j \geq u/e$.
\end{proof}

\subsection{Proof of Theorem~\ref{thm:intrinsic-volumes}}
\label{sec:two-sided-estimates}

Recall the definition in~\eqref{eq:def-effective-dim-both} (from Section~\ref{sec:equiv-fixed-points}) of the effective dimension $\ol d (K, r)$ at scale $r > 0$.

Combining Lemmas~\ref{lem:wills-fixed-eff-dim} and~\ref{lem:fixed-log-mgf-theta} yields the following characterization of intrinsic volumes in terms of the behavior of the effective dimension $\ol d (K, r)$ for $r > 0$.

\begin{lemma}
  \label{lem:intr-eff-dim}
  There exist universal constants $c_0, c \geq 1$ such that the following holds.
  For every $n \geq 1$, convex body $K \subset \R^n$, real $u > 0$ and integer $j \geq 1$:
  \begin{enumerate}
  \item if $\theta_j (K) \geq u$, then for every $\xi > 1$,
    \begin{equation}
      \label{eq:intr-vol-eff-dim}
      \ol d \parens*{K, \frac{u}{c \xi \sqrt{j}}}
      \geq \frac{j \log \xi}{c}
      \, ;
    \end{equation}
  \item conversely, if $j \geq c_0$ and, for every $\xi > 1$,
    \begin{equation}
      \label{eq:intr-vol-eff-dim-cond}
      \ol d \parens*{K, \frac{u}{\xi \sqrt{j}}}
      \geq {j \log \xi}
      \, ,     
    \end{equation}
    then $\theta_{\floor{j/c_0}} (K) \geq u/c$.
  \end{enumerate}
\end{lemma}

Theorem~\ref{thm:intrinsic-volumes} immediately follows from Lemma~\ref{lem:intr-eff-dim}, using that for $r, j_0 > 0$, the inequality $\ol d (K, r) \geq j_0$ is equivalent to
\begin{equation*}
  \log N (K, r)
  \geq j_0
  \quad \text{ or } \quad
  w \parens*{K \cap r B_2}
  \geq r \sqrt{j_0}
  \, .
\end{equation*}

\begin{proof}[Proof of Lemma~\ref{lem:intr-eff-dim}]
  Let $v_j = V_j (K/\sqrt{2\pi})$ and $\theta_j = (j! v_j)^{1/j} = \theta_j (K)$ for every $j \geq 1$.
  Let $c_1 \geq 1$ denote the universal constant $c$ from Lemma~\ref{lem:wills-fixed-eff-dim}.

  Assume first that $\theta_j \geq u$.
  Then by Lemma~\ref{lem:fixed-log-mgf-theta}, for any $\zeta > 1$ one has
  \begin{equation*}
    \log W \parens*{\frac{j \log (\zeta)}{\log (\zeta) \zeta^{-1} u}  K}
    = \phi \parens*{\frac{\zeta j}{u}} 
    \geq j \log \zeta
    \, .
  \end{equation*}
  By Lemma~\ref{lem:wills-fixed-eff-dim}, this implies that
  \begin{equation*}    
    \ol d \parens*{K, \frac{\log (\zeta) \zeta^{-1} u}{\sqrt{c_1 j \log (\zeta)}}}
    = \ol d \parens*{K, \frac{\sqrt{\log \zeta}}{\zeta} \frac{u}{\sqrt{c_1 j}}}
    \geq \frac{j \log (\zeta)}{c_1}
    \, .
  \end{equation*}
  Now, let $\xi > 1$ and set $\zeta = e \xi > e$.
  Since $\sqrt{\log \zeta} \geq 1$ and $\ol d (K, r)$ decreases in $r$, the above inequality implies that
  \begin{equation*}
    \ol d \parens*{K, \frac{u}{e \xi \sqrt{c_1 j}}}
    \geq \frac{j \log (e \xi)}{c_1}
    \geq \frac{j \log \xi}{c_1}
    \, ,
  \end{equation*}
  which proves~\eqref{eq:intr-vol-eff-dim} with $c = \max \set{ e \sqrt{c_1} , c_1}$.

  Conversely, assume that~\eqref{eq:intr-vol-eff-dim-cond} holds for every $\xi > 1$.
  Then Lemma~\ref{lem:wills-fixed-eff-dim} implies that
  \begin{equation*}
    \phi \parens*{\xi \sqrt{\log \xi} \, \frac{j}{u}}
    = \log W \parens*{\frac{j \log \xi}{\xi^{-1} \sqrt{\log \xi} u}}
    \geq \frac{j \log \xi}{c_1}
    \, .
  \end{equation*}
  Now, the map $\xi \mapsto \xi \sqrt{\log \xi}$ is a continuous and increasing bijection from $\ooint{1}{+\infty}$ to $\R_+^*$.
  Thus, for any $\zeta > 0$, let $\xi > 1$ be such that $\xi \sqrt{\log \xi} = \zeta$.
  In particular, one has $\zeta \leq \xi \sqrt{\xi}$ and thus $\log \zeta \leq \frac{3}{2} \log \xi$.
  The above inequality then implies that
  \begin{equation*}
    % \label{eq:proof-intr-vol-dim-cond-phi}
    \phi \parens*{\frac{\zeta j}{u}}
    \geq \frac{2 j \log \zeta}{3 c_1}
    \, .
  \end{equation*}
  Assume now that $j \geq 3 c_1/2$ and let $j' = \floor{2 j / (3 c_1)}$.
  We thus have $j' \geq 1$ and
  % $2 j/(3 c_1) \geq j' \geq \frac{1}{2} (2 j / (3 c_1)) = j / (3 c_1)$.
  % $j' \leq 2 j/(3 c_1)$.
  $j \geq 3 c_1 j'/2$.
  Since $\phi$ is decreasing, we deduce that for every $\zeta > 0$,
  \begin{equation*}
    \phi \parens*{\frac{3 c_1 \zeta j'}{2 u}}
    \geq j' \log \zeta
    \, .
  \end{equation*}
  By Lemma~\ref{lem:fixed-log-mgf-theta}, the above condition implies that
  \begin{equation}
    \label{eq:proof-intr-vol-final}
    \theta_{j'}
    \geq \frac{2 u}{3 e c_1}
    \, .
  \end{equation}
  This establishes the second property with $c_0 = 3 c_1 / 2$ and $c = 3 e c_1 /2$.
\end{proof}

\subsection{Proof of Corollary~\ref{cor:reverse-alexandrov}
  % Reversal of Alexandrov-Fenchel%
  % Critical dimension for the decay of intrinsic volumes
  % Proof of Corollary~\ref{cor:reverse-alexandrov}
}
\label{sec:proof-reverse-alexandrov}

\begin{lemma}
  \label{lem:recov-from-mgf-constant}
  For any constant $\kappa_0 > 0$, there exists a constant $\kappa_1 > 0$ (depending only on $\kappa_0$) such that the following holds.
  Let $(v_j)_{j \geq 0}$ be an ultra-log-concave sequence with $v_0 = 1$; define $\phi : \R^+ \to \R$ as in Lemma~\ref{lem:ulc-seq-from-mgf} and $(\theta_j)_{j \geq 1}$ as in Lemma~\ref{lem:fixed-log-mgf-theta}.
  For any real $j \geq 1$, if 
  \begin{equation}
    \label{eq:cond-recov-mgf-constant}
    \phi \parens*{\frac{e j}{\theta_1}}
    \geq \frac{j}{\kappa_0}
    \, ,
  \end{equation}
  then $\theta_{\ceil{j/\kappa_1}} \geq \theta_1 / \kappa_1$.
\end{lemma}

Note that the first part of Lemma~\ref{lem:fixed-log-mgf-theta}
% (with $u = e \theta_1/\kappa_1$ and $\xi = e$) already shows that if $j \geq 1$ is an integer such that $\theta_j \geq e \theta_1 / \kappa_1$, then $\phi \parens*{\kappa_1 j / \theta_1} \geq j$.
(with $\xi = e$)
already shows that if $j \geq 1$ is an integer such that $\theta_j \geq u$, then $\phi (e j / u) \geq j$.
Lemma~\ref{lem:recov-from-mgf-constant} shows that when $u$ is of the same order of magnitude as $\theta_1$, this statement admits a converse.
This contrasts with the case where $u \ll \theta_1$, where the stronger condition of Lemma~\ref{lem:fixed-log-mgf-theta} (valid for all $\xi > 1$) is required.

\begin{proof}[Proof of Lemma~\ref{lem:recov-from-mgf-constant}]
  Let $\lambda = e j / \theta_1$.
  Since the sequence $(v_j \lambda^j)_{j \geq 0}$ is ultra-log-concave with $v_1 \lambda = e j \geq 2$, by \cite[Lemma~A.1]{mourtada2025universal} (the assumption of finite support is not used in the proof of this result), one has
  \begin{equation}
    \label{eq:proof-mgf-cst-log-max}
    \phi (\lambda)
    = \log \parens[\Bigg]{\sum_{j \geq 0} v_j \lambda^j}
    \leq 4 \log \parens*{\sup_{j \geq 1} v_j \lambda^j}
    \leq 4 \sup_{j \geq 1} j \log \parens*{\frac{e \theta_j \lambda}{j}}
    \, ,
  \end{equation}
  where we used that $v_j = \theta_j^j / j! \leq (e/j)^j \theta_j^j$.

  Now, let $\kappa_1 > 1$ to be fixed later, and $j_1 = j/\kappa_1$.
  Using that $\theta_j \leq \theta_1$ as $(v_j)_{j \geq 0}$ is ultra-log-concave with $v_0=1$, and that the map $s \mapsto s \log (e s_0/s)$ is increasing over $\ocint{0}{s_0}$ for $s_0 > 0$, we have
  \begin{align*}
    \max_{1 \leq j \leq j_1} j \log \parens*{\frac{e \theta_j \lambda}{j}}
    &\leq \max_{1 \leq j \leq j_1} j \log \parens*{\frac{e \theta_1 \lambda}{j}}
      = \max_{1 \leq j \leq j_1} j \log \parens*{\frac{e^2 j}{j}} \\
    &\leq j_1 \log \parens*{\frac{e^2 j}{j_1}}
      = \kappa_1^{-1} \log (e^2 \kappa_1) j
      \, .
  \end{align*}
  Thus, letting $\kappa_1 > 1$ be large enough such that $4 \kappa_1^{-1} \log (e^2 \kappa_1) < \kappa_0^{-1}$, we have
  \begin{equation}
    \label{eq:proof-mgf-cst-small-j}
    4 \max_{1 \leq j \leq j_1} j \log \parens*{\frac{e \theta_j \lambda}{j}}
    < \frac{j}{\kappa_0}
    \, .
  \end{equation}

  On the other hand, using again that $(\theta_{j})_{j\geq 1}$ decreases, we have for every integer $j \geq j_1$,
  \begin{equation*}
    \frac{e \theta_j \lambda}{j}
    \leq \frac{e \theta_{\ceil{j_1}}}{j_1} \frac{e j}{\theta_1}    
    = \frac{e^2 \kappa_1 \theta_{\ceil{j/\kappa_1}}}{\theta_1}
    \, .
  \end{equation*}
  Hence, if
  \begin{equation*}
    \theta_{\ceil{j/\kappa_1}}
    < \frac{\theta_1}{e^2 \kappa_1}
    \, ,
  \end{equation*}
  then
  \begin{equation*}
    4 \sup_{j \geq j_1} j \log \parens*{\frac{e \theta_j \lambda}{j}}
    \leq 0
    \, ,
  \end{equation*}
  and thus together with~\eqref{eq:proof-mgf-cst-log-max} and~\eqref{eq:proof-mgf-cst-small-j},
  \begin{equation*}
    \phi \parens*{\frac{e j}{\theta_1}}
    < \frac{j}{\kappa_0}
    \, .
  \end{equation*}
  Hence if~\eqref{eq:cond-recov-mgf-constant} holds, then $\theta_{\ceil{j/(e^2 \kappa_1)}} \geq \theta_{\ceil{j/\kappa_1}} \geq \theta_1 / (e^2 \kappa_1)$.
\end{proof}

We are now in position to prove Corollary~\ref{cor:reverse-alexandrov}.

\begin{proof}[Proof of Corollary~\ref{cor:reverse-alexandrov}]
  Again, for $j \in \N$ we let $v_j = V_j (K/\sqrt{2\pi})$, so that $\theta_j = \theta_j (K)$ for $j \geq 1$ and $\phi (\lambda) = \log W (\lambda K)$ for $\lambda \in \R^+$.
  We let $c_1 \geq 1$ denote the constant $c$ of Lemma~\ref{lem:wills-fixed-eff-dim}, $c_2 \geq 1$ the constant $c_0$ in Lemma~\ref{lem:radii-constant-factor}, and $c_3 \geq 1$ the constant $c$ in Lemma~\ref{lem:intr-eff-dim}.

  We start with the first direction.
  If
  \begin{equation*}
    j
    \leq \frac{w (K)^2}{r_* \parens*{w (K)/\kappa_0}^2}
    \, ,
  \end{equation*}
  then $w (K)/\sqrt{j} \geq r_* \parens*{w (K)/\kappa_0}$.
  Hence, by definition of $r_*$ and Lemma~\ref{lem:loc-width-bij}, one has
  \begin{equation*}
    w \parens*{K \cap \frac{w (K)}{\sqrt{j}} B_2}
    \geq \frac{w (K)}{\kappa_0}
    \, , \quad \text{ i.e.~} \quad
    w \parens*{\frac{\sqrt{j}}{w (K)} K \cap B_2}
    \geq \frac{\sqrt{j}}{\kappa_0}
    \, .
  \end{equation*}
  In particular, $\ol d \parens[\big]{K, w (K)/\sqrt{j}} \geq j/\kappa_0^2$.
  By Lemma~\ref{lem:wills-fixed-eff-dim} and since $\theta_1 = w(K)$, this gives 
  \begin{equation*}
    \phi \parens*{\frac{j}{\kappa_0 \theta_1}}
    = \log W \parens*{\frac{j/\kappa_0^2}{w (K)/\kappa_0} K}
    \geq \frac{j}{c_1 \kappa_0^2}
    \, .
  \end{equation*}
  By Lemma~\ref{lem:recov-from-mgf-constant}, this implies that (for some constant $\kappa_2 > 1$ depending on $\kappa_0, c_1$)
  \begin{equation*}
    \theta_{\ceil{j/\kappa_2}}
    \geq \frac{\theta_1}{\kappa_2}
    = \frac{w (K)}{\kappa_2}
    \, .
  \end{equation*}

  We now turn to the opposite direction.
  Let $j \geq 1$ be an integer and $\kappa_1 > 1$, and assume that $\theta_j (K) \geq w (K)/\kappa_1$.
  % , that is $\theta_j \geq \theta_1/\kappa_1$.
  % By Lemma~\ref{lem:fixed-log-mgf-theta} with $\xi = e$, this implies that
  % \begin{equation*}
  %   \phi \parens*{\frac{e j}{\theta_1/\kappa_1}}    
  %   = \log W \parens*{\frac{e \kappa_1 j}{\theta_1} K}
  %   \geq j
  %   \, .
  % \end{equation*}
  By Lemma~\ref{lem:intr-eff-dim} with $\xi = e$, this implies that
  \begin{equation*}
    \ol d \parens*{K, \frac{w (K)/\kappa_1}{c_3 e \sqrt{j}}}
    \geq \frac{j}{c_3}
    \, .
  \end{equation*}
  In other words, $\ol d \parens[\big]{K, \braces{w (K)/(e c_3^{3/2} \kappa_1)}/{\sqrt{j/c_3}}} \geq j/c_3$.
  Thus, by Lemma~\ref{lem:max-k-eff-dim-fixed} followed by Lemma~\ref{lem:link-fixed-r-k-u},
  \begin{align*}
    \frac{j}{c_3}
    &\leq \max \braces*{s_* \parens*{K, \frac{w (K)}{e c_3^{3/2} \kappa_1}}, \wt s \parens*{K, \frac{w (K)}{e c_3^{3/2} \kappa_1}}} \\
    &= \frac{w (K)^2}{e^2 c_3^3 \kappa_1^2 \min \braces*{r_* \parens*{K, {w (K)}/\bracks[\big]{e c_3^{3/2} \kappa_1}}, \wt r \parens*{K, {w (K)}/\bracks[\big]{e c_3^{3/2} \kappa_1}}}^2}
      \, .
  \end{align*}
  Finally, by Lemma~\ref{lem:radii-constant-factor},
  \begin{equation*}
    \wt r \parens*{K, {w (K)}/\bracks[\big]{e c_3^{3/2} \kappa_1}}
    \geq \frac{1}{c_2 {e c_3^{3/2} \kappa_1}} r_* \parens*{K, {w (K)}/\bracks[\big]{4 e c_2 c_3^{3/2} \kappa_1}}
    \, .
  \end{equation*}
  Combining the above two inequalities, we deduce that for some constant $\kappa_3, \kappa_4 > 1$ depending on $\kappa_1$ (and on the universal constants $c_2, c_3$), we have
  \begin{equation*}
    j
    \leq \frac{\kappa_3 \, w (K)^2}{r_* (K, w (K)/\kappa_4)^2}
    \, .
    \qedhere
  \end{equation*}
\end{proof}

\subsection{Proof of Corollary~\ref{cor:intrinsic-quantiles-gaussian}}
\label{sec:proof-intr-volum-quant}

% \begin{proof}[Proof of Corollary~\ref{cor:intrinsic-quantiles-gaussian}]
Let $c_0, c_1 \geq 1$ denote the constants $c_0, c$ of Theorem~\ref{thm:gaussian-small-ball} and $c_2, c_3 \geq 1$ the constants $c_0, c$ of Lemma~\ref{lem:intr-eff-dim}.

We start with the upper bound.
Let $j \geq 1$ be an integer and $u = \theta_j (K)/c_0$.
Since $\theta_j (K) \geq u$, the first part of Lemma~\ref{lem:intr-eff-dim} implies that, for every $\zeta > 1$,
\begin{equation*}
  \ol d \parens*{K, \frac{u}{c_3 \xi \sqrt{j}}}
  \geq \frac{j \log \xi}{c_3}
  \, .
\end{equation*}
% which also trivially holds for $\zeta \in \ocint{0}{1}$.
% Now, let $\xi > 1$, and set $\zeta > e$ be such that $\xi = \zeta/(e )$.
% The above then implies, using that $\ol d (K, \cdot)$ decreases, that
% Now, let $\xi \geq 1$, and $\zeta = e \xi$.
Assume now that $\xi \geq e$.
The above then implies, using that $\log \xi \geq 1$ and $\ol d (K, \cdot)$ decreases, that
\begin{equation*}
  \ol d \parens*{K, \frac{u/\bracks[\big]{c_3^{3/2} \xi}}{\sqrt{j \log (\xi) / c_3}}}
  \geq \ol d \parens*{K, \frac{u}{c_3 \xi \sqrt{j}}}
  \geq \frac{j \log (\xi)}{c_3}
  \, .
\end{equation*}
By Lemma~\ref{lem:max-k-eff-dim-fixed}, this implies that
\begin{equation*}
  \frac{j \log \xi}{c_3}
  \leq \max \braces*{s_* \parens*{K, \frac{u}{c_3^{3/2} \xi}}, \wt s \parens*{K, \frac{u}{c_3^{3/2} \xi}}}
  \, .
\end{equation*}
In addition, by Theorem~\ref{thm:gaussian-small-ball} we have (the inequality being valid as $u / (c_1 c_3^{3/2} \xi) < u = \theta_j (K)/c_0 \leq w (K)/c_0$)
\begin{equation*}
  \max \braces*{s_* \parens*{K, \frac{u}{c_3^{3/2} \xi}}, \wt s \parens*{K, \frac{u}{c_3^{3/2} \xi}}}
  \leq - c_1 \log \gamma_n \parens*{\frac{u}{c_1 c_3^{3/2} \xi} K^\circ}
  \, .
\end{equation*}
Combining the previous two inequalities, we deduce that for some universal constant $c_4 \geq 1$, we have for every $\xi \geq e$ that
\begin{equation*}
  \gamma_n \parens*{\frac{u}{c_4 \xi} K^\circ}
  \leq \xi^{-j/c_4}
  \, .
\end{equation*}
In other words,
\begin{equation*}
  \sup_{0 < u' \leq u/(e c_4)} (u')^{-j/c} \gamma_n (u' K^\circ)
  = \sup_{\xi \geq e} \bracks*{\parens*{\frac{u}{c_4 \xi}}^{-j/c_4} \gamma_n \parens*{\frac{u}{c_4 \xi} K^\circ}}
  \leq \parens*{\frac{u}{c_4}}^{-j/c_4}
  \, .
\end{equation*}
On the other hand,
\begin{equation*}
  \sup_{u' \geq u/(e c_4)} (u')^{-j/c_4} \gamma_n (u' K^\circ)
  \leq \parens*{\frac{u}{e c_4}}^{-j/c_4}
  \, .
\end{equation*}
Combining the previous two inequalities and using the definition~\eqref{eq:def-scaled-quantile} of $\wt \theta_{j/c_4} (K)$, we get
\begin{equation*}
  \wt \theta_{j/c_4} (K)
  \geq \frac{u}{e c_4}
  = \frac{\theta_j (K)}{e c_0 c_4}
  \, ,
\end{equation*}
which gives the second inequality in~\eqref{eq:intrinsic-quantiles-gaussian} with $c = e c_0 c_4$ as $j' \mapsto \wt \theta_{j'} (K)$ is decreasing.

We now turn to the lower bound.
Let $c_5, c_6 \geq 1$ denote the constants $c_0, c$ of Corollary~\ref{cor:link-sb-wills}.
Let $u = \wt \theta_{4 j} (K)/(16 c_5)$.
% , for some constant $c_7 \geq 1$ to be chosen later.
Since $\wt \theta_{4 j} (K) \geq 4 u$, we have for every $\xi > 0$ that
\begin{equation*}
  \gamma_n \parens*{\frac{4 u}{\xi} K^\circ}
  \leq \xi^{-4 j}
  \, , \quad \text{ that is } \quad
  - \log \gamma_n \parens*{\frac{4 u}{\xi} K^\circ}
  \geq 4 j \log \xi
  \, .
\end{equation*}
In addition, note that for $\xi > 1$, we have (using that $4 j \geq 1$ and setting $\eps = 1/2$ in~\eqref{eq:def-scaled-quantile})
\begin{equation*}
  \frac{4 u}{\xi}
  < \frac{\wt \theta_{4 j} (K)}{4 c_5}
  \leq \frac{\wt \theta_{1} (K)}{4 c_5}
  \leq \frac{2 q (1/2, K)}{4 c_5}
  \leq \frac{w (K)}{c_5}
  \, ,
\end{equation*}
where the bound $q (1/2, K) \leq 2 w (K)$ follows from Markov's inequality.
Therefore, Corollary~\ref{cor:link-sb-wills} implies that 
\begin{equation*}
  - \log \gamma_n \parens*{\frac{4 u}{\xi} K^\circ}
  \leq c_6 s_W \parens*{K, \frac{4 u}{c_6 \xi}}
  \, .
\end{equation*}
Putting things together, we obtain
\begin{equation*}
  2 j \log \xi
  < s_W \parens*{K, \frac{4 u}{c_6 \xi}}
  \, .
\end{equation*}
Now, from Lemma~\ref{lem:continuous-bij-lambda}, for any $u' > 0$ the map $j' \mapsto j'^{-1} \log W (\frac{j'}{u'} K)$ is decreasing.
Thus the above inequality together with the definition~\eqref{eq:def-eff-dim-wills} of $s_W$ imply that, for any $\xi > 1$, 
\begin{equation*}
  \log W \parens*{\frac{2 j \log \xi}{4 u/\parens{c_6 \xi}} K}
  = \log W \parens*{\frac{c_6 j \xi \log \xi}{2 u} K}
  \geq 2 j \log \xi
  \, .
\end{equation*}
Now, for any $\zeta > 1$, let $\xi = \sqrt{\zeta} > 1$, so that $\zeta = \xi^2 \geq \xi \log \xi$ and $\log \xi = \log (\zeta)/2$.
The above inequality then yields, for any $\zeta > 1$,
\begin{equation*}
  \log W \parens*{\frac{c_6 \zeta j}{2 u}}
  \geq j \log \zeta
  \, .
\end{equation*}
Thus, by Lemma~\ref{lem:fixed-log-mgf-theta}, we have
\begin{equation*}
  \theta_j (K)
  \geq \frac{2 u}{e c_6}
  = \frac{\wt \theta_{4 j} (K)}{8 e c_5 c_6}
  \, .
  \qedhere
\end{equation*}
% \end{proof}

\subsection{Extreme cases of Theorem~\ref{thm:intrinsic-quantile-exact}}
\label{sec:extr-cases-theor}

As a sanity check, let us now verify the relationship between the intrinsic widths $\theta_k (K)$ and normalized quantiles $\wt \theta_k (K)$ in the two extreme cases $k = 1$ and $k = n$.

\begin{example}
  \label{ex:example-quantiles-intrinsic-1}
  There is a constant $c \geq 1$ such that for every symmetric convex body $K$,
  \begin{equation}
    \label{eq:example-quantiles-intrinsic-1}
    \frac{w (K)}{c}
    \leq \wt \theta_1 (K)
    \leq c w (K)
    \, .
  \end{equation}
\end{example}

\begin{proof}
  On the one hand, for any $j > 0$ we have (using Markov's inequality)
  \begin{equation*}
    \wt \theta_j (K)
    \leq \parens*{\frac{1}{2}}^{-1/j} q \parens*{\frac{1}{2}, K}
    \leq 2^{1+1/j}\, w (K)
    \, .
  \end{equation*}

  On the other hand, it follows from~\cite[Theorem~1]{latala1999gaussian} that for some absolute constant $c_0 \geq 1$, for every $\eps \in \ocint{0}{1/2}$ we have
  \begin{equation*}
    q \parens*{\eps, K}
    \geq \frac{\eps w (K)}{c_0}
    \, .
  \end{equation*}
  Therefore,
  \begin{align*}
    \inf_{\eps \in \ooint{0}{1/2}} \eps^{-1} q (\eps, K)
    &\geq \frac{w (K)}{c_0}
      \, ,
      % \\
    \quad \text{ and } \quad
    \inf_{\eps \in \coint{1/2}{1}} \eps^{-1} q (\eps, K)
    \geq q (1/2, K)
      \geq \frac{w (K)}{2 c_0}
      \, .
  \end{align*}
  Hence $\wt \theta_1 (K) \geq w (K)/(2 c_0)$, and thus $\wt \theta_j (K) \geq w (K)/(2c_0)$ for every $j \in \ocint{0}{1}$.
  (We note in passing that using Lemma~\ref{lem:mini-sb} instead of the stronger result of Lata{\l}a and Oleszkiewicz~\cite{latala1999gaussian} would have given $\wt \theta_j (K) \geq w (K)/C_1(j)$ for $j \in (0, 1)$, for some function $C_1 : \ooint{0}{1} \to \R^+$.)
\end{proof}

We now turn to the opposite extreme, namely $k = n$.

\begin{example}
  \label{ex:example-quantiles-intrinsic-n}
  There exists a universal constant $c \geq 1$ such that, for every $n \geq 1$ and symmetric convex body $K \subset \R^n$ with non-empty interior,
  \begin{equation}
    \label{eq:ex-quantiles-intrinsic-n}
    c^{-1} n \vol_n (K)^{1/n}
    \leq \wt \theta_n (K)
    \leq c n \vol_n (K)^{1/n}
    \, .
  \end{equation}
  Hence, $\wt \theta_n (K) \asymp (n! \vol_n (K))^{1/n}/\sqrt{2\pi} = \theta_n (K)$.
\end{example}

\begin{proof}
  Indeed, since the density of $\gamma_n$ is upper-bounded by $(2\pi)^{-n/2}$ on $\R^n$, we have for every $u> 0$:
  \begin{equation*}
    \gamma_n (u K^\circ)
    \leq (2\pi)^{-n/2} \vol_n (u K^\circ)
    = (2\pi)^{-n/2} \vol_n (K^\circ) u^n
    \, .
  \end{equation*}
  On the other hand, since $K$ has non-empty interior, $K^\circ$ is bounded.
  Now since the density of $\gamma_n$ is lower-bounded by $(2\pi)^{-n/2} \exp (- u^2 \rad (K^\circ)^2/2)$ on $u K^\circ$, we have
  \begin{equation*}
    \gamma_n (u K^\circ)
    \geq (2\pi)^{-n/2} \vol_n (K^\circ) u^n \exp \parens*{- \frac{u^2 \rad (K^\circ)}{2}}
    \, .
  \end{equation*}
  Combining the previous two inequalities shows that
  \begin{equation*}
    \sup_{u > 0} \frac{\gamma_n (u K^\circ)}{u^n}
    = \lim_{u \to 0^+} \frac{\gamma_n (u K^\circ)}{u^n}
    = (2\pi)^{-n/2} \vol_n (K^\circ)
    \, .
  \end{equation*}
  Therefore from~\eqref{eq:def-scaled-quantile}, we have
  \begin{equation}
    \label{eq:ex-scaled-quantile-volume}
    \wt \theta_n (K)
    = \sqrt{2\pi} \vol_n (K^\circ)^{-1/n}
    \, .
  \end{equation}
  Now, it follows from the Blaschke-Santal{\'o} inequality~\cite{blaschke1917uber,santalo1949invariante} (for the upper bound) and from the converse from Bourgain and Milman~\cite{bourgain1987new} (for the lower bound) that for some universal constant $c_0 \geq 1$,
  \begin{equation}
    \label{eq:santalo-and-reverse}
    \frac{1}{c_0 n}
    \leq \vol_n (K)^{1/n} \vol_n (K^\circ)^{1/n}
    \leq \frac{c_0}{n}
    % \, ,
  \end{equation}
  (see \eg \cite[Corollary~7.2 p.~100]{pisier1999volume}).
  Combining~\eqref{eq:ex-scaled-quantile-volume} and~\eqref{eq:santalo-and-reverse} yields~\eqref{eq:ex-quantiles-intrinsic-n}.
\end{proof}

\subsection{Proof of Theorem~\ref{thm:intrinsic-quantile-exact}}
\label{sec:proof-intrinsic-quantiles-exact}

In this section, we provide the proof of the two-sided correspondence of Theorem~\ref{thm:intrinsic-quantile-exact}.

As suggested by the example of the volume (Example~\ref{ex:example-quantiles-intrinsic-n}), one direction will involve some form of the Blaschke-Santaló inequality, while the proof of the lower bound will involve the reverse inequality due to Bourgain-Milman.

\paragraph{Proof of the upper bound.}

% We start with the proof of the upper bound.

For the upper bound, we will rely on the following inequality, established by Alonso-Gutiérrez, Hernández Cifre, and Yepes Nicolás~\cite{alonso2021further}
% inside the proof of their Theorem~1.2 (specifically, in their new proof of the upper bound~\eqref{eq:wills-width-mcmullen}),
(specifically, this inequality can be extracted from the proof of the upper bound of their Theorem~1.2)
by combining the integral form of the Wills functional with the functional version of the Blaschke-Santaló inequality~\cite{klartag2005geometry}.

\begin{lemma}[\cite{alonso2021further}]
  \label{lem:ag-wills-exp}
  For every convex body $K \subset \R^n$, if $X \sim \gaussdist (0, I_n)$ then
  \begin{equation}
    \label{eq:ag-wills-exp}
    W (K) \Expect*{\exp \parens*{- \sup_{t \in K} \innerp{t}{X}}}
    \leq 1
    \, .
  \end{equation}
\end{lemma}

We now apply Lemma~\ref{lem:ag-wills-exp} to $\lambda K$, where $\lambda > 0$ is arbitrary.
Then, for any $u > 0$, we get
\begin{align*}
  1
  \geq W \parens*{\lambda K} \Expect*{\exp \parens*{- \lambda \sup_{t \in K} \innerp{t}{X}}}
  &\geq \lambda^j V_j \parens[\big]{K/\sqrt{2\pi}} \Probab*{\sup_{t \in K} \innerp{t}{X} \leq u} \exp \parens*{- \lambda u} \\
  &= \lambda^j \frac{\theta_j (K)^j}{j!} \gamma_n \parens*{u K^\circ} \exp \parens*{- \lambda u}
    \, .
\end{align*}
Upper-bounding $j! \leq j^j$, letting $\lambda = j/u$ and re-arranging, this yields
\begin{equation*}
  u \, \gamma_n \parens*{u K^\circ}^{-1/j}
  \geq e^{-1} \theta_j (K)
  \, .
\end{equation*}
Since this holds for every $u > 0$, taking the infimum yields $\wt \theta_j (K) \geq e^{-1} \theta_j (K)$.

\paragraph{Proof of the lower bound.}

The idea of the proof is to use the representation~\eqref{eq:intrinsic-tsirelson} of intrinsic volumes as expected volumes of Gaussian projections, together with the two facts below.

The first tool is the reverse Blaschke-Santaló inequality from Bourgain and Milman~\cite{bourgain1987new}.

\begin{lemma}[\cite{bourgain1987new}]
  \label{lem:bourgain-milman}
  There exists a universal constant $c_0 \geq 1$ such that for every $n \geq 1$ and symmetric convex body $K \subset \R^n$ with non-empty interior,
  \begin{equation}
    \label{eq:bourgain-milman}
    \vol_n (K) \vol_n (K^\circ)
    \geq c_0^{-n} \vol_n (B_2^n)^2
    \, .
  \end{equation}
\end{lemma}
In the remainder of the proof, we let $c_0$ denote the constant from Lemma~\ref{lem:bourgain-milman}, and $c_1 \geq 1$ a universal constant such that $(1/c_1\sqrt{n})^n \leq \vol_n (B_2^n) \leq (c_1/\sqrt{n})^n$ for all $n \geq 1$.

The second tool is the standard representation of the volume in terms of the gauge function.

\begin{fact}
  \label{fac:volume-gauge}
  For any $n \geq 1$ and symmetric convex body $K$ with non-empty interior, denoting by $\sigma_n$ the uniform (Haar) probability measure on the unit sphere $S^{n-1}$,
  \begin{equation}
    \label{eq:volume-gauge}
    \vol_n (K^\circ)
    = \vol_n (B_2^n) \int_{S^{n-1}} \parens*{\sup_{t \in K} \innerp{t}{y}}^{-n} \sigma_n (\di y)
    \, .
  \end{equation}
\end{fact}

We start with the sketch of a direct tentative approach, which does not yield the desired result but provides some intuition and motivation for the refined approach.
Combining Tsirelson's representation~\eqref{eq:intrinsic-tsirelson} with Lemma~\ref{lem:bourgain-milman} and Fact~\ref{fac:volume-gauge} applied to the Gaussian image $G_j K \subset \R^j$ (which has nonempty interior almost surely whenever the linear span of $K$ is at least $j$-dimensional), and further simplifying, we obtain that
\begin{equation*}
  \theta_j (K)
  \asymp \Expect*{\braces[\bigg]{\int_{S^{n-1}} \parens*{\sup_{t \in K} \innerp{G_j t}{y}}^{-j} \sigma_n (\di y)}^{-1} }^{1/j}
  \, .
\end{equation*}
Next, apply Jensen's inequality $\Expect*{Z^{-1}} \geq \Expect*{Z}^{-1}$ for any positive random variable $Z$, followed by Fubini's equality to swap the expectation and integration over $S^{n-1}$.
Since $\innerp{G_j t}{y} = \innerp{t}{G_j^* y}$ and since for any $y \in S^{n-1}$, one has $G_j^* y = y_1 X_1 + \dots + y_j X_j \sim \gaussdist (0, I_n)$, we deduce that, letting $X \sim \gaussdist (0, I_n)$,
\begin{equation*}
  \theta_j (K)
  \gtrsim \Expect[\bigg]{\parens*{\sup_{t \in K} \innerp{t}{X}}^{-j}}^{-1/j}
  \, .
\end{equation*}
Unfortunately, the quantity in the lower bound is not equivalent to $\wt \theta_j (K)$, and this lower bound is not sharp in general.
For instance, when $j = n$, the right-hand side equals $0$ while (by Example~\ref{ex:example-quantiles-intrinsic-n}) $\wt \theta_n (K) \asymp \theta_n (K) \asymp n \vol_n (K)^{1/n}$.

In order to refine the lower bound above, we will apply the combination of Lemma~\ref{lem:bourgain-milman} and Fact~\ref{fac:volume-gauge} to a different body than $G_j K$.
% Specifically, we in order to address integrability issues associated with negative moments, we will work with a 
First, in order to address some integrability issues associated with negative moments, we will work with a higher-dimensional Gaussian projection of $K$.
Second, in order to gain more flexibility in the optimization of the estimate, we will regularize the projection via the addition a Euclidean ball.

Specifically, let $k = j+1$ and $\rho > 0$ a positive real number, to be chosen later.
Let $\ol K_\rho = G_k K + \rho B_2^k \subset \R^k$.
We will obtain a lower and an upper estimate on $\Expect{\vol_k (\ol K_\rho)}$.

\emph{Lower estimate.}
By Lemma~\ref{lem:bourgain-milman} followed by Jensen's inequality, we have
\begin{align}
  \Expect*{\vol_k \parens*{\ol K_\rho}}
  &\geq \Expect*{c_0^{-k} \vol_k (B_2^k)^2 \vol_k \parens*{\parens[\big]{G_k K + \rho B_2^k}^\circ}^{-1}} \nonumber \\
  &\geq c_0^{-k} \vol_k (B_2^k)^2 \Expect*{\vol_k \parens*{\parens[\big]{G_k K + \rho B_2^k}^\circ}}^{-1}
    \label{eq:proof-exact-quantile-lower-bm}
    \, .
\end{align}
Now, by Fact~\ref{fac:volume-gauge}, we have
\begin{align*}
  \vol_k \parens*{\parens[\big]{G_k K + \rho B_2^k}^\circ}
  &= \vol_k (B_2^k) \int_{S^{n-1}} \parens*{\sup_{\tau \in G_k K + \rho B_2^k} \innerp{\tau}{y}}^{-k} \sigma_k (\di y) \\
  % &= \vol_k (B_2^k) \int_{S^{n-1}} \parens*{\sup_{t \in K} \innerp{G_k t}{y} + \rho}^{-k} \sigma_k (\di y) \\
  &= \vol_k (B_2^k) \int_{S^{n-1}} \parens*{\sup_{t \in K} \innerp{t}{G_k^* y} + \rho}^{-k} \sigma_k (\di y)
    \, .
\end{align*}
Taking expectations, using that $G_k^* y \sim \gaussdist (0, I_n)$ for any $y \in S^{n-1}$, and plugging into~\eqref{eq:proof-exact-quantile-lower-bm} gives
\begin{equation}
  \label{eq:proof-exact-quantile-expect-1}
  \Expect*{\vol_k \parens*{\ol K_\rho}}
  \geq c_0^{-k} \vol_k (B_2^k) \Expect*{\parens*{\sup_{t \in K} \innerp{t}{X} + \rho}^{-k}}^{-1}
  \, .
\end{equation}
Now, using that $(h + \rho)^{-k} = \int_{h}^{\infty} k (u + \rho)^{-k-1} \di u = \int_{0}^{\infty} k \indic{h \leq u} (u + \rho)^{-k-1} \di u$ for $h \in \R^+$, we get
\begin{align*}
  \Expect*{\parens*{\sup_{t \in K} \innerp{t}{X} + \rho}^{-k}}
  &= \Expect*{\int_0^\infty k \, \indic*{\sup_{t \in K} \innerp{t}{X} \leq u} \parens*{u + \rho}^{-k-1} } \\
  &= \int_0^\infty k \Probab*{\sup_{t \in K} \innerp{t}{X} \leq u} \parens*{u + \rho}^{-k-1}
    \, .
\end{align*}
Now, by definition~\eqref{eq:def-scaled-quantile} of $\wt \theta_j (K)$, for every $u > 0$ we have
\begin{equation*}
  \Probab*{\sup_{t \in K} \innerp{t}{X} \leq u}
  = \gamma_n (u K^\circ)
  \leq \parens[\bigg]{\frac{u}{\wt \theta_j (K)}}^j
  \, .
\end{equation*}
%In addition, for the rest of this proof we choose $k = j+1$.
Combining the above two inequalities yields
\begin{align}
  \Expect*{\parens*{\sup_{t \in K} \innerp{t}{X} + \rho}^{-k}}
  &\leq \wt \theta_j (K)^{-j} \int_0^\infty \frac{(j+1) u^j}{(u + \rho)^{j+2}} \di u \nonumber \\
  &= \wt \theta_j (K)^{-j} \bracks*{\frac{1}{\rho} \parens*{\frac{u}{u+\rho}}^{j+1}}_0^{\infty} 
    = \frac{\wt \theta_j (K)^{-j}}{\rho}
    \label{eq:proof-exact-quantile-expect-neg-power}
    \, .
\end{align}
Finally, plugging~\eqref{eq:proof-exact-quantile-expect-neg-power} into~\eqref{eq:proof-exact-quantile-expect-1} yields the lower estimate
\begin{equation}
  \label{eq:proof-exact-quantile-expect-lower}
  \Expect*{\vol_{j+1} \parens*{\ol K_\rho}}
  \geq c_0^{-j-1} \vol_{j+1} (B_2^{j+1}) \rho \, \wt \theta_j (K)^j
  \, .
\end{equation}

\emph{Upper estimate.}
On the other hand, Steiner's formula~\eqref{eq:steiner} yields
\begin{equation}
  \label{eq:proof-quantile-exact-up-1}
  \Expect*{\vol_{k} (\ol K_\rho)}
  = \Expect*{\vol_{k} \parens*{G_{k} K + \rho B_2^{k}}} 
    % = \Expect*{\sum_{i=0}^{k} V_i (G_{k} K) \vol_{k-i} (B_2^{k-i}) \rho^{k-i}} \\
  = \sum_{i=0}^{k} \Expect*{V_i (G_{k} K)} \vol_{k-i} (B_2^{k-i}) \rho^{k-i}
    \, .
\end{equation}
Now, Kubota's formula (in the form given in~\cite[Definition~1.1]{lotz2020concentration}) shows that, letting
% $F_i$ denote a uniformly distributed random subspace of dimension $i$ of $\R^k$, independent of $G_k$, and $P_{F_i}$ the orthogonal projection on $F_i$,
$P_{i,k}$ denote the projection $\R^k \to \R^i$ obtained by keeping the first $i$ coordinates, and letting $U_k$ denote a random matrix (independent of $G_k$) distributed according to the uniform (Haar) measure on $k \times k$ rotation matrices, we have
\begin{align*}
  \Expect*{V_i (G_k K)}
  &= \Expect*{\binom{k}{i} \frac{\vol_k (B_2^k)}{\vol_i (B_2^i) \vol_{k-i} (B_2^{k-i})} \vol_i \parens*{P_{i,k} U_k G_kK} }
    \, .
\end{align*}
Since $U_k, G_j$ are independent and $U_k$ is orthogonal, $U_k G_k$ has the same distribution as $G_k$, thus $P_{i,k} U_k G_k$ has the same distribution as $P_{i,k} G_k = G_i$.
Using this in the above equation and then using Tsirel'son's formula~\eqref{eq:intrinsic-tsirelson} gives
\begin{align}
  \Expect*{V_i (G_k K)}
  &= \binom{k}{i} \frac{\vol_k (B_2^k)}{\vol_{k-i} (B_2^{k-i})} \frac{\Expect{\vol_i \parens*{G_i K}}}{\vol_i (B_2^i)} 
  = \frac{k!}{(k-i)!} \frac{\vol_k (B_2^k)}{\vol_{k-i} (B_2^{k-i})} V_i (K/\sqrt{2\pi})
    \label{eq:expect-intr-proj}
    \, .
\end{align}
Combining~\eqref{eq:proof-quantile-exact-up-1} and~\eqref{eq:expect-intr-proj} and bounding $k!/(k-i)! \leq k^i \leq (2 j)^i$ as well as $V_i (K/\sqrt{2\pi}) \geq 0$ for $i > k$, we obtain
\begin{align}
  \label{eq:proof-quantile-exact-up-wills}
  \Expect*{\vol_{j+1} (\ol K_\rho)}
  &= {\vol_k (B_2^k)} \sum_{i=0}^{k} \frac{k!}{(k-i)!} V_i (K/\sqrt{2\pi}) \rho^{k-i} \nonumber \\
  &\leq \vol_{j+1} (B_2^{j+1}) \rho^{j+1} \sum_{i=0}^n \parens*{\frac{2 j}{\rho}}^i V_i \parens{K/\sqrt{2\pi}} \nonumber \\
  &\leq \vol_{j+1} (B_2^{j+1}) \rho^{j+1} W \parens*{\frac{2 j}{\rho} K}
    \, .
\end{align}

\emph{Optimization over $\rho$ and conclusion.}
By combining the lower estimate~\eqref{eq:proof-exact-quantile-expect-lower} with the upper estimate~\eqref{eq:proof-quantile-exact-up-wills} and simplifying, we obtain that for every $\rho > 0$,
\begin{equation*}
  c_0^{-j-1} \, \wt \theta_j (K)^j
  \leq \rho^{j} W \parens*{\frac{2 j}{\rho} K}
  \, .
\end{equation*}
Thus, letting $\lambda = 2 j/\rho$ and using that $c_0^{j+1} \leq c_0^{2j}$,
\begin{equation*}
  \wt \theta_j (K)^j
  \leq c_0^{2j} (2 j)^j \inf_{\lambda > 0} \braces*{\lambda^{-j} W (\lambda K)}
  \, .
\end{equation*}
Next, by applying Lemma~\ref{lem:ulc-seq-from-mgf} to the ultra-log-concave sequence $(V_{j} (K/\sqrt{2\pi}))_{j \geq 0}$, we deduce from the above that
\begin{equation*}
  \wt \theta_j (K)^j
  \leq c_0^{2j} (2 j)^j \cdot (e/j)^j j! V_j (K/\sqrt{2\pi})
  = \parens*{2 e c_0^2}^j \theta_j (K)^j
  \, .
\end{equation*}
Hence $\wt \theta_j (K) \leq 2 e c_0^2 \theta_j (K)$, which concludes the proof.

\subsection{Basic properties of $s_W$}
\label{sec:addit-lemm-small}

We conclude by providing additional basic properties of the Wills functional and of the quantity $s_W$ defined in Corollary~\ref{cor:link-sb-wills}, which will help apply the latter result.
% While not strictly necessary to the proof of the results of Section~\ref{sec:intr-volum-wills}, these properties may help interpret or apply Corollary~\ref{cor:link-sb-wills}.

\begin{lemma}
  \label{lem:strict-concave-wills}
  Let $K$ be a convex body which is not reduced to a single point.
  Then, the map $\lambda \mapsto \log W (\lambda K)$ is strictly concave on $\R^+$.
  % In particular, the map $\lambda \mapsto \lambda^{-1} \log W (\lambda K)$ is strictly decreasing on $\R_+^*$.
\end{lemma}

\begin{proof}[Proof of Lemma~\ref{lem:strict-concave-wills}]
  It follows from~\cite[Theorem~1.4]{alonso2021further} that the map $\lambda \mapsto \log W (\lambda K)$ is concave on $\R^+$.
  Hence, if this function were not strictly concave, it would be affine on a segment of positive length of $\R^+$.
  But this function is also analytic (as the logarithm of a polynomial) on $\R^+$; hence, the above would imply that it is affine on $\R^+$.
  But since $\log W (\lambda K) = O_K (\log \lambda)$ as $\lambda \to + \infty$ and $\log W (0 K) = 0$, one would then have $\log W (\lambda K) = 0$ for every $\lambda > 0$, thus $w (K) = 0$, which is ruled out by the assumption that $K$ is not reduced to a single point.
\end{proof}

\begin{lemma}
  \label{lem:continuous-bij-lambda}
  Let $K$ be as in Lemma~\ref{lem:strict-concave-wills}.
  The map $\lambda \mapsto \lambda^{-1} \log W (\lambda K)$ is a continuous and decreasing bijection from $\R_+^*$ to $(0, w(K))$.
\end{lemma}

\begin{proof}[Proof of Lemma~\ref{lem:continuous-bij-lambda}]
  By Lemma~\ref{lem:strict-concave-wills}, the map $\lambda \mapsto \lambda^{-1} \log W (\lambda K)$ is strictly decreasing on $\R_+^*$, and it is additionally continuous.
  Also, its limit as $\lambda \to 0^+$ equals $V_1 (K/\sqrt{2\pi}) = w(K)$, and its limit as $\lambda \to + \infty$ equals $0$, thus it is a continuous decreasing bijection $\R_+^* \to (0, w (K))$.
\end{proof}

% Finally, we have the following property on the behavior of $s_W$.

\begin{lemma}
  \label{lem:bij-decr-lambda-u}
  Let $K$ be as in Lemma~\ref{lem:strict-concave-wills}.
  The map $u \mapsto s_W (K, u)$ (defined in Corollary~\ref{cor:link-sb-wills}) is a continuous and decreasing bijection from $(0, w(K))$ to $\R_+^*$.

  In addition, for $u \in \ooint{0}{w (K)}$ and $s > 0$: if $s < s_W (K, u)$ then $\log W (\frac{s}{u} K) > s$, if $s > s_W (K, u)$ then $\log W (\frac{s}{u} K) < s$, and if $s = s_W (K, u)$ then $\log (\frac{s}{u} K) = s$.  
  % In addition, if $u \in \ooint{0}{w (K)}$ and $s > 0$, then $\log W (\frac{s}{u} K) > s$ if $s < s_W (K, u)$, while $\log W (\frac{s}{u} K) < s$ if $s > s_W (K, u)$ and $\log (\frac{s}{u} K) = s$ if $s = s_W (K, u)$.
\end{lemma}

\begin{proof}[Proof of Lemma~\ref{lem:bij-decr-lambda-u}]
  Recall from~\eqref{eq:def-eff-dim-wills}
  that, for $u \in \ooint{0}{w (K)}$,
  \begin{equation}
    \label{eq:lambda-t-repres}
    s_W (K, u)
    = \sup \braces*{s \in \R^+ : \log W \parens*{\frac{s}{u} K} \geq s}
    \, .
  \end{equation}
  By Lemma~\ref{lem:continuous-bij-lambda}, $s_W (K, u) > 0$.
  Also, the map $s \mapsto s^{-1}\log W (\frac{s}{u} K)$ is continuous and strictly decreasing on $\R_+^*$ (by Lemma~\ref{lem:strict-concave-wills}), so for $s > 0$ one has $\log W (\frac{s}{u} K) > s$ (resp.~$\log W (\frac{s}{u} K) < s$; $\log W (\frac{s}{u} K) = s$) if and only if $s < s_W (K, u)$ (resp.~$s > s_W (K, u)$; $s = s_W (K, u)$).
  
  Now, since $V_j (K/\sqrt{2\pi}) \geq 0$ for $1 \leq j \leq n$ and $V_1 (K/\sqrt{2\pi}) > 0$, for $s >0$, the map $u \mapsto \log W (\frac{s}{u} K)$ is strictly decreasing.
  When combined with the conclusion of the previous paragraph, this entails that $u \mapsto s_W (K, u)$ is strictly decreasing.
  In addition, by Lemma~\ref{lem:continuous-bij-lambda}, this map is continuous and $\lim_{u \to w(K)^{-}} s_W (K, u) = 0$.

  It thus remains to show that $\sup_{u \in (0, w(K))} s_W (K, u) = + \infty$.
  Fix $s > 0$.
  We have $\log W (\frac{s}{u} K) \geq \log \parens*{1 + \frac{s}{u} w (K)} \to + \infty$ as $u \to 0^+$, thus there exists $u \in (0, w(K))$ such that $\log W (\frac{s}{u} K) \geq s$, and thus (by~\eqref{eq:lambda-t-repres}) $s_W (K, u) \geq s$.
  Since $s > 0$ was arbitrary, this concludes the proof.  
\end{proof}

\section{Metric characterization of quantiles and comparison inequalities}
\label{sec:chain-char-quant}

In this section, we prove the results of Section~\ref{sec:chaining-comparison} concerning the relationship between small ball probabilities and the metric structure of the convexified index set (Corollary~\ref{cor:quantiles-chaining}, Corollary~\ref{cor:comparison-small-ball} and Proposition~\ref{prop:no-comparison-nonconvex}).

\subsection{Proof of Corollary~\ref{cor:quantiles-chaining}
  % Chaining characterization of lower quantiles
}
\label{sec:chain-char-lower}

The proof is based on the correspondence of Corollary~\ref{cor:link-sb-wills} between small ball probabilities and the Wills functional, together with the following metric characterization of the Wills functional from~\cite{mourtada2025universal}.

\begin{lemma}[\cite{mourtada2025universal}, Corollary~4.1]
  \label{lem:mourt-wills-chaining}
  There exists a universal constant $c>0$ such that the following holds.
  % For every $n \geq 1$ and every convex body $K \subset \R^n$,
  % \begin{equation}
  %   \label{eq:wills-chaining}
  %   \frac{1}{c} \inf_{p \in \N} \braces*{(2^p - 1) + \gamma_2^{(p)} (K)}
  %   \leq \log W (K)
  %   \leq c \inf_{p \in \N} \braces*{(2^p - 1) + \gamma_2^{(p)} (K)}
  %   \, .
  % \end{equation}
  For every $n \geq 1$ and every convex body $K \subset \R^n$, define
  \begin{equation}
    \label{eq:def-p-star-trunc-chaining}
    p^* (K)
    = \sup \braces*{p \in \N : \gamma_2^{(p)} (K) \geq 2^p }
    \in \N \cup \set{\pm\infty}    
    \, ,
  \end{equation}
  with the convention that $p^* (K) = - \infty$ if the set is empty, namely if $\gamma_2 (K) < 1$.
  In addition, let $\gamma_* (K) = \gamma_2 (K)$ if $p^* (K) = -\infty$, and $\gamma_* (K) = 2^{p^* (K)}$ otherwise.
  Then:
  \begin{equation}
    \label{eq:wills-chaining}
    \frac{1}{c} \gamma_* (K)
    \leq \log W (K)
    \leq c \gamma_* (K)
    \, .
  \end{equation}
\end{lemma}

We will also use the following lemma, which follows from~\cite[Theorem~5.2]{mourtada2025universal}.
\begin{lemma}
  \label{lem:wills-difference-body}
  For any convex body $K \subset \R^n$, one has
  $\log W (K) \leq \log W (K - K) \leq 2 \log W (K)$.
\end{lemma}

\begin{proof}[Proof of Lemma~\ref{lem:wills-difference-body}]
  The lower bound is due to the fact that the Wills functional is translation-invariant and monotone under inclusion (\eg, using~\eqref{eq:intrinsic-tsirelson}).
  The upper bound follows from the inequality $\log W (K + K') \leq \log W (K) + \log W (K')$ for any convex body $K' \subset \R^n$ (\cite[Theorem~5.2]{mourtada2025universal}) and the fact that $W (-K) = W (K)$.
\end{proof}

We now turn to the proof.
First, if $K$ is reduced to a single point, then the inequalities~\eqref{eq:quantiles-chaining} trivially hold as all quantities equal $0$.
We thus assume that $K$ is not reduced to a point.
Thus, $K-K$ is a nontrivial symmetric convex body, with $w (K-K) = 2 w (K)$.

By Corollary~\ref{cor:link-sb-wills}, for some universal constants $c_0, c_1 > 1$, for any $u \in (0, w(K)/c_0)$,
\begin{equation}
  \label{eq:proof-quantiles-sb-lambda}
  \frac{1}{c_1} s_W (K - K, c_1 u)
  \leq - \log \gamma_n (u (K-K)^\circ)
  \leq c_1 s_W (K-K, u/c_1)
  \, .
\end{equation}
Also, let $c_2 > 1$ denote the constant $c$ from the Majorizing Measures theorem (Theorem~\ref{thm:talagrand-maj-meas}) and $c_3 > 1$ the constant from Lemma~\ref{lem:mourt-wills-chaining}, so that for every $\lambda \in \R^+$,
\begin{align}
  \label{eq:proof-quantiles-maj-meas}
  \frac{1}{c_2} \gamma_2 (K)
  &\leq w (K)
    \leq c_2 \gamma_2 (K)
    \, ; \\
  \frac{1}{c_3} \gamma_* (\lambda K)
  &\leq \log W (\lambda K)
    \leq c_3 \gamma_* (\lambda K)
    \label{eq:proof-quantiles-wills-chaining}
    \, .
\end{align}

We start with the upper bound in~\eqref{eq:quantiles-chaining}.
Let $c > 1$ to be chosen later and $u = \gamma_2^{(p)} (K)/c$ (note that $\gamma_2^{(p)} (K) > 0$ since $K$ is infinite).
We want to show that
\begin{equation*}
  u
  \leq Q \parens*{2^p/c, \sup_{t, t' \in K} \abs*{\innerp{t - t'}{X}}}
  \, ,
\end{equation*}
which amounts to saying that (using~\eqref{eq:low-tail-gauss-polar})
\begin{equation*}
  \frac{2^p}{c}
  \leq - \log \Probab*{\sup_{t, t' \in K} \abs*{\innerp*{t-t'}{X}} \leq u}
  = - \log \gamma_n \parens*{u (K - K)^\circ}
  \, .
\end{equation*}
Now, from the lower bound in~\eqref{eq:proof-quantiles-sb-lambda}, it suffices to show that
\begin{equation*}
  \frac{c_1 2^p}{c}
  \leq s_W \parens*{K-K, c_1 u}
  \, .
\end{equation*}
By Lemma~\ref{lem:bij-decr-lambda-u}, this holds provided that
\begin{equation*}
  \log W \parens*{\frac{c_1 2^p/c}{c_1 u} (K-K)}
  = \log W \parens*{\frac{2^p}{c u} (K-K)}
  \geq \frac{c_1 2^p}{c}
  \, .
\end{equation*}
By the lower bound in Lemma~\ref{lem:wills-difference-body}, it suffices that $\log W (\frac{2^p}{c u} K) \geq c_1 2^p/c$, and by the lower bound in~\eqref{eq:proof-quantiles-wills-chaining}, it suffices that
\begin{equation*}
  \gamma_* \parens*{\frac{2^p}{c u} K}
  \geq \frac{c_1 c_3 2^p}{c}
  \, .
\end{equation*}
Let $c = c_1 c_3$.
Since $c u = \gamma_2^{(p)} (K)$, the condition above amounts to
\begin{equation*}
  \gamma_* \parens[\bigg]{\frac{2^p}{\gamma_2^{(p)} (K)} K}
  \geq 2^p
  \, ,
\end{equation*}
which indeed holds by definition of $\gamma_* (K)$ and since
\begin{equation*}
  \gamma_2^{(p)} \parens[\bigg]{\frac{2^p}{\gamma_2^{(p)} (K)} K}
  = \frac{2^p}{\gamma_2^{(p)} (K)} \gamma_2^{(p)} (K)
  = 2^p
  \, .
\end{equation*}

We now turn to the lower bound in~\eqref{eq:quantiles-chaining}.
Let $u = 2 c_1 c_3 \gamma_2^{(p)} (K)$.
We will show that, for $c' = 8 c_0 c_1 c_3 \geq 2$,
\begin{equation}
 \label{eq:proof-quantiles-chaining-upper-todo}
  Q \parens*{c' 2^p, \sup_{t, t' \in K} \abs*{\innerp*{t-t'}{X}}}
  \leq c' \gamma_2^{(p)} (K)
  \, .
\end{equation}

Assume first that $u \geq w (K)/c_0$.
Since $c' \geq 2$, we have $\exp (- c' 2^p) \leq e^{-2} \leq 1/2$.
Markov's inequality then implies that~\eqref{eq:proof-quantiles-chaining-upper-todo} holds, as
\begin{equation*}
  Q \parens*{c' 2^p, \sup_{t, t' \in K} \abs*{\innerp*{t-t'}{X}}}
  \leq 2 \Expect*{\sup_{t, t' \in K} \abs*{\innerp*{t-t'}{X}}}
  = 4 w (K)
  \leq 4 c_0 u
  = c' \gamma_2^{(p)} (K)
  \, .
\end{equation*}

Assume now that $u < w (K)/c_0$, and let $\lambda = 2^p/\gamma_2^{(p)} (K)$.
For any $q \geq p+1$, we have $\gamma_2^{(q)} (\lambda K) = \frac{2^p}{\gamma_2^{(p)} (K)} \gamma_2^{(q)} (K) \leq 2^{p} < 2^{q}$, thus $p^*(\lambda K) \leq p$ and $\gamma_* (\lambda K) \leq 2^{p}$ (when $p^*(\lambda K) = -\infty$, we indeed have $\gamma_* (\lambda K) = \gamma_2 (\lambda K) < 1 \leq 2^p$).
Using~\eqref{eq:proof-quantiles-wills-chaining} and the upper bound in Lemma~\ref{lem:wills-difference-body}, this gives
\begin{align*}
  \log W \parens*{\frac{2 c_3 2^p}{2 c_3 \gamma_2^{(p)} (K)} (K-K)}
  = \log W (\lambda (K-K))
  \leq 2 \log W (\lambda K)
  \leq 2 c_3 \gamma_* (\lambda K)
  \leq 2 c_3 2^p
  \, .
\end{align*}
By Lemma~\ref{lem:bij-decr-lambda-u}, this shows that
\begin{equation*}
  2 c_3 2^p
  \geq s_W \parens*{K-K, 2 c_3 \gamma_2^{(p)} (K)}
  = s_W \parens*{K - K, u/c_1}
  \, .
\end{equation*}
Together with~\eqref{eq:low-tail-gauss-polar} and~\eqref{eq:proof-quantiles-sb-lambda}, this gives
\begin{equation*}
  - \log \Probab*{\sup_{t, t' \in K} \abs*{\innerp*{t-t'}{X}} \leq u}
  = - \log \gamma_n \parens*{u (K-K)^\circ}
  \leq c_1 s_W \parens*{K - K, u/c_1}
  \leq 2 c_1 c_3 2^p
  \, ,
\end{equation*}
which means that
\begin{equation*}
  Q \parens*{2 c_1 c_3 2^p, \sup_{t, t' \in K} \abs*{\innerp*{t-t'}{X}}}
  \leq u
  = 2 c_1 c_3 \gamma_2^{(p)} (K)
  \, .
\end{equation*}
Since $2 c_1 c_3 \leq c'$, this implies that~\eqref{eq:proof-quantiles-chaining-upper-todo} also holds in this case.

Combining the upper and lower bounds, and using monotonicity in $c$, shows that Corollary~\ref{cor:quantiles-chaining} holds with $c = \max \set{c_1 c_3, c'} = 8 c_0 c_1 c_3$.

\subsection{Proof of Corollary~\ref{cor:comparison-small-ball}%
  %Comparison inequality for convex sets
}
\label{sec:comp-ineq-conv}

We note that the main part of the proof, namely the bound for the lower tail, could have been established by combining Corollary~\ref{cor:link-sb-wills} with a comparison theorem for the Wills functional~\cite[Theorem~4.1]{mourtada2025universal}.
However, we prefer to deduce it from the chaining characterization of small quantiles of Corollary~\ref{cor:quantiles-chaining}, as truncated generic chaining functionals are more transparently metric quantities than the Wills functional.

%\begin{proof}[Proof of Corollary~\ref{cor:comparison-small-ball}]
We will prove the inequality between quantiles:
for every $s \in \R_+^*$,
\begin{equation}
  \label{eq:comparison-quantiles}
  Q \parens[\bigg]{c s, \sup_{t,t' \in K_1} \abs{\innerp{t-t'}{X}}}
  \leq  c \, Q \parens[\bigg]{s, \sup_{t,t' \in K_2} \abs{\innerp{t-t'}{X}}}
  % \, ,
\end{equation}
whenever $K_1, K_2$ are nontrivial (not reduced to single points),
which will convert to the desired inequality between cumulative distribution functions.

First, we claim that $\gamma_2^{(p)} (K_1) \leq \gamma_2^{(p)} (K_2)$ for every $p \in \N$.

Indeed, let $(\A_j)_{j \geq 0}$ be an admissible sequence of partitions of $K_2$.
For every $j \geq 0$, define $\B_j = \set{\phi^{-1} (A) : A \in \A_j}$, so that $(\B_j)_{j \geq 0}$ is an admissible sequence of partitions of $K_1$.
Let $t_1 \in K_1$ be arbitrary, and set $t_2 = \phi (t_1) \in K_2$.
For any $j \geq 0$, since $\phi (t_1) \in A_j (t_2)$, we have $B_j (t_1) = \phi^{-1} (A_j (t_2))$, thus $\phi (B_j (t_1)) \subset A_j (t_2)$ and thus (given that $\phi$ is expansive) $\diam B_j (t_1) \leq \diam \phi (B_j (t_1)) \leq \diam A_j (t_2)$.
Since this holds for every $t_1 \in K_1$ and $j \in \N$, we deduce that
\begin{equation*}
  \gamma_2^{(p)} (K_1)
  \leq \sup_{t_1 \in K_1} \sum_{j \geq p} 2^{j/2} \diam B_j (t_1)
  \leq \sup_{t_2 \in K_2} \sum_{j \geq p} 2^{j/2} \diam A_j (t_2)
  \, .
\end{equation*}
Since this holds for every admissible sequence $(\A_j)_{j \geq 0}$, we conclude that $\gamma_2^{(p)} (K_1) \leq \gamma_2^{(p)} (K_2)$.

Let $c_0 \geq 1$ denote the constant from the Majorizing Measures theorem (Theorem~\ref{thm:talagrand-maj-meas}), and $c_1 \geq 1$ the constant from Corollary~\ref{cor:quantiles-chaining}.
Also, from Corollary~\ref{cor:upper-quantiles-gaussian} and Lemma~\ref{lem:quantiles-constant} (together with the trivial bound $\diam(T) \leq \gamma_2 (T)$ and the fact that $Q (s, \cdot) = Q_+ (s', \cdot)$ with $s' = - \log(1-e^{-s}) \asymp_{c_1} \log(e \vee 1/s)$ for $s \leq c_1$), there exists a constant $c_2 \geq 1$ such that for any $n \geq 1$, compact set $T \subset \R^n$, and $s \leq c_1$,
\begin{equation}
  \label{eq:proof-quantiles-up}
  \frac{1}{c_2} \braces*{\gamma_2 (T) + \diam (T) \sqrt{\log \parens*{e \vee \frac{1}{s}}}}
  \leq Q \parens[\bigg]{s, \sup_{t,t' \in T} \abs{\innerp{t-t'}{X}}}
  \leq c_2 \braces*{\gamma_2 (T) + \diam (T) \sqrt{\log \parens*{e \vee \frac{1}{s}}}}
  .
\end{equation}

Let $s > c_1$, and $p \in \N$ such that $c_1 2^{p} < s \leq 2 c_1 2^{p}$.
We have from Corollary~\ref{cor:quantiles-chaining} and the above:
\begin{align*}
  Q \parens[\bigg]{s, \sup_{t,t' \in K_1} \abs{\innerp{t-t'}{X}}}
  &\leq Q \parens[\bigg]{c_1 2^p, \sup_{t,t' \in K_1} \abs{\innerp{t-t'}{X}}}
    \leq c_1 \gamma_2^{(p)} (K_1)
    \leq c_1 \gamma_2^{(p)} (K_2) \\
  &\leq c_1^2 Q \parens[\bigg]{2^p/c_1, \sup_{t,t' \in K_2} \abs{\innerp{t-t'}{X}}}
    \leq c_1^2 Q \parens[\bigg]{\frac{s}{2 c_1^2}, \sup_{t,t' \in K_2} \abs{\innerp{t-t'}{X}}}
    \, .
\end{align*}

On the other hand, if $s \leq c_1$, from~\eqref{eq:proof-quantiles-up},
\begin{align*}
  Q \parens[\bigg]{s, \sup_{t,t' \in K_1} \abs{\innerp{t-t'}{K_1}}}
  &\leq c_2 \braces*{\gamma_2 (K_1) + \diam (K_1) \sqrt{\log \parens*{e \vee \frac{1}{s}}}} \\
  &\leq c_2 \braces*{\gamma_2 (K_2) + \diam (K_2) \sqrt{\log \parens*{e \vee \frac{1}{s}}}} 
    \leq c_2^2 Q \parens*{s, \sup_{t \in K_2} \innerp{t}{X}}
    \, .
\end{align*}

The above two inequalities imply that~\eqref{eq:comparison-quantiles} holds with $c = \max\braces*{2 c_1^2, c_2^2}$ for any $s > 0$.

We conclude by deducing the formulation~\eqref{eq:comparison-cdf} in terms of cumulative distribution functions.
If $K_1$ contains a single point, the inequality is trivial.
Assume now $K_1$ is nontrivial; then, the existence of $\phi$ implies that $K_2$ is nontrivial as well.
In this case, the inequality is trivial if $u = 0$, since both sides of the inequality are equal to $0$.
Assume now that $u > 0$.
We use that for $K \neq \set{0}$ a symmetric convex body in $\R^n$, the map $u \mapsto \Probab*{\sup_{t \in K} \innerp{t}{X} \leq u} = \gamma_n (u K^\circ)$ is a continuous and strictly increasing bijection from $\R_+^*$ to $(0, 1)$,
%%% small justification? $K^\circ$ has non-empty interior so strictly increasing, and $\partial K^\circ$ has $0$ Lebesgue measure and thus Gaussian measure
and that for $s > 0$ we have $\Probab*{\sup_{t \in K} \innerp{t}{X} \leq u} = \exp(-s)$ if and only if $u = Q (s, \sup_{t \in K} \innerp{t}{X})$.
Thus, letting $s = - \log \Probab*{\sup_{t \in K_2-K_2} \innerp{t}{X} \leq u} = - \log \Probab*{\sup_{t,t' \in K_2} \abs{\innerp{t-t'}{X}} \leq u}$, we have from~\eqref{eq:comparison-quantiles}:
\begin{align*}
  c u
  &= c \, Q \parens[\bigg]{s, \sup_{t,t' \in K_2} \abs{\innerp{t-t'}{X}}}
    \geq Q \parens[\bigg]{c s, \sup_{t,t' \in K_1} \abs{\innerp{t-t'}{X}}}
    \, ,
\end{align*}
and therefore
\begin{align*}
  \Probab[\bigg]{\sup_{t,t' \in K_1} \abs{\innerp{t-t'}{X}} \leq c u}
  \geq \exp (- c s)
  = \Probab[\bigg]{\sup_{t,t' \in K_2} \abs{\innerp{t-t'}{X}} \leq u}^c
  \, .
\end{align*}
This is precisely the claimed inequality.
%\end{proof}

\subsection{Proof of Proposition~\ref{prop:no-comparison-nonconvex}%
}
\label{sec:lack-comp-ineq}

% It is worth noting that the comparison inequality of Corollary~\ref{cor:comparison-small-ball} for general quantiles holds under the assumption of an expansive embedding between the \emph{convex} sets $K_1, K_2$.
% In contrast, the typical and large values of the supremum of a Gaussian process indexed by a (not necessarily convex) set $T \subset \R^n$ are governed by the quantities $\gamma_2 (T)$ and $\diam (T)$, thus satisfies a comparison inequality for general sets $T_1, T_2$.

% The following result shows that the assumption of convexity in Corollary~\ref{cor:comparison-small-ball} cannot be removed for the lower tails, as such a comparison principle fails to hold for general non-convex sets.
% This shows a critical difference between the upper and lower tails of suprema of Gaussian processes: while the upper tail of the supremum of a Gaussian process indexed by a general set $T$ may be characterized (in a monotonic way, via Corollary~\ref{cor:upper-quantiles-gaussian}) in terms of the metric structure of $T$, the lower tail can only be understood in terms of the metric structure of the \emph{convexified} set $K = \conv (T-T)$.

%\begin{proof}[Proof of Proposition~\ref{prop:no-comparison-nonconvex}]
  Let $T_{1,n} = \set{-1, 1}^n \subset \R^n$ and
  % $T_{2,n} = \set{2\sqrt{n} \bracks*{k - (2^n + 1)/2} : k \in \N, 1 \leq k \leq 2^n} \subset \R \subset \R^n$
  \begin{equation*}
    T_{2,n} = \set*{2\sqrt{n} \parens*{k - \frac{2^n + 1}{2}} :
      % k \in \N, 1 \leq k \leq 2^n
      k = 1, 2, \dots, 2^n
    } \subset \R \subset \R^n
    \, ,
  \end{equation*}
  via the embedding $\R \to \R^n$, $x \mapsto (x, 0)$.
  Both sets $T_{1,n}, T_{2,n}$ are symmetric and have cardinality $2^n$.
  Let $\phi_n : T_{1,n} \to T_{2,n}$ be an arbitrary bijection.
  Note that for every $t_1, t_1' \in T_{1,n}$, one has $\norm{t_1 - t_1'} \leq 2 \sqrt{n}$, while for every $t_2, t_2' \in T_{2,n}$ with $t_2 \neq t_2'$ one has $\norm{t_2 - t_2'} \geq 2 \sqrt{n}$.
  Since $\phi_n$ is injective, we deduce that it is also expansive.

  Now, denoting $X = (X_1, \dots, X_n) \sim \gaussdist (0, I_n)$, one has
  \begin{align*}
    \sup_{t \in T_{1,n}} \innerp{t}{X}
    = \max_{1 \leq i \leq n} \abs*{X_i} \, ,
    \quad \text{ while } \quad
    \sup_{t \in T_{2,n}} \innerp{t}{X}
    = \sqrt{n} \parens*{2^n - 1} \abs*{X_1}
    \leq \sqrt{n} 2^{n} \abs*{X_1}
    \, .
  \end{align*}
  Thus, using that the density of $X_1$ is upper bounded by $1/\sqrt{2\pi}$ on $\R$ and lower-bounded by $e^{-1/8}/\sqrt{2\pi}$ on $[-1/2, 1/2]$, we deduce that for $0<\eps \leq 1$,
  \begin{align*}
    \Probab*{\sup_{t \in T_{1,n}} \innerp{t}{X} \leq \eps}
    &= \Probab*{\abs*{X_1} \in \bracks*{-\eps, \eps}}^n
      \leq \parens*{\frac{2 \eps}{\sqrt{2 \pi}}}^n
      \leq \eps^n
      \, ; \\
    \Probab*{\sup_{t \in T_{2,n}} \innerp{t}{X} \leq \eps}
    &\geq \Probab*{\abs*{X_1} \in \bracks*{- \frac{\eps}{\sqrt{n} 2^n}, \frac{\eps}{\sqrt{n} 2^n}}}
      \geq \frac{e^{-1/8}}{\sqrt{2\pi}} \frac{2 \eps}{\sqrt{n} 2^n}
      > \frac{\eps}{2 \sqrt{n} 2^n}
      \geq \frac{\eps}{4^n}
      \, .
  \end{align*}
  In particular,
  \begin{equation}
    \label{eq:no-comparison-example}
    \Probab*{\sup_{t \in T_{1,n}} \innerp{t}{X} \leq \frac{1}{n}}
    \leq \frac{1}{n^n}
    \, , \quad \text{while} \quad
    \Probab*{\sup_{t \in T_{2,n}} \innerp{t}{X} \leq \frac{1}{2^n}}
    > \frac{1}{8^n}
    \, .
  \end{equation}
  We then let $c_n = \log(n)/\log(8) = \log_8 n \to \infty$ as $n \to \infty$, and $u_n = 1/(n c_n) > 0$.
  For some $n_0 \geq 1$, we have $u_n \geq 1/2^n$ for all $n \geq n_0$.
  Thus~\eqref{eq:no-comparison-example} implies that~\eqref{eq:no-comparison-nonconvex} holds for all $n \geq n_0$.
  (To handle the case $1 \leq n < n_0$, we may take $T_{1,n} = T_{2,n} = \set{-1, 1} \subset \R \subset \R^n$, $u_n = 1$ and $c_n = 1/2$, which trivially works as
  $\Probab{\abs*{X_1} \leq 1} > \Probab*{\abs*{X_1} \leq 1}^2 \geq \Probab*{\abs*{X_1} \leq 1/2}^2$.)
%\end{proof}

\section{Low-dimensional sections and projections of convex sets}
\label{sec:low-dimens-sect}

We now turn to the proofs of the results of Section~\ref{sec:low-dim-sections} on the inradius of Gaussian projections (Theorem~\ref{thm:inner-radius-gaussian-proj}) and on the diameter of random sections (Corollary~\ref{cor:diameter-section}).
For completeness, we also include a proof of the formulation we provided of the classical estimate on the diameter of Gaussian projections (Proposition~\ref{prop:diameter-gaussian-proj}).

\subsection{Proof of Theorem~\ref{thm:inner-radius-gaussian-proj}: lower bound}
\label{sec:proof-lower-bound-inradius}

The proof of the lower bound~\eqref{eq:inrad-gaussian-lower} on $\inrad (G_k K)$ relies on well-established tools in proofs of Dvoretzky's theorem following Milman's approach, as described in~\cite[Chapter~4]{pisier1999volume}, \cite[Chapter~5]{artstein2015asymptotic} and~\cite[Chapter~11]{vershynin2018high}, together with a simple inclusion argument.

The following lemma is essentially standard from the Gaussian approach to Dvoretzky's theorem.

\begin{lemma}
  \label{lem:dvo-lower-inrad}
  For any symmetric convex body $K$ in $\R^n$, any $k \geq 1$ and $s \in \R^+$ we have
  \begin{equation}
    \label{eq:high-prob-lower-dvo}
    \Probab*{\inrad (G_k K) \leq w (K) - \rad (K) \parens[\big]{\sqrt{k} + \sqrt{2 s}}}
    \leq \exp \parens*{- s}
    \, .
\end{equation}  
\end{lemma}

\begin{proof}[Proof of Lemma~\ref{lem:dvo-lower-inrad}]
  We first note that for $r \in \R^+$, we have $\inrad (G_k K) \geq r$ if and only if (as $K$ is closed) $r B_2^k \subset G_k K$, which amounts to saying that for every $v \in \R^k$,
  \begin{equation*}
    \sup_{t \in K} \innerp{G_k t}{v}
    \geq \sup_{\theta \in r B_2^k} \innerp{\theta}{v}
    = r \norm{v}
    \, .
  \end{equation*}
  Thus, letting $G_k^* : \R^k \to \R^n$ denote the adjoint map of $G_k$, given by $G_k^* (v) = \sum_{j=1}^k v_j X_j$ for $v = (v_1, \dots, v_k) \in \R^k$, and $h_K (x) = \sup_{t \in K} \innerp{t}{x}$ for $x \in \R^n$, we have
  \begin{equation}
    \label{eq:inrad-inf-sup}
    \inrad (G_k K)
    = \inf_{v \in S^{k-1}} \sup_{t \in K} \innerp{t}{G_k^* v}
    = \inf_{v \in S^{k-1}} h_K \parens[\bigg]{\sum_{j=1}^k v_j X_j}
    \, .
  \end{equation}
  Now, by~\cite[Theorem~9.2.4]{artstein2015asymptotic} (a consequence of Gordon's min-max comparison inequality~\cite{gordon1985some}), we have
  \begin{equation}
    \label{eq:comparison-gordon-dvoretzky}
    \Expect[\bigg]{\inf_{v \in S^{k-1}} h_K \parens[\bigg]{\sum_{j=1}^k v_j X_j}}
    \geq \Expect*{h_K (X_1)} - \sup_{t \in K} \norm{t} \Expect*{\norm{X_1}}
    \geq w (K) - \rad (K) \sqrt{k}
    \, .
  \end{equation}
  Now, endowing $\R^{k \times n}$ with the Hilbert-Schmidt scalar product and its associated norm $\norm{\cdot}_2$ (so that $G_k \sim \gaussdist(0, I_{k \times n})$), the map which associates to $g_k \in \R^{k \times n}$ the quantity $\inrad (g_k K)$ is $\rad (K)$-Lipschitz, as from~\eqref{eq:inrad-inf-sup} we have
  for $g_k, g_k' \in \R^{k \times n}$,
  \begin{align*}
    \abs*{\inrad (g_k K) - \inrad(g_k' K)}
    &= \abs*{\inf_{v \in S^{k-1}} \sup_{t \in K} \innerp{g_k t}{v} - \inf_{v \in S^{k-1}} \sup_{t \in K} \innerp{g'_k t}{v}} \\
    &\leq \sup_{t \in K} \sup_{v \in S^{k-1}} \abs*{\innerp{(g_k - g_k') t}{v}} \\
    &= \sup_{t \in K} \norm{(g_k - g_k') t}
      \leq \norm{t} \norm{g_k - g_k'}_2 \\
    &\leq \rad (K) \norm{g_k - g_k'}_2
      \, .
  \end{align*}
  Thus, Gaussian concentration (Theorem~\ref{thm:gaussian-concentration}) and together with~\eqref{eq:inrad-inf-sup} and~\eqref{eq:comparison-gordon-dvoretzky} concludes the proof.

  We note that a similar inequality, but with $\rad (K)$ replaced by $c \rad (K)$ for some universal constant $c>0$ (which would have sufficed for the proof of Theorem~\ref{thm:inner-radius-gaussian-proj}), could have been obtained by using the sub-Gaussian property~\cite[Theorem~11.1.6]{vershynin2018high}, as opposed to Gordon's min-max comparison inequality.
\end{proof}

We now proceed to the proof of the lower bound on the inradius.
Let $u \in \ooint{0}{w(K)}$ and $\eps \in (0, 1)$.
Denote $r = r_* (K, u)$ (defined in~\eqref{eq:def-radius-widths}), so that $w (K \cap r B_2) = u$ and $\rad (K \cap r B_2) \leq r$.
In addition, since $G_k (K \cap r B_2) \subset G_k K$, we have $\inrad (G_k K) \geq \inrad (G_k (K \cap r B_2))$.
Thus, by Lemma~\ref{lem:dvo-lower-inrad} applied to $K \cap r B_2$, with probability at least $1-e^{-k/2}$ we have
\begin{equation*}
  \inrad (G_k K)
  \geq \inrad (G_k (K \cap r B_2))
  > w (K \cap r B_2) - r \parens[\big]{\sqrt{k} + \sqrt{k}}
  = u - 2 r \sqrt{k}
  \, .
\end{equation*}
Thus, if
\begin{equation*}
  k
  \leq \frac{\eps^2 u^2}{4 r^2}
  = \frac{\eps^2 u^2}{4 r_* (K, u)^2}
  \, ,
\end{equation*}
then the previous inequality implies that $\inrad (G_k K) \geq (1-\eps) u$.
Taking $\eps = 3/4$ and bounding $(3/4)^2/4 \geq 1/8$, we obtain that 
\begin{equation}
  \label{eq:proof-lower-inrad}
  \Probab*{\inrad (G_k u) > \frac{u}{4}}
  \geq 1- e^{-k/2}
  \quad \text{ whenever } \quad
  k
  \leq \frac{u^2}{8 r_* (K, u)^2}
  \, .
\end{equation}

\subsection{Proof of Theorem~\ref{thm:inner-radius-gaussian-proj}: upper bound}
\label{sec:proof-upper-inrad}

We now turn to the main part of the proof, namely that of the upper bound on the inradius.

\begin{lemma}
  \label{lem:norm-preserv}
  Let $n, k \geq 1$, $T \subset \R^n$ a closed set, and $G_k : \R^n \to \R^k$ beX the Gaussian linear map given by~\eqref{eq:def-gaussian-matrix}.
  For every $s \in \R^+$, we have
  \begin{equation}
    \label{eq:norm-preserv}
    \Probab*{\inf_{t \in T} \norm{G_k t} \leq \sqrt{\frac{k}{2}} \inf_{t \in T} \norm{t} - w (T) - \rad (T) \sqrt{2 s}}
    \leq \exp (-s)
    \, .
  \end{equation}
\end{lemma}

\begin{proof}[Proof of Lemma~\ref{lem:norm-preserv}]
  By~\cite[Theorem~9.3.3]{artstein2015asymptotic} letting $X \sim \sgn, Y \sim \gaussdist (0, I_k)$ we have
  \begin{equation*}
    \Expect*{\inf_{t \in T} \norm{G_k t}}
    \geq \Expect*{\norm{Y}} \inf_{t \in T} \norm{t} - \Expect*{\sup_{t \in T} \innerp{t}{X}}
    \geq \sqrt{\frac{k}{2}} \inf_{t \in T} \norm{t} - w (T)
    \, ,
  \end{equation*}
  where we used from~\cite{haagerup1981best} that $\Expect{\norm{Y}} \geq \sqrt{k/2}$.
  In addition, the map $g \in \R^{k \times n} \mapsto \inf_{t \in T} \norm{g t}$ is $\rad(T)$-Lipschitz with respect to the Hilbert-Schmidt norm $\norm{\cdot}_2$ on $\R^{k \times n}$.
  The result then follows from Gaussian concentration (Theorem~\ref{thm:gaussian-concentration}).
\end{proof}

The following fact is key to the proof of Theorem~\ref{thm:inner-radius-gaussian-proj}: the Gaussian width of a Gaussian projection of a convex body is controlled (with high probability) in terms of the Gaussian width of the original body.

\begin{lemma}
  \label{lem:gaussian-width-proj}
  Let $n, k\geq 1$, $K \subset \R^n$ a symmetric convex body and $G_k : \R^n \to \R^k$ the Gaussian linear map given by~\eqref{eq:def-gaussian-matrix}.
  For every $s \in \R^+$, we have
  \begin{equation}
    \label{eq:gaussian-width-proj}
    \Probab*{w \parens*{G_k K} \geq \parens[\big]{\sqrt{k} + \sqrt{2 s}} w (K)}
    \leq \exp \parens*{- s}
    \, .
  \end{equation}
\end{lemma}

\begin{proof}[Proof of Lemma~\ref{lem:gaussian-width-proj}]
  We first bound the expectation of $w (G_k K)$ and then control its fluctuations using Gaussian concentration.

  Let $X \sim \gaussdist (0, I_n)$ be independent of $G_k$.
  We may write $X = \norm{X} V$, where $V$ is uniformly distributed over $S^{n-1}$ and independent of $\norm{X}$.
  Next, note that for every $v = (v_1, \dots, v_k) \in S^{n-1}$, we have $G_k^* v = v_1 X_1 + \dots + v_k X_k \sim \gaussdist (0, I_n)$.
  Since $V$ is independent of $G_k^*$, we deduce that $Y = G_k^* V$ has distribution $\gaussdist (0, I_n)$.
  Finally, $Y$ is independent of $\norm{X}$ as $G_k, V$ are.
  With these prerequisites in place, we have
  \begin{equation*}
    \Expect*{w (G_k K)}
    = \Expect*{\sup_{t \in K} \innerp{G_k t}{X}}
    = \Expect*{\sup_{t \in K} \innerp{t}{G_k^* X}}
    = \Expect*{\norm{X} \sup_{t \in K} \innerp{t}{Y}} 
    = \Expect*{\norm{X}} \Expect*{\sup_{t \in K} \innerp{t}{Y}}
    .
  \end{equation*}
  Recalling that $Y \sim \gaussdist (0, I_n)$ and using that $\Expect{\norm{X}} \leq \Expect{\norm{X}^2}^{1/2} = \sqrt{k}$, we deduce that
  \begin{equation}
    \label{eq:expected-width-proj}
    \Expect*{w (G_k K)}
    \leq \sqrt{k} \, w (K)
    \, .
  \end{equation}

  We now show that the map $\R^{k \times n} \to \R^+$, $g \mapsto w (g K)$ is $w(K)$-Lipschitz.
  Let $g, g' \in \R^{k \times n}$ and $X \sim \gaussdist (0, I_n)$.
  We have
  \begin{align*}
    \abs*{w (g K) - w (g' K)}
    &\leq \Expect*{\abs*{\sup_{t \in K} \innerp{g t}{X} - \sup_{t \in K} \innerp{g' t}{X}}}
      \leq \Expect*{\sup_{t \in K} \abs*{\innerp*{(g - g') t}{X}}}  \\
    &= \Expect*{\sup_{t \in K} \innerp*{(g - g') t}{X}}
  \end{align*}
  as $K$ is symmetric.
  Now, for any $t, t' \in K$ we have $\norm{(g - g') t - (g - g') t'} \leq \norm{g-g'}_2 \norm{t - t'}$.
  Thus, by the Slepian-Sudakov-Fernique comparison inequality (\eg \cite[Theorem~7.2.11]{vershynin2018high}),
  \begin{equation*}
    \Expect*{\sup_{t \in K} \innerp*{(g - g') t}{X}}
    \leq \norm{g - g'}_2 \Expect*{\sup_{t \in K} \innerp*{t}{X}}
    \, .
  \end{equation*}
  Combining the previous inequalities then yields $\abs{w (g K) - w (g' K)} \leq \norm{g - g'}_2 w (K)$; namely, the map $g \mapsto w (g K)$ is $w (K)$-Lipschitz.
  (We note that using Talagrand's majorizing measures theorem, Theorem~\ref{thm:talagrand-maj-meas}, instead of the comparison inequality would have shown that this map is $c w (K)$-Lipschitz for some universal constant $c > 0$, which would suffice for our purposes.)

  Finally, combining the expectation bound~\eqref{eq:expected-width-proj} with Gaussian concentration (Theorem~\ref{thm:gaussian-concentration}) concludes the proof of Lemma~\ref{lem:gaussian-width-proj}.
\end{proof}

We are now in position to prove the upper bound on the inner radius.

Let $u \in \ooint{0}{w(K)}$ and $r = r_* (K, u)$.
Let $k$ be an integer satisfying
\begin{equation}
  \label{eq:cond-k-large-preserv}
  k
  \geq \frac{144 u^2}{r^2}
  \, .
\end{equation}

Applying Lemma~\ref{lem:norm-preserv} to $T = K \cap r S^{n-1}$ and using that $w (T) \leq w (K \cap rB_2) = u$, we obtain that with probability at least $1 - e^{-k/10}$,
\begin{equation*}
  \inf_{t \in K \cap r S^{n-1}} \norm{G_k t}
  > \sqrt{\frac{k}{2}} \cdot r - u - r \sqrt{\frac{k}{5}}
  \geq \frac{r \sqrt{k}}{3} - u
  \, .
\end{equation*}
Since
$r \sqrt{k} \geq 12 u$
from the condition~\eqref{eq:cond-k-large-preserv}, we deduce that $\Probab*{E_1} \geq 1 - e^{-k/10}$, where $E_1$ denotes the event
\begin{equation}
  \label{eq:proof-upper-inrad-event-1}
  E_1
  = \set*{\inf_{t \in K \cap r S^{n-1}} \norm{G_k t} > 3 u % \frac{r \sqrt{k}}{6}
  }
  \, .
\end{equation}
From now on, we work under the event $E_1$.
In this case, for any $t \in K$ with $\norm{t} > r$, since $K$ is symmetric convex the point $t' = \frac{r}{\norm{t}} t$ belongs to $K \cap r S^{n-1}$, thus
\begin{equation*}
  \norm{G_k t}
  = \frac{\norm{t}}{r} \norm{G_k t'}
  > 3 u
  \, .
\end{equation*}
In particular, $G_k (K \setminus r B_2) \cap 3 u B_2^k = \varnothing$.
Thus, if $\inrad (G_k K) \geq 3 u$, namely $3 u B_2^k \subset G_k K$, then $3 u B_2^k \subset G_k (K \cap r B_2)$;
in particular, one would have
\begin{equation}
  \label{eq:proof-upper-inrad-ineq-width}
  3 u \sqrt{\frac{k}{2}}
  \leq w (u B_2^k)
  \leq w \parens*{G_k (K \cap r B_2)}
  \, .
\end{equation}

On the other hand, Lemma~\ref{lem:gaussian-width-proj} applied to $K \cap r B_2$ (with $w (K \cap r B_2) = u$) shows that
\begin{equation}
  \label{eq:proof-upper-inrad-event-2}
  \Probab*{E_2}
  \geq 1 - e^{-k/2} \, ,
  \quad \text{ where } \quad
  E_2
  = \set*{w \parens*{G_k (K \cap r B_2)} \leq 2 u \sqrt{k}}
  \, .
\end{equation}
Thus, under $E_2$, one has
\begin{equation*}
  w \parens*{G_k (K \cap r B_2)}
  \leq 2 u \sqrt{k}
  < 3 u \sqrt{\frac{k}{2}}
  \, ,
\end{equation*}
and thus from the criterion~\eqref{eq:proof-upper-inrad-ineq-width}, under $E_1 \cap E_2$ (such that $\Probab{E_1 \cap E_2} \geq 1 - 2 e^{-k/10}$) one has $\inrad (G_k u) < 3 u$.

In conclusion, we have shown that
\begin{equation}
  \label{eq:proof-upper-inrad}
  \Probab[\big]{\inrad (G_k u) < 3 u}
  \geq 1- 2 e^{-k/10}
  \quad \text{ whenever } \quad
  k
  \geq \frac{144 u^2}{r_* (K, u)^2}
  \, .
\end{equation}
This concludes the proof of Theorem~\ref{thm:inner-radius-gaussian-proj}, which holds with $c = 144$.

\subsection{Proof of Corollary~\ref{cor:diameter-section}}
\label{sec:proof-radius-section}

The argument is standard, and we include it merely for completeness.

Let $K = L^\circ$, which is a symmetric convex body in $\R^n$ with $K \neq \set*{0}$.

Let $F = (\ker G_k)^\perp = \im (G_k^*)$ denote the orthogonal of the kernel of $G_k$, which is the image of $G_k^* : \R^k \to \R^n$.
Then $F$ follows the uniform distribution over $k$-dimensional subspaces of $\R^n$.

By \eg~\cite[Theorem~2.13 p.~353]{davidson2001local}, denoting by $\smin (G_k^*), \smax (G_k^*)$ the smallest and largest singular values of $G_k^*$, one has for every $1 \leq k \leq n$ and $s \in \R^+$:
\begin{align}
  \label{eq:smax-gaussian-matrix}
  \Probab*{\smax (G_k^*) \geq \sqrt{n} + \sqrt{k} + \sqrt{2 s}}
  &\leq \exp (-s)
    \, ; \\
  \label{eq:smin-gaussian-matrix}
  \Probab*{\smin (G_k^*) \leq \sqrt{n} - \sqrt{k} - \sqrt{2 s}}
  &\leq \exp (-s)
    \, .
\end{align}
From~\eqref{eq:smax-gaussian-matrix} and~\eqref{eq:smin-gaussian-matrix}, with probability at least $1 - 2 \exp (-k/100)$ the following two conditions hold:
\begin{align*}
  \smin (G_k^*)
  &> \sqrt{n} - \sqrt{k} - \sqrt{2 \cdot \frac{k}{100}}
  \geq \sqrt{n} - \sqrt{\frac{n}{2}} - \sqrt{\frac{n}{100}}
  \geq \frac{\sqrt{n}}{6}
  \, ; \\
  \smax (G_k^*)
  &< \sqrt{n} + \sqrt{k} + \sqrt{2 \cdot \frac{k}{100}}
  \leq 2 \sqrt{n}
  \, .
\end{align*}
We now work under the event $E_1 = \set{\sqrt{n}/6 \leq \smin (G_k^*) \leq \smax (G_k^*) \leq 2 \sqrt{n}}$, with $\Probab{E_1} \geq 1 - 2\exp (-k/100)$.
Under $E_1$, the restriction $G_k|_F : F \to \R^k$ is a linear bijection, with
\begin{equation*}
  \frac{\sqrt{n}}{6} \norm{t}
  \leq \norm{G_k t}
  \leq 2 \sqrt{n} \norm{t}
\end{equation*}
for every $t \in F$.
Therefore, for any $v \in \R^k$ we have
\begin{equation}
  \label{eq:proof-gauss-orth-proj}
  \frac{\norm{v}}{2 \sqrt{n}}
  \leq \norm{(G_k|_F)^{-1} v}
  \leq \frac{6 \norm{v}}{\sqrt{n}}
  \, .
\end{equation}

Let $P_F : \R^n \to F$ denote the orthogonal projection on $F$. 
Since $F^\top = \ker G_k$, we have $G_k K = G_k (P_F K)$.
Thus $P_F K = (G_k|_F)^{-1} (G_k K)$, and~\eqref{eq:proof-gauss-orth-proj} then implies that, denoting $\inrad_F (A) = \sup \set*{r \geq 0 : (r B_2^n) \cap F \subset A}$ for $A \subset F$,
\begin{equation}
  \label{eq:proof-inrad-gauss-proj}
  \frac{\inrad (G_k K)}{2 \sqrt{n}}
  \leq \inrad (P_k K)
  \leq \frac{6 \inrad (G_k K)}{\sqrt{n}}
  \, .
\end{equation}

Now, let $x \in F$.
We have $x \in (P_F K)^\circ$ if and only if $\sup_{t \in K} \innerp{P_F t}{x} \leq 1$.
But since $\innerp{P_F t}{x} = \innerp{t}{x}$ for any $t \in \R^n$ as $x \in F$, the latter condition amounts to $\sup_{t \in K} \innerp{t}{x} \leq 1$, namely $x \in K^\circ$.
We therefore have $(P_F K)^\circ \cap F = K^\circ \cap F$.
Thus $\rad (K^\circ \cap F) = \inrad_F (P_F K)^{-1}$.

It follows from the above (and the fact that $L = K^\circ$ and $\diam (L \cap F) = 2 \rad (L \cap F)$) that
\begin{equation}
  \label{eq:proof-inrad-diam-sec}
  \Probab*{
    \frac{\sqrt{n}}{3 \inrad (G_k L^\circ)}
    \leq \diam (L \cap F)
    \leq \frac{4 \sqrt{n}}{\inrad (G_k L^\circ)}
  }
  \geq 1 - 2 \exp \parens*{- \frac{k}{100}}
  \, .
\end{equation}

Corollary~\ref{cor:diameter-section} then follows from~\eqref{eq:proof-inrad-diam-sec} and Theorem~\ref{thm:inner-radius-gaussian-proj}.

\subsection{Proof of Proposition~\ref{prop:diameter-gaussian-proj}}
\label{sec:proof-diam-proj}

We now provide a proof of the well-known estimate on the diameter of low-dimensional projections, in the form given in Proposition~\ref{prop:diameter-gaussian-proj}.
We start with the lower bound on the diameter.

On the one hand, using that $\norm{G_j t} \geq \abs{\innerp{t}{X_1}}$ for any $t \in \R^n$, we have
\begin{equation*}
  \diam (G_j T)
  = \sup_{t, t' \in T} \norm*{G_j t - G_j t'}
  \geq \sup_{t, t' \in T} \abs*{\innerp*{t - t'}{X_1}}
  \, .
\end{equation*}
Recalling that $X_1 \sim \gaussdist (0, I_n)$ and from~\eqref{eq:def-gaussian-width}, we have
\begin{equation*}
  \Expect*{\sup_{t, t' \in T} \abs*{\innerp*{t - t'}{X_1}}}
  = 2  w (T)
  \, .
\end{equation*}
Since in addition the map $x \mapsto \sup_{t, t' \in T} \abs{\innerp{t - t'}{x}}$ is $\diam (T)$-Lipschitz on $\R^n$, the Gaussian concentration inequality (Theorem~\ref{thm:gaussian-concentration}) together with the previous two inequalities shows that
% for every $s \in \R^+$,
\begin{equation}
  \label{eq:proof-diam-proj-lower-width}
  \Probab*{\diam (G_j T) \geq 2 w (T) - \diam (T) \sqrt{j/8}}
  \geq 1 - e^{-j/16}
  \, .
\end{equation}

On the other hand, let $t, t' \in T$ such that $\norm{t-t'} = \diam (T)$ (which exist by compactness of $T$).
Then $\diam (G_j T) \geq \norm{G_j (t - t')}$.
Since $G_j (t - t') \sim \gaussdist (0, \norm{t-t'}^2 I_j)$, a standard tail bound on the $\chi^2$ distribution (\eg,~\cite[Lemma~1]{laurent2000adaptive}) shows that
\begin{equation}
  \label{eq:proof-diam-proj-lower-dim}
  \Probab*{\diam (G_j T) \geq \diam (T) \sqrt{j/2}}
  \geq \Probab*{\norm{G_j (t-t')}^2 \geq \frac{\norm{t-t'}^2 j}{2}}
  \geq 1 - e^{-j/16}
  \, .
\end{equation}

Thus, taking the best of the two lower bounds~\eqref{eq:proof-diam-proj-lower-width} and~\eqref{eq:proof-diam-proj-lower-dim} (depending on the values of $w (T), \diam (T), j$), we deduce that with probability at least $1-e^{-j/16}$,
\begin{align}
  \diam (G_j T)
  &\geq \max \braces*{2 w (T) - \diam (T) \sqrt{j/8}, \diam (T) \sqrt{j/2}} \nonumber \\
  &\geq \frac{1}{2} \braces*{2 w (T) - \diam (T) \sqrt{j/8} + \diam (T) \sqrt{j/2}} \nonumber \\
  &= w (T) + \frac{1}{4 \sqrt{2}} \diam (T) \sqrt{j}
    \label{eq:proof-diam-proj-lower}
    \, .
\end{align}

We next provide the upper bound on the diameter, following~\cite{artstein2015asymptotic,vershynin2018high}.
By a standard $\eps$-net argument (see~\cite[pp.~171--172]{vershynin2018high}), there exists a subset $N_j$ of the unit sphere $S^{j-1}$ with $|N_j| \leq 5^j$ such that $\norm{\theta} \leq 2 \sup_{v \in N_j} \innerp{\theta}{v}$ for any $\theta \in \R^j$.
We thus have
\begin{equation}
  \label{eq:proof-diam-proj-upper-net}
  \diam (G_j T)
  = \sup_{t, t' \in T} \norm*{G_j (t-t')}
  \leq 2 \sup_{t, t' \in T, v \in N_j} \innerp*{G_j (t-t')}{v}
  = 2 \sup_{v \in N_j} \sup_{t, t' \in T} \innerp*{t-t'}{G_j^* v}
  .
\end{equation}
Now, for any $v \in N_j \subset S^{j-1}$, one has $G_j^* v \sim \gaussdist (0, I_n)$.
Thus as above,
\begin{equation*}
  \Probab[\bigg]{\sup_{t, t' \in T} \innerp*{t-t'}{G_j^* v} \geq 2 w (T) + \diam (T) \sqrt{2 j \log (10)}}
  \leq 10^{-j}
  \, .
\end{equation*}
Plugging this bound into~\eqref{eq:proof-diam-proj-upper-net} and using a union bound yields:
\begin{align}
  &\Probab*{\diam (G_j T) \geq 4 w (T) + 2 \diam (T) \sqrt{2 j \log (10)}} \nonumber \\
  &\leq \sum_{v \in N_j} \Probab[\bigg]{\sup_{t, t' \in T} \innerp*{t-t'}{G_j^* v} \geq 2 w (T) + \diam (T) \sqrt{2 j \log (10)}} \nonumber \\
  &\leq 5^j 10^{-j}
    = 2^{-j}
    \leq e^{-j/16}
    \label{eq:proof-diam-proj-upper}
    \, .
\end{align}

Combining the lower and upper bounds~\eqref{eq:proof-diam-proj-lower} and~\eqref{eq:proof-diam-proj-upper} with a union bound concludes.

\section{Basic properties of suprema of Gaussian processes}
\label{sec:basic-prop-supr}

In this section, we collect proofs of several essentially standard estimates for Gaussian processes, in the precise forms used elsewhere in the paper.

\subsection{A simple lower tail estimate}
\label{sec:simple-lower-tail}

We start with a simple lower tail estimate, which was used in particular in Section~\ref{sec:chain-char-quant}.

\begin{lemma}
  \label{lem:mini-sb}
  For any symmetric compact set $T \subset \R^n$ with $T \neq \set{0}$ and $\eps \in \ocint{0}{1/2}$, 
  \begin{equation}
    \label{eq:mini-sb}
    \Probab*{\sup_{t \in T} \innerp{t}{X} \leq \frac{\sqrt{\pi}}{4} \frac{\eps}{\sqrt{\log (1/\eps)}} w (T)}
    \leq \eps
    \, .
  \end{equation}
\end{lemma}

We note that the $\sqrt{\log (1/\eps)}$ factor in~\eqref{eq:mini-sb} can be replaced by a constant, by an elegant theorem of Lata{\l}a and Oleszkiewicz~\cite{latala1999gaussian} (confirming a conjecture of L.~A.~Shepp).
However, we do not need this refinement here, as we will only use Lemma~\ref{lem:mini-sb} when $\eps$ is of constant order.

\begin{proof}[Proof of Lemma~\ref{lem:mini-sb}]
  Let $\eta \leq 1/2$.
  On the one hand, by Gaussian concentration (Theorem~\ref{thm:gaussian-concentration}), denoting $d = (w(T)/\rad (T))^2$ we have
  \begin{align}
    \label{eq:proof-mini-1}
    \Probab*{\sup_{t \in T} \innerp{t}{X} \leq \eta w (T)}
    &\leq \Probab*{\sup_{t \in T} \innerp{t}{X} - w (T) \leq - (1-\eta) w (T)} \nonumber \\
    &\leq \exp \parens*{- \frac{(1-\eta)^2 w (T)^2}{2 \rad (T)^2}}
      \leq \exp \parens*{- \frac{d}{8}}
      \, .
  \end{align}
  On the other hand, letting $t_0 \in T$ such that $\norm{t_0} = \rad (T)$ and letting $Z \sim \gaussdist (0, 1)$, using that the density of $Z$ is upper bounded by $1/\sqrt{2\pi}$ we obtain
  \begin{align}
    \label{eq:proof-mini-2}
    \Probab*{\sup_{t \in T} \innerp{t}{X} \leq \eta w (T)}
    &\leq \Probab*{\abs*{\innerp{t_0}{X}} \leq \eta w(T)}
    = \Probab*{\rad (T) \abs*{Z} \leq \eta w (T)} \nonumber \\
    &\leq \frac{1}{\sqrt{2\pi}} \frac{2 \eta w (T)}{\rad (T)}
      = \sqrt{{2}/{\pi}} \eta \sqrt{d}
      \, .
  \end{align}
  Finally, let $\eta = \frac{\sqrt{\pi}}{4} \frac{\eps}{\sqrt{\log (1/\eps)}}$.
  Since $\eps \leq 1/2$, we have $\eta \leq \frac{\sqrt{\pi}}{4 \sqrt{\log 2}} \frac{1}{2} \leq 1/2$.

  Now, if $d \geq 8 \log (1/\eps)$, inequality~\eqref{eq:proof-mini-1} implies that
  \begin{equation*}
    \Probab*{\sup_{t \in T} \innerp{t}{X} \leq \eta w (T)}
    \leq \eps
    \, .
  \end{equation*}
  On the other hand, if $d < 8 \log (1/\eps)$, inequality~\eqref{eq:proof-mini-2} yields
  \begin{equation*}
    \Probab*{\sup_{t \in T} \innerp{t}{X} \leq \eta w (T)}
    \leq \sqrt{\frac{2}{\pi}} \frac{\sqrt{\pi}}{4} \frac{\eps}{\sqrt{\log (1/\eps)}} \sqrt{8 \log (1/\eps)}
    = \eps
    \, .
  \end{equation*}
  This concludes the proof.
\end{proof}

As a direct consequence and together with Talagrand's Majorizing Measures theorem (Theorem~\ref{thm:talagrand-maj-meas}), we deduce the following:

\begin{lemma}
  \label{lem:quantiles-constant}
  For any $\kappa_1 > 1$, there exists a constant $\kappa_2 \geq 1$ depending only on $\kappa_1$ such that the following holds.
  For any $n \geq 1$, symmetric and compact subset $T \subset \R^n$, and $s \in \ccint{1/\kappa_1}{\kappa_1}$,
  \begin{equation}
    \label{eq:quantiles-constant}
    \frac{1}{\kappa_2} \gamma_2 (T)    
    \leq Q \parens*{s, \sup_{t \in T} \innerp{t}{X}}
    \leq \kappa_2 \gamma_2 (T)
    \, .
  \end{equation}
\end{lemma}

\begin{proof}[Proof of Lemma~\ref{lem:quantiles-constant}]
  Let $c_0 \geq 1$ denote the constant in Theorem~\ref{thm:talagrand-maj-meas}.
  The upper bound in~\eqref{eq:quantiles-constant} follows from Markov's inequality (as $\sup_{t \in T} \innerp{t}{X} \geq 0$ since $T$ is symmetric), which yields:
  \begin{equation*}
    Q \parens*{s, \sup_{t \in T} \innerp{t}{X}}
    \leq Q \parens*{\kappa_1^{-1}, \sup_{t \in T} \innerp{t}{X}}
    \leq \frac{w (T)}{1 - \exp(-\kappa_1^{-1})}
    \leq \frac{c_0}{1- \exp (-\kappa_1^{-1})} \gamma_2 (T)
    \, .
  \end{equation*}
  The lower bound in~\eqref{eq:quantiles-constant} follows from Lemma~\ref{lem:mini-sb} with $\eps = \exp (-\kappa_1)$, which yields
  \begin{equation*}
    Q \parens*{s, \sup_{t \in T} \innerp{t}{X}}
    \geq Q \parens*{\kappa_1, \sup_{t \in T} \innerp{t}{X}}
    \geq \frac{\sqrt{\pi}}{4} \frac{\exp (-\kappa_1)}{\sqrt{\kappa_1}} w (T)
    \geq \frac{\sqrt{\pi}}{4 c_0} \frac{\exp (-\kappa_1)}{\sqrt{\kappa_1}} \gamma_2 (T)
    \, .
  \end{equation*}
  This concludes the proof.
\end{proof}

\subsection{Proof of Corollary~\ref{cor:upper-tail-gaussian}}
\label{sec:proof-cor-upper-tail}

We first prove the upper bound in~\eqref{eq:upper-tail-gaussian}.
Since the function $f : \R^n \to \R$ defined by $f (x) = \sup_{t \in K} \innerp{t}{x}$ is $\rad (K)$-Lipschitz and $\Expect*{f (X)} = w (K)$, by Theorem~\ref{thm:gaussian-concentration}:
\begin{align*}
  \Probab*{\sup_{t \in K} \innerp{t}{X} \geq u}
  &= \Probab*{f (X) - \Expect*{f (X)} \geq u - w (K)} \\
  &\leq \Probab*{f (X) - \Expect*{f (X)} \geq \frac{u}{2}}
    \leq \exp \parens*{- \frac{u^2}{8 \rad (K)^2}}
    \, .
\end{align*}
For the lower bound, by compactness, there exists $t_0 \in K$ such that $\norm{t_0} = \rad (K)$.
Since $K$ is symmetric, $\set{t_0, - t_0} \subset K$ and thus, denoting $Z = \innerp{t_0}{X}/\rad (K) \sim \gaussdist (0, 1)$,
\begin{equation*}
  w (K)
  \geq \Expect*{\abs*{\innerp{t_0}{X}}}
  = \rad (K) \Expect*{\abs*{Z}}
  = \sqrt{\frac{2}{\pi}} \rad (K)
  \, ,
\end{equation*}
and therefore $u \geq 2 w (K) \geq 2 \sqrt{2/\pi} \rad (K) \geq \rad(K)$.
Using Fact~\ref{fac:single-gauss-upper-tail} below, we deduce that
\begin{equation*}
  \Probab*{\sup_{t \in K} \innerp{t}{X} \geq u}
  \geq \Probab*{\innerp{t_0}{X} \geq u}
  = \Probab*{Z \geq \frac{u}{\rad (K)}}
  \geq \exp \parens*{- \frac{3 u^2}{\rad (K)^2}}
  \, .
\end{equation*}
This proves~\eqref{eq:gaussian-large-balls}.
We deduce~\eqref{eq:upper-tail-gaussian} by letting $K = \conv (T-T)$, since $\rad (K) = \diam (T)$, $w (K) = 2 w (T)$ and $\sup_{t \in K} \innerp{t}{X} = \sup_{t,t' \in T} \abs{\innerp{t-t'}{X}}$.

\begin{fact}
  \label{fac:single-gauss-upper-tail}
  Let $Z \sim \gaussdist (0, 1)$.
  For every $u \geq 1$, one has
  \begin{equation}
    \label{eq:single-gauss-upper-tail}
    \Probab*{Z \geq u}
    \geq e^{- 3 u^2}
    \, .
  \end{equation}
\end{fact}

\begin{proof}[Proof of Fact~\ref{fac:single-gauss-upper-tail}]
  Using that $\sqrt{2\pi} \leq e \leq e^{u^2}$ as $u \geq 1$, we obtain
  \begin{equation*}
    \Probab*{Z \geq u}
    = \int_{u}^{\infty} \frac{1}{\sqrt{2\pi}} e^{-v^2/2} \di v
    \geq \frac{1}{\sqrt{2\pi}} \int_{u}^{2u} e^{-v^2/2} \di v
    \geq e^{-u^2} u e^{-(2u)^2/2}
    \geq e^{-3 u^2}
    \, .
    \qedhere
  \end{equation*}
\end{proof}

\subsection{Proof of Corollary~\ref{cor:upper-quantiles-gaussian}}
\label{sec:proof-cor-upper-quantiles}

The inequalities trivially hold if $T$ is reduced to a point, we thus assume that $T$ is nontrivial.

Let $c_0$ be the constant $c$ of Theorem~\ref{thm:talagrand-maj-meas}.
First, by Theorem~\ref{thm:gaussian-concentration}, one has
\begin{equation*}
  \Probab*{\sup_{t,t' \in T} \abs*{\innerp{t-t'}{X}} \geq 2 w (T) + \diam (T) \sqrt{2 s}}
  \leq e^{-s}
  \, ,
\end{equation*}
and therefore, using Theorem~\ref{thm:talagrand-maj-meas},
\begin{equation*}
  Q_+ \parens[\bigg]{s, \sup_{t,t' \in T} \abs*{\innerp{t-t'}{X}}}
  \leq 2 w (T) + \diam (T) \sqrt{2 s}
  \leq 2 c_0 \gamma_2 (T) + \diam (T) \sqrt{2 s}
  \, .
\end{equation*}
This establishes the upper bound in~\eqref{eq:upper-quantiles-gaussian}.

We now turn to the lower bound.
On the one hand, from Lemma~\ref{lem:mini-sb} with $\eps= 1/2$ (as $1-e^{-3} \geq 1/2$) and Theorem~\ref{thm:talagrand-maj-meas},
\begin{equation*}
  Q_+ \parens[\bigg]{s, \sup_{t,t' \in T} \abs*{\innerp{t-t'}{X}}}
  \geq Q_+ \parens[\bigg]{3, \sup_{t \in T-T} \innerp{t}{X}}
  \geq \frac{\sqrt{\pi}}{8 \sqrt{\log 2}} w (T-T)
  \geq \frac{\sqrt{\pi}}{4 c_0 \sqrt{\log 2}} \gamma_2 (T)
  \, .
\end{equation*}
On the other hand, letting $t_1, t_2 \in T$ such that $\norm{t_1-t_2} = \diam (T)$,
from Fact~\ref{fac:single-gauss-upper-tail}, letting $Z \sim \gaussdist (0, 1)$ we have
\begin{equation*}
  Q_+ \parens[\bigg]{s, \sup_{t,t' \in T} \abs*{\innerp{t-t'}{X}}}
  \geq Q_+ \parens*{s, \abs*{\innerp{t_1-t_2}{X}}}
  = Q_+ \parens*{s, \diam (T) Z}
  \geq \diam (T) \sqrt{\frac{s}{3}}
  \, .
\end{equation*}
Combining the previous two inequalities, we obtain
\begin{equation*}
  Q_+ \parens[\bigg]{s, \sup_{t,t' \in T} \abs*{\innerp{t-t'}{X}}}
  \geq \frac{1}{2} \parens*{\frac{\sqrt{\pi}}{4 c_0 \sqrt{\log 2}} \gamma_2 (T) + \rad (T) \sqrt{\frac{s}{3}}}
  \, ,
\end{equation*}
which concludes the proof.

\section{Examples}
\label{sec:examples}

We conclude with detailed examples that substantiate claims made in the Section~\ref{sec:introduction}.
Specifically, Sections~\ref{sec:lack-solution-kuelbs} and~\ref{sec:many-solut-kuelbs} illustrate the limitations of the Kuelbs-Li conditions~\eqref{eq:implicit-sb-kuelbs-li} mentioned in Section~\ref{sec:gaussian-small-ball}: the first provides an example where these conditions admit no solution, while the second exhibits a situation where they admit a range of solutions with widely different orders of magnitude.
Finally, Section~\ref{sec:exampl-oper-norm} establishes the claims made in Example~\ref{ex:operator-norm} concerning the behavior of the operator norm.

\subsection{Lack of solution to the Kuelbs-Li condition}
\label{sec:lack-solution-kuelbs}

In this section, we provide an example of a symmetric convex body $K \subset \R^n$ and a range of values of $u \in \ooint{0}{w (K)}$ for which no value of $s > 0$ satisfies the Kuelbs-Li condition~\eqref{eq:implicit-sb-kuelbs-li}.

Specifically, we will show that for any $\kappa_0 > 1$, there exists $n \geq 1$, a nontrivial symmetric convex body $K \subset \R^n$ and $u \in \ooint{0}{w(K)/\kappa_0}$ such that no $s > 0$ satisfies
\begin{equation}
  \label{eq:cond-kuelbs-li-imposs}
  \frac{s}{\kappa_0}
  \leq \log N \parens*{K, \frac{u}{\sqrt{s}}}
  \leq \kappa_0 s
  \, .
\end{equation}

A simple obstruction to the existence of solutions to~\eqref{eq:cond-kuelbs-li-imposs} is the gap in the Sudakov minoration inequality (recalled in Lemma~\ref{lem:sudakov}).
Indeed, if $r = u/\sqrt{s}$ then the conditions~\eqref{eq:cond-kuelbs-li-imposs} are equivalent to
\begin{equation}
  \label{eq:equiv-kuelbs-li-radius}
  \frac{u}{\sqrt{\kappa_0}}
  \leq r \sqrt{\log N (K, r)}
  \leq \sqrt{\kappa_0} \, u
  \, .
\end{equation}
However, it is well-known~\cite[Chapter~2]{talagrand2021upper} for certain high-dimensional convex bodies $K$, we have $\sup_{r > 0} r \sqrt{\log N (K, r)} \ll w (K)$.
In this case, any $u>0$ with $\sup_{r > 0} r \sqrt{\log N (K, r)} \ll u \ll w (K)$ provides a counterexample.

For a concrete example, consider the ellipsoid $K \subset \R^n$, with half-lengths $(1/\sqrt{i})_{1 \leq i \leq n}$:
\begin{equation*}
  K
  = \set*{t = (t_i)_{1 \leq i \leq n} \in \R^n : \sum_{i=1}^n i t_i^2 \leq 1}
  \, .
\end{equation*}
By a standard moment equivalence for norms of Gaussian vectors (or \eg~\cite[Proposition 2.13.1]{talagrand2021upper}), we have
\begin{equation*}
  w (K)
  \geq \frac{1}{c_0} \parens*{\sum_{i=1}^n \frac{1}{i}}^{1/2}
  \geq \frac{\sqrt{\log n}}{c_1}
\end{equation*}
for some universal constants $c_0, c_1 > 1$.
On the other hand, by standard entropy estimates for ellipsoids (\eg, \cite[Corollary~2.13.8 p.~78]{talagrand2021upper} or~\cite[Lemma~7.3]{mourtada2025universal}), we have
\begin{equation*}
  \sup_{r > 0} \set*{r \sqrt{\log N (K, r)}}
  \leq c_2
\end{equation*}
for some absolute constant $c_2 > 1$.
Thus, if $u = (\log n)^{1/4}$, for $n$ large enough depending on $\kappa_0$, we have $u < w (K)/\kappa_0$ and $r \sqrt{\log N (K, r)} < u/\sqrt{\kappa_0}$ for any $r > 0$.
Thus no $s > 0$ satisfies the condition~\eqref{eq:cond-kuelbs-li-imposs}.

\subsection{Large range of solutions to the Kuelbs-Li condition}
\label{sec:many-solut-kuelbs}

We now present an example of a different configuration, where the Kuelbs-Li condition~\eqref{eq:kuelbs-li} holds for a wide range of values of $s > 0$, and thus does not yield a characterization of the order of magnitude of the small ball probability.

Specifically, we will show that for every $\kappa_0 > 1$ and $\eps \in (0, 1)$, there exists a constant $\kappa_1 = \kappa_1 (\kappa_0, \eps)$ such that the following holds.
For every integer $n \geq \kappa_1$, there exists a nontrivial symmetric convex body $K \subset \R^n$ and $1 < u \leq w (K)/\kappa_0$ such that, for every $s \in \ccint{4 u^2}{n^{1 - \eps} u^2}$,
\begin{equation}
  \label{eq:condition-kuelbs-li-kappa}
  \kappa_1^{-1} s
  \leq \log N \parens*{K, \frac{u}{\sqrt{s}}}
  \leq \kappa_1 s
  \, .
\end{equation}
Thus, the Kuelbs-Li condition~\eqref{eq:condition-kuelbs-li-kappa} admits a range of solutions that is nearly linear (in a multiplicative sense) in the dimension $n$.

Specifically, for $n \geq 2$ we let
\begin{equation*}
  K = B_1^n = \set*{t = (t_i)_{1 \leq i \leq n} \in \R^n : \sum_{i=1}^n \abs{t_i} \leq 1}
\end{equation*}
be the unit $\ell^1$ ball in $\R^n$.
By standard estimates on an expected maximum of independent Gaussian variables,
\begin{equation}
  \label{eq:width-l1}
  \frac{\sqrt{\log n}}{c_0}
  \leq w (K)
  = \Expect*{\max_{1 \leq i \leq n} \abs{X_i}}
  \leq c_0 \sqrt{\log n}
\end{equation}
for some absolute constant $c_0 > 1$.
We thus take $u = \sqrt{\log n}/(c_0 \kappa_0) \leq w (K)/\kappa_0$.
In addition, estimates of Schütt~\cite{schutt1984entropy} on covering numbers of $B_1^n$ show that for $n^{-1/2} \leq r \leq 1/2$,
\begin{equation}
  \label{eq:covering-l1}
  \frac{\log \parens{1 + n r^2}}{c_1 r^2}
  \leq \log N (K, r)
  \leq \frac{c_1 \log \parens{1 + n r^2}}{r^2}
  \, .
\end{equation}
Now if $n^{-(1-\eps)/2} \leq r \leq 1/2$, then $n^\eps \leq 1 + n r^2 \leq 1+n/4 \leq n$ and thus the above yields
\begin{equation*}
  \sqrt{\frac{\eps \log n}{c_1}}
  \leq r \sqrt{\log N (K, r)}
  \leq \sqrt{c_1 \log n}
  \, .
\end{equation*}
Hence, recalling that $u = \sqrt{\log n}/(c_0\kappa_0)$ and letting $s = u^2/r^2$, we deduce that for every
$s \in \ccint{4 u^2}{n^{1-\eps} u^2}$, 
\begin{equation*}
    \frac{\eps c_0^2\kappa_0^2}{c_1}\,s
    \leq \log N\parens*{K,\frac{u}{\sqrt{s}}}
    \leq c_1c_0^2\kappa_0^2\,s.
\end{equation*}

\subsection{Upper and lower deviations for the operator norm}
\label{sec:exampl-oper-norm}

We now justify the claims of Example~\ref{ex:operator-norm}.

Endow the space $\R^{n \times n}$ of $n \times n$ matrices with the Hilbert-Schmidt scalar product defined by $\innerp{A}{B} = \tr \bracks{A^* B}$, for $A, B \in \R^{n \times n}$, with $A^*$ the adjoint of $A$; and denote $\hsnorm{A} = \sqrt{\innerp{A}{A}}$ the associated Hilbert-Schmidt norm.
The unit ball induced by this scalar product is the Hilbert-Schmidt (Schatten-$2$) ball $S_2^n = \set{A \in \R^{n \times n} : \hsnorm{A} \leq 1}$.

Next, let $\norm{A}_1 = \tr \bracks{\parens{A^* A}^{1/2}}$ denote the trace norm (Schatten-$1$ norm) of a matrix $A \in \R^{n \times n}$, and let
% $K = S_1^n = \set*{A \in \R^{n \times n} : \norm{A}_1 \leq 1}$
\begin{equation}
  \label{eq:def-k-schatten-1-ball}
  K = S_1^n = \set*{A \in \R^{n \times n} : \norm{A}_1 \leq 1}
\end{equation}
denote the Schatten-$1$ unit ball.
We then have, for any matrix $G \in \R^{n \times n}$,
\begin{equation}
  \label{eq:op-norm-sup-trace-ball}
  \opnorm{G}
  = \sup_{A \in K} \innerp{A}{G}
  \, .
\end{equation}

Finally, recall that $G_n$ is a random $n \times n$ matrix with \iid standard Gaussian entries.
Thus $\Expect{G_n} = 0$ and for any $A \in \R^{n \times n}$, one has $\Expect{\innerp{A}{G_n}^2} = \hsnorm{A}^2$; that is, $G_n$ is a standard Gaussian vector in the Hilbert-Schmidt space.

The upper bound $\Expect*{\opnorm{G_n}} \leq 2 \sqrt{n}$ follows from Slepian's inequality and can be found, \eg, in~\cite[Theorem~2.13 p.~353]{davidson2001local}.
For the lower bound, we use that $\opnorm{G_n} \geq \hsnorm{G_n} / \sqrt{n}$ followed by the moment comparison from~\cite[Corollary~3]{latala1999gaussian}, which gives
\begin{equation*}
  \Expect*{\opnorm{G_n}}
  \geq \frac{1}{\sqrt{n}} \Expect*{\hsnorm*{G_n}}
  \geq \frac{1}{\sqrt{n}} \sqrt{\frac{2}{\pi}} \Expect*{\hsnorm*{G_n}^2}^{1/2}
  = \sqrt{\frac{2 n}{\pi}}
  \, .
\end{equation*}

The estimates~\eqref{eq:operator-norm-upper-dev} on the upper deviation probability come from the bound $w (K) = \Expect{\opnorm{G_n}} \leq 2 \sqrt{n}$ and the fact that $\rad (K) = 1$ (as $\hsnorm{A} \leq \norm{A}_1$ for $A \in \R^{n \times n}$, with equality if $A = t t^\top$ for $t \in S^{n-1}$), together with Corollary~\ref{cor:upper-tail-gaussian}.

We now turn to the lower deviation estimates~\eqref{eq:operator-norm-lower-dev}.
The inequality $\norm{A}_1 \leq \sqrt{n} \hsnorm{A}$ implies that $n^{-1/2} S_2^n \subset S_1^n = K$, and thus $K \cap r S_2^n = r S_2^n$ for any $r \in \ccint{0}{n^{-1/2}}$.
Thus for any such $r$, by Cauchy-Schwarz and again by~\cite[Corollary~3]{latala1999gaussian},
\begin{align*}
  w (K \cap r S_2^n)
  &= r \Expect*{\hsnorm*{G_n}}
    \leq r \Expect*{\hsnorm*{G_n}^2}^{1/2}
    = r n \, ; \\
  w \parens*{K \cap r S_2^n}
  &\geq \sqrt{\frac{2}{\pi}} r \Expect*{\hsnorm*{G_n}^2}^{1/2}
    = \sqrt{\frac{2}{\pi}} r n
    \, .
\end{align*}
These inequalities imply that for any $u \in \ocint{0}{\sqrt{n}}$, one has
\begin{equation}
  \label{eq:r-star-schatten-1}
  \frac{u}{n}
  \leq r_* \parens*{K, u}
  \leq \sqrt{\frac{\pi}{2}} \frac{u}{n}
  \, .
\end{equation}
In particular, for
% $u \leq \sqrt{n}$ with $u \asymp \sqrt{n}$, one has $r_* (K, u) \asymp 1/\sqrt{n}$.
$\sqrt{n}/c_0 \leq u \leq \sqrt{n}$ with $c_0 \geq 1$, one has $1/\parens*{c_0 \sqrt{n}} \leq r_* (K, u) \leq \sqrt{\pi/(2 n)}$.
The lower deviation estimates~\eqref{eq:operator-norm-lower-dev} then follow from Corollary~\ref{cor:small-ball-constant-regime}.

% {
%   \bibliography{biblio-small-ball}

\begin{thebibliography}{AGHCYN21}

\bibitem[AAGM15]{artstein2015asymptotic}
Shiri Artstein-Avidan, Apostolos Giannopoulos, and Vitali~D. Milman.
\newblock {\em Asymptotic Geometric Analysis, Part I}.
\newblock American Mathematical Society, 2015.

\bibitem[AAGM21]{artstein2021asymptotic}
Shiri Artstein-Avidan, Apostolos Giannopoulos, and Vitali~D. Milman.
\newblock {\em Asymptotic Geometric Analysis, Part {II}}.
\newblock American Mathematical Society, 2021.

\bibitem[AGHCYN21]{alonso2021further}
David Alonso-Guti{\'e}rrez, Mar{\'\i}a~A. Hern{\'a}ndez~Cifre, and Jes{\'u}s
  Yepes~Nicol{\'a}s.
\newblock Further inequalities for the (generalized) {Wills} functional.
\newblock {\em Communications in Contemporary Mathematics}, 23(03):2050011,
  2021.

\bibitem[AMS04]{artstein2004duality}
Shiri Artstein, Vitali Milman, and Stanis{\l}aw~J. Szarek.
\newblock Duality of metric entropy.
\newblock {\em Annals of Mathematics}, 159(3):1313--1328, 2004.

\bibitem[And55]{anderson1955integral}
Theodore~W. Anderson.
\newblock The integral of a symmetric unimodal function over a symmetric convex
  set and some probability inequalities.
\newblock {\em Proceedings of the American Mathematical Society},
  6(2):170--176, 1955.

\bibitem[AT07]{adler2007random}
Robert~J. Adler and Jonathan~E. Taylor.
\newblock {\em Random Fields and Geometry}.
\newblock Springer, 2007.

\bibitem[Bla17]{blaschke1917uber}
Wilhelm Blaschke.
\newblock Uber affine geometrie {VII}: {Neue} extremeingenschaften von ellipse
  und ellipsoid.
\newblock {\em Ber. Verh. S{\"a}chs. Akad. Wiss., Math. Phys. Kl}, 69:412--420,
  1917.

\bibitem[BM87]{bourgain1987new}
Jean Bourgain and Vitaly~D. Milman.
\newblock New volume ratio properties for convex symmetric bodies in
  $\mathbb{R}^n$.
\newblock {\em Inventiones Mathematicae}, 88(2):319--340, 1987.

\bibitem[Bog98]{bogachev1998gaussian}
Vladimir~I. Bogachev.
\newblock {\em Gaussian measures}.
\newblock Number~62 in Mathematical Surveys and Monographs. American
  Mathematical Society, 1998.

\bibitem[Bor75]{borell1975brunn}
Christer Borell.
\newblock The {Brunn-Minkowski} inequality in {Gauss} space.
\newblock {\em Inventiones Mathematicae}, 30(2):207--216, 1975.

\bibitem[BZ88]{burago1988geometric}
Yurii~D. Burago and Viktor~A. Zalgaller.
\newblock {\em Geometric inequalities}.
\newblock Springer-Verlag, Berlin, 1988.

\bibitem[CEFM04]{cordero2004b}
Dario Cordero-Erausquin, Matthieu Fradelizi, and Bernard Maurey.
\newblock The {(B)} conjecture for the {Gaussian} measure of dilates of
  symmetric convex sets and related problems.
\newblock {\em Journal of Functional Analysis}, 214(2):410--427, 2004.

\bibitem[Che76]{chevet1976processus}
Simone Chevet.
\newblock Processus {Gaussiens} et volumes mixtes.
\newblock {\em Zeitschrift f{\"u}r Wahrscheinlichkeitstheorie und verwandte
  Gebiete}, 36(1):47--65, 1976.

\bibitem[DS01]{davidson2001local}
Kenneth~R. Davidson and Stanis{\l}aw~J. Szarek.
\newblock Local operator theory, random matrices and {Banach} spaces.
\newblock In {\em Handbook of the geometry of {Banach} spaces}, volume~1, pages
  317--366. Elsevier, 2001.

\bibitem[Dud67]{dudley1967sizes}
Richard~M. Dudley.
\newblock The sizes of compact subsets of {Hilbert} space and continuity of
  {Gaussian} processes.
\newblock {\em Journal of Functional Analysis}, 1(3):290--330, 1967.

\bibitem[Dvo59]{dvoretzky1959theorem}
Aryeh Dvoretzky.
\newblock A theorem on convex bodies and applications to {Banach} spaces.
\newblock {\em Proceedings of the National Academy of Sciences},
  45(2):223--226, 1959.

\bibitem[Fat03]{fatalov2003constants}
Vadim~R. Fatalov.
\newblock Constants in the asymptotics of small deviation probabilities for
  {Gaussian} processes and fields.
\newblock {\em Russian Mathematical Surveys}, 58(4):725--772, 2003.

\bibitem[Fed59]{federer1959curvature}
Herbert Federer.
\newblock Curvature measures.
\newblock {\em Transactions of the American Mathematical Society},
  93(3):418--491, 1959.

\bibitem[Fer75]{fernique1975regularite}
Xavier Fernique.
\newblock Regularit{\'e} des trajectoires des fonctions al{\'e}atoires
  gaussiennes.
\newblock In {\em Ecole d'Et{\'e} de Probabilit{\'e}s de Saint-Flour
  IV---1974}, pages 1--96. Springer, 1975.

\bibitem[Fer76]{fernique1976corps}
Xavier Fernique.
\newblock Corps convexes et processus gaussiens de petit rang.
\newblock {\em Zeitschrift f{\"u}r Wahrscheinlichkeitstheorie und Verwandte
  Gebiete}, 35(4):349--353, 1976.

\bibitem[GLTV26]{guedon2026random}
Olivier Gu{\'e}don, Alexander~E. Litvak, Kateryna Tatarko, and Beatrice-Helen
  Vritsiou.
\newblock On random diameters of convex bodies.
\newblock {\em Preprint arXiv:2609.24131}, 2026.

\bibitem[GMT05]{giannopoulos2005asymptotic}
Apostolos Giannopoulos, Vitali~D. Milman, and Antonis Tsolomitis.
\newblock Asymptotic formulas for the diameter of sections of symmetric convex
  bodies.
\newblock {\em Journal of Functional Analysis}, 223(1):86--108, 2005.

\bibitem[Gor85]{gordon1985some}
Yehoram Gordon.
\newblock Some inequalities for {Gaussian} processes and applications.
\newblock {\em Israel Journal of Mathematics}, 50:265--289, 1985.

\bibitem[Gor88]{gordon1988escape}
Yehoram Gordon.
\newblock On {Milman's} inequality and random subspaces which escape through a
  mesh in $\mathbb{R}^n$.
\newblock In {\em Geometric Aspects of Functional Analysis}, volume 1317 of
  {\em Lecture Notes in Mathematics}, pages 84--106. Springer Berlin
  Heidelberg, 1988.

\bibitem[Haa81]{haagerup1981best}
Uffe Haagerup.
\newblock The best constants in the {Khintchine} inequality.
\newblock {\em Studia Mathematica}, 70:231--283, 1981.

\bibitem[Had75]{hadwiger1975wills}
Hugo Hadwiger.
\newblock Das {Will'sche} funktional.
\newblock {\em Monatshefte f{\"u}r Mathematik}, 79(3):213--221, 1975.

\bibitem[HJSD79]{hoffmann1979lower}
J{\o}rgen Hoffmann-J{\o}rgensen, Lawrence~A. Shepp, and Richard~M. Dudley.
\newblock On the lower tail of {Gaussian} seminorms.
\newblock {\em The Annals of Probability}, 7(2):319--342, 1979.

\bibitem[KL93]{kuelbs1993metric}
James Kuelbs and Wenbo~V. Li.
\newblock Metric entropy and the small ball problem for {Gaussian} measures.
\newblock {\em Journal of Functional Analysis}, 116(1):133--157, 1993.

\bibitem[KM05]{klartag2005geometry}
Boaz Klartag and Vitali~D. Milman.
\newblock Geometry of log-concave functions and measures.
\newblock {\em Geometriae Dedicata}, 112(1):169--182, 2005.

\bibitem[KR97]{klain1997introduction}
Daniel~A. Klain and Gian-Carlo Rota.
\newblock {\em Introduction to Geometric Probability}.
\newblock Cambridge University Press, 1997.

\bibitem[KV07]{klartag2007small}
Bo'az Klartag and Roman Vershynin.
\newblock Small ball probability and {Dvoretzky's} theorem.
\newblock {\em Israel Journal of Mathematics}, 157(1):193--207, 2007.

\bibitem[Lat02]{latala2002some}
Rafa{\l} Lata{\l}a.
\newblock On some inequalities for {Gaussian} measures.
\newblock {\em Proceedings of the ICM}, II:813--822, 2002.

\bibitem[Led96]{ledoux1996isoperimetry}
Michel Ledoux.
\newblock Isoperimetry and {Gaussian} analysis.
\newblock In {\em Lectures on Probability Theory and Statistics}, volume 1648
  of {\em Lecture Notes in Mathematics}, pages 165--294. Springer, 1996.

\bibitem[Lif]{lifshits2023biblio}
Mikhail~A. Lifshits.
\newblock Bibliography of small deviation probabilities.
\newblock \url{https://airtable.com/shrMG0nNxl9SiGxII/tbl7Xj1mZW2VuYurm}; last
  checked: 2026--09--09.

\bibitem[Lif99]{lifshits1999asymptotic}
Mikhail~A. Lifshits.
\newblock Asymptotic behavior of small ball probabilities.
\newblock In {\em {Probability Theory and Mathematical Statistics. Proceedings
  of the 7th international Vilnius conference}}, pages 453--468, 1999.

\bibitem[Lif12]{lifshits2012lectures}
Mikhail~A. Lifshits.
\newblock {\em {Lectures on Gaussian Processes}}.
\newblock Springer, 2012.

\bibitem[LL99]{li1999approximation}
Wenbo~V. Li and Werner Linde.
\newblock Approximation, metric entropy and small ball estimates for {Gaussian}
  measures.
\newblock {\em The Annals of Probability}, 27(3):1556--1578, 1999.

\bibitem[LM00]{laurent2000adaptive}
Beatrice Laurent and Pascal Massart.
\newblock Adaptive estimation of a quadratic functional by model selection.
\newblock {\em The Annals of Statistics}, 28(5):1302--1338, 2000.

\bibitem[LMN{\etalchar{+}}20]{lotz2020concentration}
Martin Lotz, Michael~B. McCoy, Ivan Nourdin, Giovanni Peccati, and Joel~A.
  Tropp.
\newblock Concentration of the intrinsic volumes of a convex body.
\newblock In B.~Klartag and E.~Milman, editors, {\em Geometric Aspects of
  Functional Analysis}, volume 2266 of {\em Lecture Notes in Mathematics},
  pages 139--167, 2020.

\bibitem[LO99]{latala1999gaussian}
Rafa{\l} Lata{\l}a and Krzysztof Oleszkiewicz.
\newblock Gaussian measures of dilatations of convex symmetric sets.
\newblock {\em The Annals of Probability}, 27(4):1922--1938, 1999.

\bibitem[LO05]{latala2005small}
Rafa{\l} Lata{\l}a and Krzysztof Oleszkiewicz.
\newblock Small ball probability estimates in terms of width.
\newblock {\em Studia Mathematica}, 169(3):305--314, 2005.

\bibitem[LR10]{ledoux2010beta}
Michel Ledoux and Brian Rider.
\newblock Small deviations for {Beta} ensembles.
\newblock {\em Electronic Journal of Probability}, 15:1319 -- 1343, 2010.

\bibitem[LS01]{li2001gaussian}
Wenbo~V. Li and Qi-Man Shao.
\newblock Gaussian processes: {Inequalities}, small ball probabilities and
  applications.
\newblock {\em Handbook of Statistics}, 19:533--597, 2001.

\bibitem[LT91]{ledoux1991probability}
Michel Ledoux and Michel Talagrand.
\newblock {\em Probability in Banach Spaces: Isoperimetry and processes}.
\newblock Springer-Verlag, 1991.

\bibitem[McM91]{mcmullen1991inequalities}
Peter McMullen.
\newblock Inequalities between intrinsic volumes.
\newblock {\em Monatshefte f{\"u}r Mathematik}, 111(1):47--53, 1991.

\bibitem[Mil71]{milman1971new}
Vitali~D. Milman.
\newblock New proof of the theorem of {A. Dvoretzky} on intersections of convex
  bodies.
\newblock {\em Functional Analysis and its Applications}, 5(4):288--295, 1971.

\bibitem[Mil85]{milman1985random}
Vitali~D. Milman.
\newblock Random subspaces of proportional dimension of finite dimensional
  normed spaces: {Approach} through the isoperimetric inequality.
\newblock In {\em Lecture notes in mathematics}, pages 106--115, 1985.

\bibitem[Mou25]{mourtada2025universal}
Jaouad Mourtada.
\newblock Universal coding, intrinsic volumes, and metric complexity.
\newblock {\em Journal of the European Mathematical Society}, 2025.

\bibitem[MS97]{milman1997global}
Vitali~D. Milman and Gideon Schechtman.
\newblock Global versus local asymptotic theories of finite-dimensional normed
  spaces.
\newblock {\em Duke Mathematical Journal}, 90(1):73--93, 1997.

\bibitem[Pis86]{pisier1986probabilistic}
Gilles Pisier.
\newblock Probabilistic methods in the geometry of {Banach} spaces.
\newblock In G.~Letta and M.~Pratelli, editors, {\em Probability and Analysis},
  volume 1206 of {\em Lecture Notes in Mathematics}, pages 167--241. Springer,
  1986.

\bibitem[Pis99]{pisier1999volume}
Gilles Pisier.
\newblock {\em The volume of convex bodies and {Banach} space geometry}.
\newblock Cambridge University Press, 1999.

\bibitem[PP13]{paouris2013small}
Grigoris Paouris and Peter Pivovarov.
\newblock Small-ball probabilities for the volume of random convex sets.
\newblock {\em Discrete \& Computational Geometry}, 49(3):601--646, 2013.

\bibitem[PP26]{paouris2026optimal}
Grigoris Paouris and Reese Pathak.
\newblock Optimal mean width and metric entropy estimates for convex bodies.
\newblock {\em Preprint arXiv:2607.29522; v2, October 2026}, 2026.

\bibitem[PPV19]{paouris2019quantitative}
Grigoris Paouris, Peter Pivovarov, and Petros Valettas.
\newblock On a quantitative reversal of {Alexandrov's} inequality.
\newblock {\em Transactions of the American Mathematical Society},
  371(5):3309--3324, 2019.

\bibitem[PTJ86]{pajor1986subspaces}
Alain Pajor and Nicole Tomczak-Jaegermann.
\newblock Subspaces of small codimension of finite-dimensional {Banach} spaces.
\newblock {\em Proceedings of the American Mathematical Society},
  97(4):637--642, 1986.

\bibitem[PTV22]{paouris2022hypercontractivity}
Grigoris Paouris, Konstantin Tikhomirov, and Petros Valettas.
\newblock Hypercontractivity and lower deviation estimates in normed spaces.
\newblock {\em The Annals of Probability}, 50(2):688--734, 2022.

\bibitem[RS57]{rogers1957difference}
Claude~A. Rogers and Geoffrey~C. Shephard.
\newblock The difference body of a convex body.
\newblock {\em Archiv der Mathematik}, 8(3):220--233, 1957.

\bibitem[San49]{santalo1949invariante}
Luis~A. Santal{\'o}.
\newblock Un invariante afin para los cuerpos convexos del espacio de $n$
  dimensiones.
\newblock {\em Portugaliae Mathematica}, 8(4):155--161, 1949.

\bibitem[Sch84]{schutt1984entropy}
Carsten Sch{\"u}tt.
\newblock Entropy numbers of diagonal operators between symmetric {Banach}
  spaces.
\newblock {\em Journal of Approximation Theory}, 40(2):121--128, 1984.

\bibitem[Sch13]{schneider2013convex}
Rolf Schneider.
\newblock {\em Convex bodies: The Brunn-Minkowski theory}.
\newblock Cambridge University Press, second edition, 2013.

\bibitem[Sle62]{slepian1962one}
David Slepian.
\newblock The one-sided barrier problem for {Gaussian} noise.
\newblock {\em Bell System Technical Journal}, 41(2):463--501, 1962.

\bibitem[ST78]{sudakov1978extremal}
Vladimir~N. Sudakov and Boris~S. Tsirel'son.
\newblock Extremal properties of half-spaces for spherically invariant
  measures.
\newblock {\em Journal of Soviet Mathematics}, 9(1):9--18, 1978.

\bibitem[Sud71]{sudakov1971gaussian}
Vladimir~N. Sudakov.
\newblock Gaussian random processes and measures of solid angles in {Hilbert}
  space.
\newblock {\em Soviet Mathematics. Doklady}, 12:412--415, 1971.

\bibitem[Sud76]{sudakov1976geometric}
Vladimir~N. Sudakov.
\newblock Geometric problems of the theory of infinite-dimensional probability
  distributions.
\newblock {\em Trudy Mat. Inst. Steklov}, 141:3--191, 1976.

\bibitem[SW08]{schneider2008stochastic}
Rolf Schneider and Wolfgang Weil.
\newblock {\em Stochastic and Integral Geometry}.
\newblock Springer Berlin Heidelberg, 2008.

\bibitem[Tal87]{talagrand1987regularity}
Michel Talagrand.
\newblock Regularity of {Gaussian} processes.
\newblock {\em Acta Mathematica}, 159:99--149, 1987.

\bibitem[Tal93]{talagrand1993new}
Michel Talagrand.
\newblock New {Gaussian} estimates for enlarged balls.
\newblock {\em Geometric and Functional Analysis}, 3(5):502--526, 1993.

\bibitem[Tal94]{talagrand1994small}
Michel Talagrand.
\newblock The small ball problem for the {Brownian} sheet.
\newblock {\em The Annals of Probability}, 22(3):1331--1354, 1994.

\bibitem[Tal21]{talagrand2021upper}
Michel Talagrand.
\newblock {\em Upper and lower bounds for stochastic processes: Decomposition
  Theorems}.
\newblock Ergebnisse der Mathematik und ihrer Grenzgebiete. Springer, 2nd
  edition, 2021.

\bibitem[TIS76]{tsirelson1976gaussian}
Boris~S. Tsirel'son, Ildar~A. Ibragimov, and Vladimir~N. Sudakov.
\newblock Norms of {Gaussian} sample functions.
\newblock In {\em Proceedings of the {3rd Japan--USSR Symposium on Probability
  Theory}}, pages 20--41, 1976.

\bibitem[Tsi86]{tsirelson1986geometric}
Boris~S. Tsirel'son.
\newblock A geometric approach to maximum likelihood estimation for
  infinite-dimensional {Gaussian} location. {II}.
\newblock {\em Theory of Probability and its Applications}, 30(4):820--828,
  1986.

\bibitem[Ver18]{vershynin2018high}
Roman Vershynin.
\newblock {\em High-dimensional probability: An introduction with applications
  in data science}, volume~47.
\newblock Cambridge University Press, 2018.

\bibitem[Vit96]{vitale1996wills}
Richard~A. Vitale.
\newblock The {Wills} functional and {Gaussian} processes.
\newblock {\em The Annals of Probability}, 24(4):2172--2178, 1996.

\bibitem[Vit01]{vitale2001intrinsic}
Richard~A. Vitale.
\newblock Intrinsic volumes and {Gaussian} processes.
\newblock {\em Advances in Applied Probability}, 33(2):354--364, 2001.

\bibitem[Wil73]{wills1973gitterpunktanzahl}
J{\"o}rg~M. Wills.
\newblock Zur gitterpunktanzahl konvexer mengen.
\newblock {\em Elemente der Mathematik}, 28:57--63, 1973.

\end{thebibliography}
%   \bibliographystyle{alpha}
% }

% \begingroup
% \small

\newcommand{\etalchar}[1]{$^{#1}$}

% \endgroup

\end{document}